\documentclass[english]{smfart}

\usepackage{smfthm}
\usepackage[dvips]{epsfig}
\usepackage{amsmath}
\usepackage{amssymb}
\usepackage[all]{xy}
\usepackage[germanb,francais,english]{babel}
\usepackage{smfhs}

\theoremstyle{plain}

\newcommand{\M}{{\mathcal M}}
\renewcommand{\O}{{\mathcal O}}
\newcommand{\F}{{\mathcal F}}
\newcommand{\R}{{\mathbb R}}
\newcommand{\Z}{{\mathbb Z}}
\newcommand{\C}{{\mathbb C}}
\newcommand{\Q}{{\mathbb Q}}

\renewcommand{\L}{{\mathcal L}}
\renewcommand{\P}{{\mathbb P}}
\newcommand{\vt}{{X_\Sigma}}
\newcommand{\primcol}{{\mathcal P}}
\newcommand{\primcolset}{{\mathfrak P}}

\newcommand{\proofend}{\hfill $\square$}

\newcommand{\we}{{\omega}}      
\newcommand{\totwe}[1]{\we_{\text{total}}^{#1}}    
\newcommand{\cocl}{{\mbox{$\lambda$}}}   
\newcommand{\gmodsp}[1][m]{\M_{g,#1}^A(X)}  
\newcommand{\modsp}[1][m]{\M_{0,#1}^A(\vt)} 
\newcommand{\ucgmodsp}[1][m]{{\mathcal C}_{g,#1}^A(X)}  
\newcommand{\ucmodsp}[1][m]{{\mathcal C}_{0,#1}^A(\vt)}  
\newcommand{\mygraph}{{\Gamma}}   
\newcommand{\gtype}[1][{}]{{\mygraph_{#1}^{\text{type}}}}   
\newcommand{\gtop}[1][{}]{{\mygraph_{#1}^{\text{top}}}}    
\newcommand{\gdmsp}[1][m]{\overline{\M}_{g,#1}}    
\newcommand{\dmsp}[1][m]{\overline{\M}_{0, #1}}   
\newcommand{\ucgdmsp}[1][m]{\overline{\mathcal C}_{g,#1}}  

\newcommand{\virt}{{\text{vir}}}
\newcommand{\virgg}[1]{{[#1]^\virt}}
\newcommand{\viry}{{\virgg{Y}}}

\newcommand{\lotimes}{\stackrel{L}{\otimes}}
\newcommand{\face}{\diamond}

\newcommand{\mycomment}[1]{}

\newcommand{\mylabel}[1]{\label{#1}\mycomment{#1}}

\newcommand{\rk}{\operatorname{rk}}

\newcommand{\Ext}{\operatorname{Ext}}
\newcommand{\Hom}{\operatorname{Hom}}

\newcommand{\Spec}{\operatorname{Spec}}

\newcommand{\id}{\operatorname{id}}
\newcommand{\Vertex}{\operatorname{Vert}}
\newcommand{\Edge}{\operatorname{Edge}}
\newcommand{\val}{\operatorname{val}}
\newcommand{\V}{\Vertex}
\newcommand{\E}{\Edge}
\newcommand{\Aut}{\operatorname{Aut}}

\newcommand{\ev}{\operatorname{ev}}

\newcommand{\topmap}{\operatorname{top}}
\newcommand{\typemap}{\operatorname{type}}
\newcommand{\Lin}{\operatorname{Lin}}

\newcommand{\fix}{{\text{fix}}}
\newcommand{\move}{{\text{move}}}

\newcommand{\vdim}{{\dim_\virt}}

\newcounter{rastera}\newcounter{rasterb}

\newcounter{tempcounter}

\author{Holger Spielberg}
\address{Max Planck Institute for Mathematics in the Sciences\\Inselstra{\ss}e 22--26\\04109 Leipzig\\Allemagne}
\email{Holger.Spielberg@mis.mpg.de}
\urladdr{http://personal-homepages.mis.mpg.de/spielb}

\title{The Gromov--Witten invariants of symplectic toric manifolds}
\alttitle{Les invariants de Gromov--Witten des vari\'et\'es toriques symplectiques}

\begin{document}

\begin{abstract}
We study the fix point components of the big torus action on the
moduli space of stable maps into a smooth projective toric
variety, and apply Graber and Pandharipande's localization formula
for the virtual fundamental class to obtain an explicit formula

for the Gromov--Witten invariants of toric varieties. Using this
formula we compute all genus--$0$ $3$--point invariants of the
Fano manifold $\P(\O_{\P^2}(2) \oplus 1)$, and we show for 
the (non--Fano) manifold $\P(\O_{\P^2}(3) \oplus 1)$
that its quantum cohomology ring does not
correspond to Batyrev's ring defined in \cite{bat93}.
\end{abstract}

\selectlanguage{francais}\begin{altabstract}
Nous \'etudions les composantes des points fixes de l'action du gros
tore sur l'espace de modules des applications stables dans une
vari\'et\'e torique projective lisse. En appliquant la formule de
localisation de Graber et Pandharipande \`a la classe fondamentale
virtuelle de l'espace de modules, nous obtenons une formule
explicite pour les invariants de Gromov--Witten des vari\'et\'es
toriques. \`A l'aide de cette nouvelle formule, nous calculons tous
les invariants de genre 0 \`a trois points de la vari\'et\'e de
Fano $\P(\O_{\P^2}(2) \oplus 1)$ et montrons que pour la vari\'et\'e
$\P(\O_{\P^2}(3) \oplus 1)$, qui n'est pas de Fano, l'anneau de la
cohomologie quantique ne correspond pas \`a l'anneau d\'efini par
Batyrev dans \cite{bat93}.
\end{altabstract}
\selectlanguage{english}


\thanks{The author gladly acknowledges support by the German Academic
Exchange Office (DAAD) through a graduate scholarship (Hochschulsonderprogramm
II), financed by the German federal state and the German {\it L\"ander}.}

\maketitle

\section{Introduction}

The aim of this article is to give a formula that computes the Gromov--Witten invariants
of symplectic toric manifolds.

\subsection*{Gromov--Witten invariants}
Gromov--Witten invariants and quantum cohomology express essentially the same
symplecto--topological data\footnote{Care has to be taken since there exist several
different versions of quantum cohomology: the big quantum cohomology ring indeed contains
the same data as the genus--0 Gromov--Witten invariants of a symplectic manifold; the small
quantum cohomology ring contains much less data, and in particular not all Gromov--Witten
invariants are needed for its definition. When we refer to the quantum cohomology ring
we usually mean the small version.} first studied by Witten in theoretical
physics (\cite{wit91}). In fact, he looked at quantum cohomology as an example
of a topological $\sigma$--model where what we now call Gromov--Witten invariants
are basically the correlation functions. This lead to the interpretation of these
invariants as counting certain \mbox{(pseudo--)}ho\-lo\-morphic curves in a symplectic manifold.

Let $(M,\omega)$ be a compact symplectic manifold, and $J$ be a compatible almost--complex
structure on $(M,\omega)$. A map $f:(\Sigma_g,j) \longrightarrow (M,J)$ from a genus--$g$ curve
$(\Sigma_g, j)$ to $M$ is called $J$--holomorphic if $f$ is $\C$--linear, namely if
\[ \bar\partial_J f := \frac{1}{2}(df + J\circ df\circ j)= 0. \]
For K\"ahler manifolds $(M,\omega,J)$, these are exactly the
holomorphic maps. Now we fix an integral degree--$2$ homology
class $A\in H_2(M, \Z)$, and only look at $J$--holomorphic maps
such that $f_*[\Sigma_g]=A$. For some classes $A$, there will be
only a finite number of such curves up to reparametrization, and
this number will be, under certain genericity assumptions, one of
the Gromov--Witten invariants of the manifold $(M,\omega)$.

This number, though, is not a priori a symplectic invariant: the construction above
strongly depends on the chosen compatible almost--complex structure $J$. In fact, even
the dimension of the space of $J$--holomorphic maps $f:(\Sigma_g,j)\longrightarrow (M,J)$
with $f_*[\Sigma_g]=A$ might change for different almost--complex structures $J$, that is,
the above number might be defined for some $J$, but not for some others. This phenomenon of a
``moduli space of $J$--holomorphic maps'' being too big comes from the unpleasant property
of the $\bar\partial_J$--operator of  not always being transversal to the zero section
in the infinite dimensional vector bundle
\[ \cE \mapright \operatorname{Map}(\Sigma_g, M) \]
whose fiber at $f\in \operatorname{Map}(\Sigma_g,M)$ is the space
$\cE_f=\Omega^{0,1}(f^*TM)$. In fact, $\bar\partial_J$ is a
Fredholm operator, and its index can be computed using
Riemann--Roch arguments. We will usually refer to this index as
the virtual dimension of the corresponding moduli space, since the
index is equal to the actual dimension of the moduli space when
$\bar\partial_J$ is indeed transversal (to the above mentioned
zero section). Note that being a Fredholm operator in particular
includes the property of the index being finite.

There is, however, another important problem of such a ``definition'' of an invariant:
the moduli space of $J$--holomorphic curves in a degree--$2$ homology class $A$ is in general
not compact. Take for example the family of conics that is given  by the equation
$xy=\varepsilon$. For $\varepsilon>0$, these conics are all smooth, but in the limit $\varepsilon \longrightarrow0$
we obtain a singular conic with a node. In fact, Gromov has proven in \cite{gro85} that this
is all that can happen: a series of $J$--holomorphic maps
converges to a $J$--holomorphic map with singularities at worst nodes, \ie. where the underlying curve $\Sigma_g$ might
have nodes. So to compactify the moduli
space of $J$--holomorphic curves it suffices to add these curves with nodes, an approach that eventually
lead to Kontsevich's space of stable maps. This strategy, though, has one big
disadvantage: the dimensions of the boundary components that we have to add for this compactification can be
bigger than the dimension of the moduli space we started with, even the virtual dimension of
the boundary components might get bigger. So we might end up counting $J$--holomorphic curves
with nodes instead of smooth curves.

In the past years, the above mentioned difficulties have been resolved by different means,
keeping more or less the intuitive idea of the invariant counting certain curves. Ruan and Tian
(\cite{rt95}) were the first who rigorously defined Gromov--Witten invariants in a mathematical
context. They restricted themselves to weakly monotone symplectic manifolds. These manifolds have the nice
property that the virtual dimension of the boundary components is always smaller than the virtual
dimension of the moduli space of smooth curves. Moreover, they were able to show that for a generic
almost--complex structure $J$, the operator $\bar\partial_J$ is transversal for all components
of the compactified moduli space. So, in the case of weakly monotone symplectic manifolds,
the invariant still counts $J$--holomorphic curves. However, the description of all $J$--holomorphic
curves in a symplectic manifold for an arbitrary almost--complex structure $J$ (compatible
with the symplectic structure) remains an unsolved problem.

Later, several successful attempts were undertaken to define Gromov--Witten
invariants for all symplectic manifolds (for example \cite{sie96,lt96,fo96}),
as well as for projective complex varieties (for example \cite{bf97,lt98a}). All constructions
in both categories of varieties follow basically the same principle: instead of
trying to obtain a moduli space of the expected dimension with a fundamental class,
they take any compatible (respectively the given) almost--complex structure $J$ and construct a
virtual fundamental class in the moduli space corresponding to $J$. The
virtual fundamental class so defined is then supposed to behave as the
fundamental class of a generic moduli space (if it existed at all).

Although the constructions in both categories are technically quite different,
the Gromov--Witten invariants obtained  are the same (see \cite{sie98,lt98b}).
Actually, even the main idea for the construction of the virtual fundamental class
is the same in both approaches: they both  use excess intersection
theory to ``slice out'' a cycle in exactly the right dimension, being led by the
observation that the operator $\bar\partial_J$ is not transversal.
In the algebro--geometric construction this is done by using a particular
{\em tangent obstruction theory} $E^\bullet$, that is a two--term complex of locally free
sheaves on the moduli space $\cM$ and a morphism (in the derived category)
\[ \phi: E^\bullet \mapright L_{\cM}^\bullet\]
to the cotangent complex $L^\bullet_\cM$ of the moduli space, such that the rank
$\rk E^\bullet = \rk(E^0-E^{-1})$ of the complex $E^\bullet$ is constant and equal
to the virtual dimension of the moduli space $\cM$. Roughly speaking, one can say that this
obstruction theory $\phi: E^\bullet \longrightarrow L_\cM^\bullet$ encodes how
the virtual moduli cycle has to be cut out off the moduli space $\cM$.

The above mentioned equivalence of the definitions in the two different categories
opens an interesting opportunity for manifolds that are symplectic and
complex varieties at the same time, K\"ahler manifolds: one could try to
use the rather developed machinery of algebraic geometry to finally obtain
symplectic invariants!

\subsection*{Toric manifolds}
Toric manifolds, \ie. those which contain an algebraic torus as an open and dense subset
and whose action on itself extends to the entire manifold, are an important set of examples
to consider here because many are in fact K\"ahler. Moreover, although they include
representatives of many classes of manifolds so far looked at in the
context of Gromov--Witten invariants (complex projective space;
Fano and weakly monotone manifolds), most toric manifolds do not fit into any of these groups.
In spite of this diversity, all toric manifolds are combinatorically classified with
the help of fans that basically describe the intersection pattern of the
divisors of the toric variety.

However, what makes toric manifolds particularly nice to us is the action of the
``big'' torus on them. This action has only finitely many stable submanifolds
which again can be easily derived from the fan description of the toric manifold.
In addition, the action on the toric manifold $X$ naturally induces a torus action
on the moduli spaces of stable maps to $X$, the fixed point components of which
can be described combinatorically  in terms of the zero and one dimensional stable
submanifolds in $X$, hence again by fan data.
This opens to us the possibility to apply equivariant theory to our problem.

\subsection*{Equivariant theory}
In \cite{gp97} Graber and Pandharipande have proven a
localization formula for algebraic stacks $Y$ with a
$\C^*$--action and a $\C^*$--equivariant perfect obstruction
theory that can be $\C^*$--equivariantly embedded into a
non--singular Deligne--Mumford stack. Similarly to the classical
localization formula, they look, on a fixed point component $Y_i$
of the action on the stack $Y$, at a decomposition of the
obstruction theory $E_i^\bullet$ restricted to $Y_i$ into the part
that is fixed by the action, and the moving part:
\[ E_i^\bullet = E_i^{\bullet,\fix} \oplus E_i^{\bullet,\move}.\]
Their main observation is that the fixed part $E_i^{\bullet,\fix}$ is again an obstruction
theory for the fixed point component $Y_i$, and that the role of the normal bundle is
taken by the moving part $E^{\bullet,\move}$, accordingly called {\em virtual normal
bundle}: $N_i^\virt = E_{i,\bullet}^\move$, where $E_{i,\bullet}$ is the dual complex to
$E^\bullet_i$.

To be precise, let $Y$ be an algebraic stack with a $\C^*$--action
that can be $\C^*$--equivariantly embedded into a non-singular
Deligne--Mumford stack. Let $\phi:E^\bullet \rightarrow
L^\bullet_Y$ be a $\C^*$--equivariant perfect obstruction theory
for $Y$, $[Y,E^\bullet]$ and $[Y_i, E_i^\bullet]$ be the virtual
fundamental classes of $Y$ and $E^\bullet$, and of the fixed point
components $Y_i$ and the induced perfect obstruction theories
$E_i^\bullet$, respectively. Then they have shown the following
{\em localization formula} \cite{gp97}:
\[
[Y, E^\bullet] = \iota_* \sum_i \frac{[Y_i,E_i^\bullet]}{e^{\C^*}(N_i^{\virt})}.
\]

In particular, this localization formula holds for the moduli
stacks $\cM^A_{g,m}(\vt)$ of stable maps to a smooth projective
toric variety\footnote{In fact, the theorem holds for all moduli
stacks of stable maps into a non--singular variety.}. Furthermore,
let $G$ be a $\C^*$--equivariant bundle with rank $\rk G=\deg
[Y,E^\bullet]$. Denote by $G_i$ its restriction to the fixed point
components $Y_i$ of $Y$. Then the localization formula immediately
implies the following ``Bott residue formula'' \cite{gp97} which
we will use for the computation of the (algebraic) genus--zero
Gromov--Witten invariants of a smooth projective toric variety
$\vt$:
\begin{equation}\label{eq:bottformula}
\int_{[Y,E^\bullet]} e(G) = \sum_i \int_{[Y_i, E_i^\bullet]}
\frac{e^{\C^*}(G_i)}{e^{\C^*}(N_i^{\virt})},
\end{equation}
an equation that holds in the localized ring
$A^{\C^*}(Y)\otimes\Q[\mu,\frac{1}{\mu}]$.

\noindent Note that since $\rk G=\deg [Y,E^\bullet]$ we actually have
\[ \int_{[Y,E^\bullet]} e(G) = \int_{[Y,E^\bullet]} e^{\C^*}(G).\]
In particular, the right hand side of (\ref{eq:bottformula}) takes values in $\Q$,
and not just in a polynomial ring over $\Q$.

\subsection*{Gromov--Witten invariants of symplectic toric manifolds}
The Bott resi\-due formula is indeed very helpful for resolving our
initial problem of calculating the Gromov--Witten invariants of
symplectic toric manifolds. Remember that the original idea of
Gromov--Witten invariants was that they count certain
holomorphic\footnote{Or, in the general set--up,
pseudo--holomorphic.} curves. In a generalized version and in the
set--up of virtual fundamental classes, these invariants are
defined by integration over the virtual fundamental class:
\mycomment{eq:gwiformel}
\begin{equation}\label{eq:gwiformel}
\Psi^A_{g,m}(\beta;\alpha_1,\ldots,\alpha_m):=\int_{[\gmodsp,E^\bullet]}
\ev^*(\alpha_1\otimes\ldots\otimes\alpha_m)\wedge\pi^*\beta,
\end{equation}
where $\alpha_1,\ldots,\alpha_m\in H^*(X;\Z)$, $\beta\in
H^*(\gdmsp)$, $\ev:\gmodsp\rightarrow X^m$ is the $m$--point
evaluation map, and $\pi:\gmodsp\rightarrow \gdmsp$ is the natural
forgetting (and stabilization) morphism to the Deligne--Mumford
space of stable curves.

Now let $X=\vt$ be a $d$--dimensional smooth projective toric
variety. Then the cohomology of $\vt$ is generated by its
$(\C^*)^d$--invariant divisors. Therefore the classes $\alpha_i\in
H^*(\vt, \Z)$ can be expressed as the Euler classes of some
$(\C^*)^d$--equivariant bundles on $\vt$, and since the action on
the moduli space $\gmodsp$ is the pull back action, the same
applies to the class
$\ev^*(\alpha_1\otimes\ldots\otimes\alpha_m)$. If we restrict to
the case where the class $\beta\in H^*(\gdmsp)$ is
trivial\footnote{Note that this is no restriction to the class
$\beta$ if we only look at genus--zero three--point Gromov--Witten
invariants, \ie. when $g=0$ and $m=3$, since the moduli space
$\dmsp[3]$ consists of just a single point.}, \ie.
$\beta=1=P.D.([\gdmsp])$, we can apply\footnote{Although they
proved their Localization Theorem only for $(\C^*)$--actions, it
obviously generalizes to (diagonal) torus actions: we just
``decompose'' the $(\C^*)^d$--action into $d$ commutative
$(\C^*)$--actions, and apply their localization formula $d$
times.} Graber and Pandharipande's Bott residue formula
(\ref{eq:bottformula}) to compute the above integral
(\ref{eq:gwiformel}).

Hence to eventually obtain the values of these Gromov--Witten invariants, we have to study the
objects on the right hand side of equation (\ref{eq:bottformula}), \ie. the fixed point components
in $\gmodsp$, their virtual fundamental class and their virtual normal bundle, and the restrictions
to the fixed point components of the equivariant bundles corresponding to the classes $\alpha_i$.
In the rest of this section we will restrict ourselves to genus--zero maps, \ie. the moduli spaces
$\modsp$.

\noindent{\bf Fixed point components in $\modsp$:} To describe the
fixed point components in the moduli space of stable maps
$\modsp$, we have generalized Kontsevich's graph approach
\cite{kon95} that he uses in the case of $X=\C\P^n$. The main
observation is that the irreducible components of a stable map
$(C;\underline{x};f)$ that is fixed by the $(\C^*)^d$--action have
to be mapped either to a fixed point of the action in $\vt$ or to
an irreducible one--dimensional $(\C^*)^d$--invariant subvariety
of $\vt$. Moreover, the irreducible components of $C$ that are not
mapped to a point are rigid in each fixed point component. Hence
the fixed point components are essentially products of
Deligne--Mumford spaces of stable curves, a fact that makes it
particularly easy to compute their virtual fundamental class: for
the Deligne--Mumford spaces of stable curves $\dmsp$, it is just
the usual fundamental class, $[\dmsp]^\virt=[\dmsp]$.

\noindent{\bf The virtual normal bundle:} For the study of the virtual normal bundle, or the moving
part of the obstruction theory $E^\bullet$, we consider a $(\C^*)^d$--equivariant long exact sequence
derived from a the pull back to the fixed point components of a distinguished triangle containing $E^\bullet$
(see Section \ref{sec:virtnormbdl}). This way we can reduce the computation of the equivariant Euler class
of the virtual normal  bundle to the computation of the equivariant Euler classes of bundles such as
$R^i\pi_*\underline\Hom(f^*\Omega^1_{\vt}, \O_\cC)$ and $R^i\pi_*\underline\Hom(\Omega^1_{\cC/\cM}(D), \O_\cC)$,
where $\pi: \cC \longrightarrow \cM$ is a $(\C^*)^d$--fixed stable map to $\vt$, and $f:\cC \longrightarrow \vt$
is the universal map to $\vt$.

The main result of this thesis is Theorem \ref{thm:maintheorem} giving an explicit formula
for the genus--zero Gromov--Witten invariants
\mycomment{eq:gwiform}
\begin{equation}\label{eq:gwiform}
\Psi^A_{0,m}(1; \alpha_1, \ldots, \alpha_m)
\end{equation}
of a smooth projective toric variety. This
formula gives in particular all genus--zero three--point
Gromov--Witten invariants of a smooth projective toric variety.

Gromov--Witten invariants and the quantum cohomology of toric
varieties have already been studied by various authors. First
claims on the structure of the quantum cohomology ring were made
by Batyrev in \cite{bat93}, though without the rigorous framework
of the subject that is now available. Givental has computed the
quantum cohomology of weakly monotone toric varieties using
``mirror techniques'' and equivariant methods
(\cite{giv96,giv97}). By using the generalized
Vafa--Intriligator formula, certain Gromov--Witten invariants can
be obtained using a presentation of the quantum cohomology ring
coming from a presentation of the ordinary cohomology ring
(\cite{sie97}). Recently, Qin and Ruan (\cite{qr98}) have
studied the quantum cohomology ring and some of the Gromov--Witten
invariants of certain projective bundles over $\C\P^n$. In
particular they verify Batyrev's conjecture for a small class of
such bundles (Theorem 5.21); our example $\P_{\C\P^2}(\O(2)\oplus
1)$, however, is not treated by their theorem. Moreover, we can show that
the quantum cohomology ring of $\P_{\C\P^2}(\O(3)\oplus 1)$
does not coincide with Batyrev's ring (\cite{bat93}).
Lian, Liu and Yau
\cite{lly97} have also studied the quantum cohomology ring of
complex projective space in an equivariant setting, however so
far they have not yet generalized their results to a bigger class of
manifolds.

\subsection*{Contents}
The paper is structured as follows. In Section \ref{sec:stablemc}
we will recall the definition of the moduli spaces of stable curves
and maps, and give some of their properties. In Section \ref{sec:defgwi} 
we will describe
the construction of the virtual fundamental class in the sense of
Behrend and Fantechi (\cite{bf97,beh97}), and will describe the
obstruction theory used for the Gromov--Witten invariants. 
Graber and Pandharipande's localization formula
will be discussed in Section \ref{sec:gploc}. In Section
\ref{sec:torvar} we will recall the definition and some properties
of toric manifolds. Torus actions on toric
varieties and their moduli spaces of stable maps will be discussed
in Section \ref{sec:tvaction}. In
Section \ref{sec:virtnormbdl} we will determine for an arbitrary
projective toric manifold the virtual normal bundle to the fixed
point components of the moduli space of stable maps to $\vt$ for
the induced $(\C^*)^d$--action. This leads to an explicit formula
for all genus--0 Gromov--Witten invariants of the form
(\ref{eq:gwiform}) for any smooth projective toric variety. In Section
\ref{sec:sumgraph} we will prove some useful lemmata on the 
combinatorics in our formula, thus improving it slightly for practical
computations. As an
application and example, we show how to derive the Gromov--Witten
invariants and the quantum cohomology of projective space $\P^n$
and  the Fano threefold $\P(\O_{\P^2}(2)\oplus1)$ in Section
\ref{sec:examples}. In this Section we also prove the Proposition
on the quantum cohomology ring of $\P_{\C\P^2}(\O(3)\oplus 1)$. 

A big part of this article comes from the author's thesis \cite{spi99}.
The main theorem (Theorem \ref{thm:maintheorem}) as well as its application
to the quantum cohomology ring of $\P_{\C\P^2}(\O(3)\oplus 1)$ was announced
in \cite{spi99b}.

The author wants to thank Mich\`ele Audin, Olivier Debarre, Emmanuel Peyre,
Claude Sabbah, Bernd Siebert and Tilmann Wurzbacher for discussions and
many useful remarks, as well as IRMA Strasbourg and the MPI for Mathematics
in the Sciences, Leipzig, for their hospitality while preparing this article.
\subsection*{General conventions}
In the algebro--geometric category, we always work over the field
of complex numbers $k=\C$, unless otherwise mentioned. Accordingly, dimensions of varieties
are given as complex dimensions.

Although we mostly work in the algebro--geometric category, we prefer to use homology and
cohomology instead of Chow groups.

\newpage

\section{Stable curves and maps, and their moduli spaces}\mylabel{sec:stablemc}

Prestable and stable curves have been intensively studied since Deligne and
Mumford's first paper \cite{dm69} on the moduli space $\cM_g$ for $g\ge2$ (and no
marked points). Later, their results have been extended by Knudsen 
(\cite{knu76a,knu83b,knu83c})
to marked stable curves. In \cite{kee92}, Keel has given a different description of the
genus--0 moduli spaces $M_{0,m}$ as subsequent blow ups.

The notion of a stable map to a smooth variety $X$ is a generalization of
stable curves that is due to Kontsevich. In fact it turned out that the space
of stable maps is the ``right'' compactification of the space of (J--)holomorphic
maps in view of Gromov compactness.

The recently published book \cite{hm98} by Harris and Morrison
collects many of the results known about stable curves and maps, and their 
moduli spaces, and gives many references to the literature.

For the reader's convenience, we will repeat their definition and some 
of their properties that we will use later on in the paper.

\subsection{Prestable and stable curves}
Let $S$ be a scheme, and $g, m\ge0$ be some non--negative integers.
\begin{defi}\mylabel{def:prestablecurve}
A {\em genus--$g$ prestable curve with $m$ marked points} is a flat
and proper morphism $\pi: C \rightarrow S$ together with $m$ distinct sections
$x_1,\ldots, x_m: S \rightarrow C$ such that:

\begin{enumerate}
\item the geometric fibers $C_s=\pi^{-1}(s)$ of $\pi$ are reduced and connected
curves with at most ordinary double points;
\item $C_s$ is smooth at $P_i:=x_i(s)$ $(1\leq i\leq m)$;
\item $P_i\neq P_j$ for $i\neq j$.
\item the algebraic genus of the fibers is $g$: $\dim H^1(C_s, \O_{C_s})=g$.
\setcounter{tempcounter}{\value{enumi}}
\end{enumerate}

\noindent Such a prestable curve is called {\em stable} if it
fulfills in addition the following {\em stability condition}:

\begin{enumerate}\setcounter{enumi}{\value{tempcounter}}
\item\mylabel{def:stabcond} The number of points where a non--singular rational component $E$
of $C_s$ meets the rest of $C_s$ plus the number of marked points $P_i$
on $E$ is at least three.
\end{enumerate}
\end{defi}

\begin{defi}
Let ${\mathfrak M}_{g,m}$ and $\gdmsp[m]$ be the categories
of $m$--pointed pre\-stable respectively stable curves. {\em
Morphisms} in these categories are diagrams of the
form
\[\xymatrix{
C' \ar[r]^\phi \ar[d]^{\pi'} & C \ar[d]^{\pi}\\
S' \ar@(ul,ur)[u]^-{x_i'} \ar[r]^\psi & S \ar@(ul,ur)[u]^-{x_i}
}\]
where
\begin{enumerate}
\item $\phi\circ x_i'=x_i\circ \psi$ for $1\leq i\leq m$,
\item $\phi$ and $\pi'$ induce an isomorphism $C' \stackrel{\thicksim}{\rightarrow} C\times_S S'$.
\end{enumerate}
If the morphism of schemes $\psi: S' \longrightarrow S$ is an isomorphism, we call the morphism
between the two curves an
{\em isomorphism}.
\end{defi}

\begin{theo}[\protect{\cite[Theorem 2.7]{knu83b}}] For all relevant $g$ and $m$,
$\gdmsp[m]$
is a separated algebraic stack, proper and smooth over $\Spec(\Z)$ of dimension
$\dim \gdmsp[m]=3(g-1)+m$.
\end{theo}

\begin{rema} We will not give the definition of a stack and refer the reader
for example to \cite{vis89} or \cite{lm92}.
\end{rema}

\begin{rema}
In the genus--$0$ case, $\dmsp[m]$ is in fact a fine moduli space and a non--singular
variety. Although our applications later on will only involve genus--$0$ curves
and maps we have nonetheless chosen to introduce $\dmsp[m]$ as stacks, since
the corresponding moduli problem for stable maps will no longer admit a fine
moduli space (even for genus--$0$ maps).
\end{rema}

\subsection{The universal curve of the moduli stack of stable curves}

The moduli stack of stable maps $\gdmsp$ admits a universal curve
$\ucgdmsp\MRE{\pi} \gdmsp$, that is for a stable curve
$ C \longrightarrow S$ and its map $S \longrightarrow \gdmsp$ to the moduli stack
there is a map $C \longrightarrow \ucgdmsp$ such that the following diagram
is commutative:
\[
\xymatrix{
C \ar[r] \ar[d] &
{\ucgdmsp} \ar[d]_{\pi}\\
S \ar[r] &
{\gdmsp}.
}
\]
Moreover, the description of the universal curve stack is
particularly easy: it is just the moduli stack of stable curves
with one extra marked point: $\ucgdmsp = \gdmsp[{m+1}]$. The map
$\pi: \ucgdmsp \longrightarrow \gdmsp$ is the forgetting morphism, \ie. the 
natural morphism that forgets the extra marked point and stabilizes:
\begin{alignat*}{2}
\gdmsp[m+1]&\mapright&& \gdmsp[m]\\
(C; p_1, \ldots, p_m, p_{m+1})&\DMRI{}&& (\tilde{C}; p_1, \ldots, p_m),
\end{alignat*}
where $\tilde{C}$ is the curve resulting from $C$ after
stabilization (if necessary).

\subsection{The universal cotangent lines on $\dmsp$}\mylabel{sec:lbdms}

Consider the universal curve ${\mathcal C}_{0,m} \longrightarrow \dmsp[m]$
and the $m$ sections $x_1, \ldots, x_m$ given by the marked
points. Let $K_{{\mathcal C}/\M}$ be the cotangent bundle to the
fibers of ${\mathcal C}_{0,m} \longrightarrow \dmsp[m]$. Then the {\em $i^{th}$
universal cotangent line} is defined to be $\L_i:= x_i^*(K_{{\mathcal
C}/\M})$. In other words, over a stable curve $(C; x_1, \ldots,
x_m)\in \dmsp[m]$ the fiber of the universal cotangent line bundle
${\mathcal L}_i$ is just the cotangent space $T_{x_i}^* C$ of $C$ at
the point $x_i$.

For a tuple $(d_1, \ldots, d_m)$ of non--negative integers satisfying the condition
$ \sum_i d_i = \dim  \dmsp[m] = m - 3$, 
define the number (\cf. \cite{wit91})
\mycomment{eq:dmintersection}
\begin{equation}\label{eq:dmintersection}
\langle \tau_{d_1}\tau_{d_2}\cdots\tau_{d_m}\rangle :=
\int_{\dmsp[m]} c_1(\L_1)^{d_1}\wedge\ldots\wedge c_1(\L_m)^{d_m}.
\end{equation}
If the $d_i$ do not satisfy the dimension equation
$\sum_i d_i = m - 3$, or if one of the $d_i<0$, we set $\langle
\tau_{d_1}\tau_{d_2}\cdots\tau_{d_m}\rangle :=0$.

\begin{rema}
Note that these integrals are obviously symmetric in the tuple
$(d_1, \ldots d_m)$. Therefore we can abbreviate $\langle \tau_{d_1}
\cdots \tau_{d_m}\rangle$ by using exponents, that is for example
$\langle \tau_1\tau_1\tau_0\tau_0\tau_0\rangle$ simply becomes
$\langle\tau_1^2\tau_0^3\rangle$, as does 
$\langle\tau_1\tau_0\tau_1\tau_0\tau_0\rangle$. Remark
that the sum of the exponents still gives the number of marked
points, that is the Deligne--Mumford space of stable curves we
are working on.
\end{rema}

It was conjectured by Witten \cite{wit91} and later proven by
Kontsevich \cite{kon92} that these intersection numbers fulfill
the so--called {\em string equation}:
\[
\langle\tau_0\prod_{i=1}^m \tau_{k_i}\rangle =
\sum_{j=1}^m\langle \tau_{k_j-1}\prod_{i\neq j}\tau_{k_i}\rangle.
\]
With the obvious ``initial condition'' $\langle\tau_0^3\rangle=1$ we can thus obtain
the following explicit formula for these products, the proof of which is 
a straightforward computation.

\begin{coro}[\protect{\cite[Exercise 2.63]{hm98}}]\mylabel{cor:dmexplicit}
The intersection numbers (\ref{eq:dmintersection}) on the Deligne--Mum\-ford space
of stable curves are given by:
\[ \langle \prod_{i=1}^m \tau_{k_i}\rangle = \frac{(m-3)!}{\prod_{i=1}^m k_i!}.
\] \proofend
\end{coro}

\subsection{The moduli space of stable maps}\mylabel{subsec:defmodspm}

Stable maps are a generalization of stable curves that one can retrieve in the
following definition simply by taking the manifold $X$ to be a point:

\begin{defi}\mylabel{def:stablemap}
Let $m,g\ge0$ and $X$ be a smooth variety. A {\em genus--$g$
stable map to $X$ with $m$ marked points}
is given by a genus--$g$ prestable curve $\pi: C \longrightarrow
S$ with marked point sections $x_i: S \longrightarrow C$, and a
morphism $f: C \longrightarrow X$ such that for each geometric fiber $C_s$, 
the non--singular components $E$ of $C_s$ that are mapped to a point by $f$ satisfy
the stability condition \ref{def:stabcond} of definition
\ref{def:prestablecurve}.

A {\em morphism of stable maps} $(\pi: C\rightarrow S; \underline{x}; f)$ and
$(\pi':C'\rightarrow S'; \underline{x}'; f')$ is a morphism
of the two prestable curves commuting with the morphism $f$ and $f'$: $f=f'\circ\phi$:
\[\xymatrix{
C' \ar[r]^\phi \ar[d]^{\pi'} \ar@(ur,ul)[rr]^{f'} & C \ar[d]^{\pi} \ar[r]^f&X\\
S' \ar@(ul,ur)[u]^-{x_i'} \ar[r]^\psi & S. \ar@(ul,ur)[u]^-{x_i}
}\]
Such a morphism is an {\em isomorphism} if the underlying morphism of prestable maps is one.
\end{defi}

\begin{defi}
Let $A \in H_2(X,\Z)$ be an integral degree--2 homology class  of
$X$. We denote by $\gmodsp[m]$ the category of genus--$g$ stable
maps to $X$ with $m$ marked points, such that the push forward by
$f$ of the fundamental class $[C_s]$ of the fibers is
$f_*[C_s]=A$. The morphisms in this category are the morphisms
between stable maps.
\end{defi}

The dimension of the moduli stack is a priori not known. However, by Riemann-Roch
arguments, one finds that the {\em virtual dimension}
of the moduli stack of stable
maps is given by
\[ \vdim \gmodsp[m] = (1-g)(\dim X-3) + \langle c_1(X), A\rangle +m.\]
Unfortunately, even if the moduli stack is not empty altogether, the virtual dimension
and the actual dimension of the moduli stack almost never coincide.

\begin{exem}
A rather classical example for when the virtual dimension of the
moduli space does not coincide with the actual dimension is the
following (see \eg. \cite{aud96}). Let $X=\widetilde{\P^2}$ be
the two dimensional complex projective space blown up at one
point, and let $A=2E$ be twice the class of the exceptional
divisor. The virtual dimension of $\M_{0,0}^{2E}(\widetilde{\P^2})$ 
is equal to $1$. However, since maps in the
class $2E$ have to lie in the exceptional fiber, this moduli stack
is equal to $\M_{0,0}^{2H}(\P^1)$, 
where $H$ is the fundamental class of $\P^1$. The virtual dimension of the latter
moduli stack is two, which is in fact equal to the factual dimension since $\P^1$
is a convex variety (see example \ref{ex:convexvar}).
\end{exem}

\begin{exem}\mylabel{ex:convexvar}
Convex varieties are among the few exceptions where the Rie\-mann--Roch formula actually
gives the accurate dimension of the moduli stack of genus zero stable maps.
A smooth projective variety $X$ is
called {\em convex} if for every morphism $f: \P^1 \rightarrow X$,
\[ H^1(\P^1, f^*TX) = 0.\]
Examples of convex spaces include all homogeneous spaces $G/P$ where $G$ is
a semi--simple Lie group and $P$ is a parabolic subgroup. Hence, projective
spaces, Grassmannians, smooth quadrics, flag varieties, and products of
such spaces are all convex. The beautiful paper of Fulton and Pandharipande
\cite{fp97} gives a very detailed account of genus zero stable maps
to convex manifolds.
\end{exem}

The following well--known lemma provides us with an equivalent criterion
for stability that we will use later on.

\begin{lemm}\mylabel{lem:stableinj}
Let $C$ be a marked rational curve with singularities (over $S=\Spec\C$) that are at worst
double points, and let $D$ be the divisor given by the marked points.
Further, let $X$ be a smooth variety and $f: C \longrightarrow X$ be a map.

Then the map $f$ is stable (with respect to the given marked points)
if and only if the following map induced by the natural map
$f^*\Omega^1_X \stackrel{\phi}{\longrightarrow} \Omega^1_C \longrightarrow \Omega^1_C(D)$ is injective:
\[ \Hom(\Omega^1_C(D), \O_C) \MRE{\Phi} \Hom(f^*\Omega^1_X, \O_C). \]
\end{lemm}

\begin{rema}
The above lemma generalizes directly to any pre--stable curve
$\pi: C \longrightarrow S$ with marked point sections
$x_i:S\longrightarrow C$ and a morphism $f:C \longrightarrow X$:
the tuple $(C\rightarrow S; \underline{x}; f)$ is a stable map if
and only if the morphism
\[ R^0\pi_*\underline\Hom (\Omega^1_{C/S}(D),\O_C) \longrightarrow R^0\pi_*\underline\Hom
(f^*\Omega^1_X, \O_C)\]
is injective. This follows directly from the fact that a morphism of sheaves is injective
if and only if it is injective on each stalk, and from the property that
\begin{align*}
\phantom{\text{and}\quad}R^0\pi_*\underline\Hom (\Omega^1_{C/S}(D),\O_C)_s&= \Hom(\Omega^1_{C_s}(D_s), \O_{C_s}) \quad \text{and}\\
R^0\pi_*\underline\Hom(f^*\Omega^1_X, \O_C)_s&= \Hom(f_s^*\Omega^1_X, \O_{C_s}).
\end{align*}
The latter is implied by Grauert's continuity theorem (see for example 
\cite[Th\'eo\-r\`eme 4.12(ii)]{bs77}).
\end{rema}

\section{Gromov--Witten invariants}\mylabel{sec:defgwi}

Gromov--Witten invariants of a symplectic manifold $X$ are defined
using intersection theory on the moduli space of stable (holomorphic or
pseudo--holomorphic) maps to $X$. They are invariants of the deformation
class of the symplectic structure $\omega$ of $X$, so in particular they ought to be
independent of the (pseudo--)complex structure $J$ compatible with $\omega$.

Unfortunately, even the dimension of the moduli spaces of stable maps 
can vary with the (pseudo--)complex structure. However, these moduli spaces
are the pre--image of zero under the $\bar\partial_J$ operator 
and we would have $\dim \gdmsp = \dim_\virt \gdmsp$
if $\bar\partial_J$ were transversal to the zero section of
$\Omega_C^{0,1}(f^*TX)$ at each stable $J$--holomorphic map $f:C\rightarrow X$ in $\gmodsp$.

Two different approaches have been developed to solve this problem:
one is to try to make $\bar\partial_J$ transversal to the zero section, the other
is to use principles of excess intersection theory to obtain a cycle in
$H_*(\gmodsp)$ of degree equal to the virtual dimension $\dim_{\virt}\gmodsp$ of the moduli space.
The former has been pursued  by Ruan and Tian (\cite{rt95}) for weakly monotone
symplectic manifolds.

The latter has been developed by Behrend and Fantechi as well as Li and Tian
(\cite{bf97,beh97,lt98a}) for all smooth projective complex varieties, and by Fukaya
and Ono, Li and Tian, Ruan, and Siebert (\cite{fo96,lt96,rua96,sie96}) for all smooth
symplectic manifolds\footnote{Of course, this class of manifolds includes the smooth projective
complex varieties, though the constructions by Behrend and Fantechi, and Lian and Tian
are entirely in the algebro--geometric category. In particular, Behrend and Fantechi construct
a cycle in the Chow ring $A_*(\gmodsp)$ of the moduli space.}. The  basic idea of the
construction is as follows: Consider a smooth variety $W$, two smooth subvarieties
$X,Y$ of $W$, and their intersection $Z$:
\[ 
\xymatrix{
Z \ar[r] \ar[d] &
X \ar[d]^f \\
Y \ar[r]^g &
W.
}
\]

Now, if $X$ and $Y$ intersect properly, \ie. if $\dim Z=
\dim X + \dim Y - \dim W$ then the fundamental cycle of Z is the
intersection of the fundamental cycles of $X$ and $Y$: $[Z]=[X]\cdot[Y]$.
Otherwise, using excess intersection theory we can find a cycle in the Chow ring $A_*(Z)$
representing $[X]\cdot[Y]$, the {\em virtual cycle} of $Z$: $[Z]^{\virt}$.
Let $s:Z\rightarrow C_{Y/W}\times_Y Z$ be the zero section of the normal cone to
$Y$ in $W$ pulled back to $Z$. Then $[Z]^{\virt}$ is the intersection
of the zero section $s$ with the normal cone $C_{Z/X}$ to $X$ in $Z$:
\[ [Z]^{\virt}= s^*(C_{Z/X}),\]
where $s^*:A_*(C_{Y/W}\times_Y Z)\rightarrow A_*Z$ is the Gysin morphism
induced by $s$.

Unfortunately, for our moduli problem such an ambient space $W$ and maps
$X,Y\rightarrow W$ do not exist naturally such that $X\times_W Y$ is the
moduli space and $[X]\cdot[Y]$ a virtual moduli cycle with the properties
we want. Instead, the construction will use an {\em obstruction theory} for
$\gmodsp$, a two--term complex $E^\bullet$ on $\gmodsp$ with $\rk E^\bullet=
\dim_{\virt}\gmodsp$.

We will sketch the definition in some generality following 
\cite{bf97,beh97}, and then apply it
to the moduli space of stable maps and Gromov--Witten invariants. 

\subsection{Perfect obstruction theory and virtual fundamental class}

Let $Y$ be a Deligne--Mumford stack, that is an algebraic stack with unramified
diagonal. Let $L^\bullet_Y$ be the cotangent complex of 
$Y$ (see for example \cite{buc81,ill71} for its
definition and properties on schemes, and \cite{lm92} for its
generalization to algebraic stacks). The {\em intrinsic
normal sheaf ${\mathfrak N}_{Y}$} is defined to be the quotient stack
\[ {\mathfrak N}_{Y} := h^1/h^0((L^\bullet_{Y})^\vee)=
[\ker (L_1\rightarrow L_2)/L_0].\]
The {\em intrinsic normal cone ${\mathfrak C}_Y$ of $Y$} is the unique
closed subcone stack ${\mathfrak C}_{Y} \hookrightarrow
{\mathfrak N}_Y$ such that for a local embedding 
\[ 
\xymatrix{
U \ar[r]^f \ar[d]_i&
M\\
Y
}
\]
of $Y$, we have ${\mathfrak C}_{Y}|_U=[C_{U/M}/f^*TM]$ 
(\cite[Corollary 3.9]{bf97}).  The intrinsic normal cone
${\mathfrak C}_Y$ is of pure dimension zero 
(\cite[Theorem 3.11]{bf97}).

\begin{defi}
Let  $Y$ be a Deligne--Mumford stack, that is,
an algebraic stack with unramified diagonal.

Let $E^\bullet=[E^{-1}\rightarrow E^0]$ be a two--term complex of
vector bundles on $Y$. Then a morphism in the derived category from
$E^\bullet$ to the cotangent complex $L^\bullet_Y$
\[ \phi: E^\bullet \rightarrow L^\bullet_Y\]
is called a {\em perfect obstruction theory} for $Y$
if $h^0(\phi)$ is an isomorphism and $h^{-1}(\phi)$ is surjective.
\end{defi}

\begin{rema}
The definition of a perfect obstruction theory in
\cite{bf97} is more general than the one given here, that is they consider
two--term complexes of locally free sheaves ${\mathcal E}^\bullet$.
A two--term complex of vector bundles $E^\bullet$ as above that is
isomorphic to ${\mathcal E}^\bullet$ in the derived category is then called
a global resolution.
\end{rema}

The morphism $\phi$ induces a closed immersion
$\phi^\vee:{\mathfrak N}_Y \rightarrow h^1/h^0((E^\bullet)^\vee)$
(Proposition 2.6 in \cite{bf97}), so $E_1={E^{-1}}^\vee$ is a
global presentation of the quotient stack $h^1/h^0((L^\bullet)^\vee)$ and ${\mathfrak
C}_Y\rightarrow {\mathfrak N}_Y$ embeds into $E_1$. Consider the
fibered product
\[
\xymatrix{
C(E^\bullet) \ar[r] \ar[d] &
E_1 \ar[d]\\
{\mathfrak C}_Y \ar[r] &
[E_1/E_0].
}
\]
Hence, $C(E^\bullet)$ is a closed subcone of the vector bundle $E_1$.
Locally, for a local embedding
$U\rightarrow M$ as above, $i^*C(E^\bullet)$ is those just given by
\[  i^*C(E^\bullet) = \left(C_{U/M}\times_U (i^*E_0)\right)/f^*TM.\]

By this construction, $C(E^\bullet)\rightarrow {\mathfrak C}_Y$ is
smooth of relative dimension $\rk E_0$. Since the intrinsic normal
cone ${\mathfrak C}_Y$ is of pure dimension zero, $C(E^\bullet)$
is thus of pure dimension $\rk E_0$.

\begin{defi}
Let $Y$, $C(E^\bullet)$ and $E_1$ be as above. Let $n=\rk E^\bullet=\rk E^0-
\rk E^{-1}=\rk E_0-\rk E_1$ be the {\em virtual dimension} of $Y$ with respect to
the obstruction theory $E^\bullet$.
The {\em virtual fundamental class}
$[Y,E^\bullet]\in H_n(Y,\Q)$ of $Y$ is the intersection of $C(E^\bullet)$ with the
zero section of $E_1$.
\end{defi}

\begin{rema}
The virtual fundamental class is independent of the choice of the perfect obstruction
theory within a quasi--isomorphism class. That is, if $F^\bullet$ is another perfect
obstruction theory and $\psi:F^\bullet\rightarrow E^\bullet$ a quasi--isomorphism, $\psi$
naturally induces the identity map for the virtual fundamental classes associated to
$F^\bullet$ and $E^\bullet$ (\cite[Proposition 5.3]{bf97}).

By abuse of notation, we will often write $[Y]^{\virt}$ for the virtual
fundamental class $[Y,E^\bullet]$ when it is understood which obstruction theory
is used.
\end{rema}

\subsection{The obstruction complex for the definition of GW invariants}
\mylabel{subsec:obstrcplxGW}

We will now describe the obstruction theory used for the definition
of the Gromov--Witten invariants of a smooth projective complex variety $X$. Moreover, if there
is an action by a torus $T_N$ on the variety $X$, this obstruction theory will be $T_N$--equivariant.

Let $X$ be a smooth projective complex variety, $A\in H_2(X;\Z)$ an integral degree--2
homology class of $X$, and $\gmodsp$ the corresponding moduli stack of stable $m$--marked
genus--$g$ maps to $X$. Let $\pi: \ucgmodsp \longrightarrow \gmodsp$ be the universal
curve, and let $x_i: \gmodsp \longrightarrow \ucgmodsp$ ($i=1,\ldots, m$) be the marked
point sections. We will denote by $D$ the divisor defined by the images of the marked point
sections $x_i$. If no confusion can arise, we will also use the notation $\cM:=\gmodsp$ and 
$\cC:=\ucgmodsp$.  We will consider the following complex:
\[ E^\bullet:=R\pi_*\left(\left(f^*\Omega^1_X[1] \oplus \Omega^1_{\cC/\cM}(D)\right)\stackrel{L}{\otimes}
\omega_{\cC/\cM}\right), \]
where $\omega_{\cC/\cM}$ is the relative dualizing sheaf.
We will first show that there is a canonical morphism $\phi: E^\bullet 
\longrightarrow L^\bullet_\cM$ and then prove that this morphism is an obstruction
theory.

Remember that Behrend and Fantechi have given an obstruction theory for the problem
relative to the stack of prestable curves:

\begin{theo}[\cite{beh97, bf97, bm96}]
Let $p: \gmodsp \rightarrow {\mathfrak M}_{g,m}$ be the ca\-noni\-cal
mor\-phism from the stack of stable maps to the stack of prestable curves
given by forgetting the map and retaining the curve without
stabilizing. Then $\gmodsp\rightarrow {\mathfrak M}_{g,m}$ is an
open substack of a relative space of morphisms, hence it has a
relative obstruction theory which is given by
\[ \psi: \left(R\pi_*f^*TX\right)^{\vee} \longrightarrow L^\bullet_{\gmodsp/{\mathfrak M}_{g,m}}.
\]
Here $\pi: \ucgmodsp \longrightarrow \gmodsp$ is the universal curve and $f:{\mathcal C}
\rightarrow X$ is the universal stable map. \proofend
\end{theo}

\begin{coro}
There exists a canonical morphism $\phi: E^\bullet \longrightarrow L^\bullet_\cM$ 
in the derived category induced by the morphism $\psi$.
\end{coro}

\begin{proof}
Consider the following cartesian diagram where $\pi: {\mathcal C} \longrightarrow {\mathcal M}$
is a stable map to $X$, and $p:{\mathcal M} \longrightarrow S$ is the forgetting
map, \ie. $Z\longrightarrow S$ is a prestable curve:
\[
\xymatrix{
Z \ar[d]^\sigma &
{\cC} \ar[l]_\tau \ar[r]^f \ar[d]^\pi &
X\\
S &
{\cM.} \ar[l]_p&
}
\]
Remember that if we have two morphisms of schemes (or stacks) $U\stackrel{h}{\longrightarrow}V\longrightarrow W$
we get a distinguished triangle of cotangent complexes:
\[ h^*L^\bullet_{V/W} \mapright L^\bullet_{U/W} \mapright L^\bullet_{U/V} \mapright h^*L^\bullet_{V/W}[1].\]
Moreover $f^*\Omega^1_X=f^*L^\bullet_X$ naturally maps to $L^\bullet_{\mathcal C}$, so we get the following
diagram:
\[
\xymatrix{
f^*\Omega^1_X \ar[r] &
L^\bullet_{\cC} \ar[r] \ar[d] & 
L^\bullet_{\cC/\cM} \ar[r]^{\sim} &
\tau^* L^\bullet_{Z/S} \ar[d]\\
&
L^\bullet_{\cC/Z} \ar[d]^{\sim}&
&
\tau^*\sigma^*L^\bullet_S[1] \ar[d]^{\sim}\\
&
\pi^*L^\bullet_{\cM/S} \ar[rr] &
&
\pi^*p^*L^\bullet_S[1].
}
\]
This diagram is in fact commutative since $\sigma$ is flat, and so by \cite[(9.2.5)]{lm92} we have
\[ \tau^* L^\bullet_{Z/S} \oplus \pi^*L^\bullet_{\cM/S} \MRE{\sim} L^\bullet_{\cC/S}, \]
and the morphisms in the diagram above are just the morphism induced by the distinguished triangle
\[ \pi^* p^* L^\bullet_S \mapright L^\bullet_\cC \mapright L^\bullet_{\cC/S} \mapright \pi^*p^*L^\bullet_S[1].\]
Applying the cut--off functor $\tau_{\ge 0}$ to $f^*L^\bullet_X \longrightarrow L^\bullet_{\cC/\cM}$
and taking the mapping cone yields the following diagram in the derived category:
\[ 
\xymatrix{
f^*\Omega^1_X \ar[r] \ar[d] &
\Omega^1_{\cC/\cM}(D) \ar[r] \ar[d] &
f^*\Omega^1_X[1] \oplus \Omega^1_{\cC/\cM}(D) \ar[d]\\
\pi^*L^\bullet_{\cM/S} \ar[r] &
\pi^*p^*L^*_S[1] \ar[r] &
\pi^* L^\bullet_\cM [1].
}
\]
The projection formula yields the desired morphism $\phi$.
\end{proof}

\begin{prop}\mylabel{claim:obstrthy}
Let $\omega_{\cC/\cM}$ be the relative dualizing sheaf of $\pi:\cC\rightarrow \cM$.
The morphism
\[ \phi: E^\bullet \mapright L^\bullet_\cM\]
is a perfect obstruction theory for the moduli stack of stable maps $\gmodsp$. If there is a torus $T_N$ acting
on $X$, this obstruction theory is $T_N$--equivariant.
\end{prop}

\begin{proof}
First we will construct 
a two--term resolution of $R\pi_*E^\bullet$ that is $T_N$--equiva\-riant if such an
action exists on $X$.  We will use similar arguments as Behrend does for $R\pi_*f^*TX$
(\cf. \cite[Proof of Proposition 5]{beh97}). Let $M$ be an ample invertible sheaf on
$X$ and let $L=\omega_{\cC/\cM}(D)\otimes f^*M^{\otimes3}$. Then by \cite[Proposition 3.9]{bm96},
for $N$ sufficiently large and $V$ a vector bundle on $\cC$ we have that
\begin{enumerate}
\item $\pi^*\pi_*(V\otimes L^{\otimes N})\longrightarrow V\otimes L^{\otimes N}$ is surjective,
\item $R^1\pi_*(V\otimes L^{\otimes N})=0$,
\item for all $s\in S$ we have that $H^0(C_s, L_s^{\otimes -N})=0$.
\end{enumerate}
Let us set 
$ F:= \pi^*\pi_*(f^*\Omega^1_X \otimes L^{\otimes N})\otimes L^{\otimes -N}$ and 
$H:=\ker (F\longrightarrow f^*\Omega^1_X)$, and 
consider the complexes (\cf. \cite[section 4]{lt98a}) indexed at $-1$ and $0$
\[ \cA^\bullet = [H\otimes\omega_{\cC/\cM} \longrightarrow 0] \quad \text{and} \quad
\cB^\bullet = [F\otimes\omega_{\cC/\cM} \longrightarrow \Omega^1_{\cC/\cM}(D)\otimes\omega_{\cC/\cM}],\]
where the morphism within the complex $\cB^\bullet$ is induced from the composition map $F \longrightarrow f^*\Omega^1_X \longrightarrow
\Omega_{\cC/\cM}(D)$. Hence there are morphisms
\[ R^1\pi_* (F\otimes\omega_{\cC/\cM}) \mapright R^1\pi_* (f^*\Omega^1_X\otimes\omega_{\cC/\cM}) \MRE{\alpha}
R^1\pi_*\cB^\bullet \]
where $\alpha$ is surjective, by lemma \ref{lem:stableinj} and duality. As before we also have
\begin{eqnarray*}
H^0(C_t, F\otimes\omega_{\cC/\cM})&=&H^0(C_t, \pi_*(f^*\Omega^1_X\otimes L^{\otimes N})_t\otimes
L_t^{\otimes-N}\otimes \O_{\cM,t})\\
&=&H^0(C_t, L_t^{\otimes -N})\otimes \pi_*(\Omega^1_X\otimes L^{\otimes N})\otimes \O_{\cM,t}=0,
\end{eqnarray*}
so $R^0\pi_* (H\otimes\omega_{\cC/\cM})=R^0\pi_*(F\otimes\omega_{\cC/\cM})=0$.
Observe that the complex $\cB^\bullet$ fits into the short exact sequence
\[ 0 \mapright \Omega_{\cC/\cM}(D)\otimes\omega_{\cC/\cM} \mapright \cB^\bullet
\mapright F\otimes\omega_{\cC/\cM}[1] \mapright 0,\]
therefore we get a corresponding long exact sequence of higher direct image sheaves:
\begin{multline*}
0 \longrightarrow R^{-1}\pi_* \cB^\bullet \longrightarrow 
\underbrace{R^0\pi_*(F\otimes\omega_{\cC/\cM})}_{=0}
\longrightarrow R^0\pi_*(\Omega^1_{\cC/\cM}(D)\otimes\omega_{\cC/\cM}) \longrightarrow\\
\longrightarrow R^0\pi_*\cB^\bullet \longrightarrow R^1\pi_*(F\otimes\omega_{\cC/\cM}) 
\MRE{\text{surj.}} R^1\pi_*(\Omega^1_{\cC/\cM}(D)\otimes\omega_{\cC/\cM})\longrightarrow\\
\longrightarrow R^1\pi_* \cB^\bullet \longrightarrow 0.
\end{multline*}
Hence $R^i\pi_* \cB^\bullet = 0$ for $i\neq 0$. Moreover, since $R^i\pi_* (H\otimes\omega_{\cC/\cM}) = 0$
for $i\neq 1$, we also get $R^i\pi_* \cA^\bullet = R^{i+1}\pi_* (H\otimes\omega_{\cC/\cM})= 0$ for $i\neq 0$.
Now note that these two complexes fit into the following short exact sequence:
\[ 0 \mapright \cA^\bullet \mapright \cB^\bullet \mapright
(f^*\Omega^1_X[1]\oplus\Omega^1_{\cC/\cM}(D))\otimes\omega_{\cC/\cM} \mapright 0,\]
yielding the long exact sequence
\[ 0 \mapright h^{-1}(E^\bullet) \mapright R^0\pi_*\cA^\bullet \mapright R^0\pi_*\cB^\bullet
\mapright h^0(E^\bullet) \mapright 0.\]
Thus we have found a two--term resolution of $E^\bullet$ by locally free sheaves:
\[ E^\bullet \cong [R^0\pi_* \cA^\bullet \longrightarrow R^0\pi_* \cB^\bullet].\]
Moreover, the entire construction is $T_N$--equivariant, so we actually have found a $T_N$--equivariant
resolution of $E^\bullet$, if such an action exists on $X$.

Finally, we observe that $\delta:R\pi_*(\Omega^1_{\cC/\cM}(D)\otimes\omega_{\cC/\cM})\cong p^*L^\bullet_{S}$ in the
derived category. Then by using the fact that $\psi:(R\pi_* f^*TX)^\vee\longrightarrow L^\bullet_{\cM/S}$ is an
obstruction theory for the relative problem, and by applying the five lemma we get that $h^0(\phi)$ is an
isomorphism and that $h^{-1}(\phi)$ is surjective:
\[
\xymatrix@C-0.5ex{
0 \ar[r] \ar[d]_{=} &
h^{-1}(E^\bullet) \ar[r] \ar[d]_{\phi^{-1}} &
\left(R^1\pi_*f^*TX\right)^\vee \ar[r] \ar[d]_{\psi^{-1}} &
R^0\pi_*(\Omega_{\cC/\cM}(D)\otimes\omega_{\cC/\cM}) \ar[r] \ar[d]_{\delta^0}^{\sim} & \\
0 \ar[r] &
h^{-1}(L^\bullet_{\cC}) \ar[r] &
h^{-1}(L^\bullet_{\cC/{\mathfrak M}}) \ar[r] &
h^0(p^*L^\bullet_{{\mathfrak M}}) \ar[r] &\\
\ar[r] &
h^0(E^\bullet) \ar[r] \ar[d]_{\phi^0} &
\left(R^0\pi_*f^*TX\right)^\vee \ar[r] \ar[d]_{\psi^0}^{\sim} &
R^1\pi_*(\Omega_{\cC/\cM}(D)\otimes \omega_{\cC/\cM}) \ar[d]_{\delta^1}^{\sim}\\
\ar[r] &
h^0(L^\bullet_{\cC}) \ar[r] &
h^0(L^\bullet_{\cC/{\mathfrak M}}) \ar[r] &
h^1(p^*L^\bullet_{{\mathfrak M}}).
}
\]
Hence $\phi:E^\bullet \longrightarrow L_{\gmodsp}^\bullet$ is indeed a ($T_N$--equivariant) perfect
obstruction theory for the moduli stack of stable maps $\gmodsp$.
\end{proof}

That is exactly the obstruction theory we will use for the definition of the
Gromov--Witten invariants. 

\begin{lemm}
The virtual fundamental class $[E^\bullet, \phi]$ of the obstruction theory
$\phi: E^\bullet \longrightarrow L^\bullet_\cM$ coincides with the virtual
fundamental class coming from Behrend's relative obstruction theory \cite{beh97}
$\psi: (R\pi_*f^*TX)^\vee \longrightarrow L^\bullet_{\cM/{\mathfrak M}}$.
\end{lemm} 

\begin{defi}
The obstruction theory $\phi: E^\bullet \longrightarrow L^\bullet_\cM$
is the {\em obstruction theory for the Gromov--Witten invariants}.
\end{defi}

\begin{rema}
Behrend has proven in \cite{beh97} that his relative obstruction
theory defines Gromov--Witten invariants, so the definition is good. Moreover,
there has been considerable progress in showing that the different versions
of Gromov--Witten invariants all coincide, \eg. \cite{lt98b,sie98}.
\end{rema}

\begin{proof}[Proof of the Lemma]
The equality of the two virtual fundamental classes can easily be seen by
looking at the complex $E_\bullet$ dual to $E^\bullet$,
\begin{align*}
E_\bullet&= R \underline\Hom (E^\bullet, \O_{\gmodsp})\\
&\cong [(R^0\pi_*\cB^\bullet)^\vee \mapright (R^0\pi_*\cA^\bullet)^\vee]\\
&\cong \left[R^1\pi_*[F\longrightarrow\Omega^1_{\cC/\cM}(D)]^\vee \mapright R^1\pi_*[H\longrightarrow 0]^\vee \right]
\quad \text{(by duality)}\\
&=\left[ {\underline\Ext}_\pi^1([F\longrightarrow \Omega^1_{\cC/\cM}(D)],\O_\cC)\mapright
{\underline\Ext}_\pi^1([H\longrightarrow 0],\O_\cC)\right]\\
&\cong R\pi_*{\underline\Hom} ([f^*\Omega^1_X \longrightarrow \Omega^1_{\cC/\cM}(D)], \O_\cC).
\end{align*}
Here we have used that by \cite[lemma II.3.1, proposition I.5.4]{har66} there exists
a morphism of functors
\[ \zeta: R(\pi_*\circ\underline\Hom(\phantom{M}, \O_\cC)) \mapright R\pi_* \circ R\underline\Hom(\phantom{M}, \O_\cC),\]
and that this morphism $\zeta$ is an isomorphism. For convenience we also use the notation
\[ \underline\Ext^i_\pi(\phantom{M}, \O_\cC) := R^i(\pi_*\circ\underline\Hom(\phantom{M}, \O_\cC)).\]

\noindent Therefore, the $E_i$'s fit into an exact sequence
\[ 0 \longrightarrow {\mathcal T}^0 \longrightarrow E_0 \longrightarrow E_1
\longrightarrow {\mathcal T}^1 \longrightarrow 0,\]
where the sheaves ${\mathcal T}^i$ are given by taking cohomology of $E_\bullet$:
\mycomment{eq:defnti}
\begin{equation}\label{eq:defnti}
{\mathcal T}^i = \underline\Ext^i_\pi([f^*\Omega_X^1 \rightarrow
\Omega_{{\mathcal C}_{g,m}^A(X)/\gmodsp}(D)], \O_{{\mathcal C}_{g,m}^A(X)}).
\end{equation}
\end{proof}

\begin{rema} Contrary to \cite{lt98b}, the complex $[f^*\Omega_X^1
\rightarrow \Omega_{\cC/\cM}(D)]$ in (\ref{eq:defnti})
is indexed at $0$ and $1$, instead of $-1$ and $0$, moving the ${\mathcal T}^i$ complex
to the left.
\end{rema}

We will end this subsection with a lemma about how this obstruction theory
behaves under base change.
This lemma will be used when we pass to the fixed point components
of the torus action on the moduli space in section \ref{sec:virtnormbdl}.

\begin{lemm}\mylabel{lem:restrobth}
Let $\pi: \cC \longrightarrow \cM$ be a stable map to $X$ that is
an atlas for $\gmodsp$. Furthermore, let $\iota:\cM_\mygraph
\longrightarrow \cM$ be a subscheme, and look at the cartesian
diagram
\[
\xymatrix{
{\cC_\mygraph} \ar[r]^{\tilde\iota} \ar[d]^{\pi_\mygraph} &
{\cC} \ar[r]^f \ar[d]^{\pi} &
X\\
{\cM_\mygraph} \ar[r]^{\iota} &
{\cM}. &
}
\]
Let $f_\mygraph:= f\circ \tilde\iota$. Then the restrictions of
the obstruction theory $E^\bullet$ and its dual $E_\bullet$ are
given by
\begin{align*}
E^\bullet|_{\cM_\mygraph}&=R\pi_*\left( \left( f_\mygraph^*
\Omega^1_X[1] \oplus
\Omega^1_{\cC_\mygraph/\cM_\mygraph}(D_\mygraph) \right) \lotimes
\omega_{\cC_\mygraph/\cM_\mygraph}\right)\\
E_\bullet|_{\cM_\mygraph}&=R\underline\Hom_\pi([f_\mygraph^*
\Omega_X^1 \longrightarrow
\Omega^1_{\cC_\mygraph/\cM_\mygraph}(D_\mygraph) ],
\O_{\cC_\mygraph}).
\end{align*}
\end{lemm}

\begin{proof}
We will prove the lemma for the obstruction complex $E^\bullet$, the arguments for the dual
complex $E_\bullet$ are similar. There is a natural morphism
\[ R\pi_*\left( \left( f_\mygraph^* \Omega^1_X[1] \oplus \Omega^1_{\cC_\mygraph/\cM_\mygraph}(D_\mygraph)
\right) \lotimes \omega_{\cC_\mygraph/\cM_\mygraph}\right)
\mapright E^\bullet|_{\cM_\mygraph}, \] and we have to show that
this morphism is a isomorphism in the derived category, \ie. a
quasi--isomorphism between complexes. Let $\cK^\bullet :=
[f^*\Omega_X^1 \otimes \omega_{\cC/\cM} \longrightarrow
\Omega^1_{\cC/\cM}(D)\otimes\omega_{\cC/\cM}]$, indexed at $-1$
and $0$. We then have to show that
\[ \left.\left(R^i \pi_* \cK^\bullet\right)\right|_{\cM_\mygraph} = R^i{\pi_\mygraph}_* (\cK^\bullet|_{\cC_\mygraph}).\]
Now $\cK^\bullet$ fits into a short exact sequence of complexes
\[ 0 \mapright \cA^\bullet \mapright \cB^\bullet \mapright \cK^\bullet \mapright 0 \]
such that $R^i\pi_* \cA^\bullet$ and $R^i\pi_* \cB^\bullet$ are locally free and
\[ R\pi_* \cK^\bullet \cong [R^0\pi_* \cA^\bullet \longrightarrow R^0\pi_* \cB^\bullet]\]
(see above). Since $\pi$ is a proper flat morphism, we have by Grauert's continuity theorem (see for example
\cite[Th\'eor\`eme 4.12(ii)]{bs77})  that
\[ \left.\left(R^i\pi_* \cA^\bullet\right)\right|_{\cM_\mygraph} = R^i{\pi_\mygraph}_* (\cA^\bullet|_{\cC_\mygraph})
\quad \text{and} \quad \left.\left(R^i\pi_*
\cB^\bullet\right)\right|_{\cM_\mygraph} = R^i{\pi_\mygraph}_*
(\cB^\bullet|_{\cC_\mygraph}).
\]
This yields the same property for the complex $\cK^\bullet$.
\end{proof}

\subsection{Definition of the Gromov-Witten invariants}

In the previous section we have constructed a ($T_N$--equivariant) perfect obstruction theory
for the moduli stack $\gmodsp$.
Hence we get a virtual fundamental class $[\gmodsp]^{\virt}:=
[\gmodsp, E^\bullet]\in H_n(\gmodsp,\Q)$, where
$n=\rk E^\bullet$ is equal to the virtual dimension of $\gmodsp$: $n=(1-g) (\dim X -3)
+ \langle c_1(X), A\rangle +m$. So for cohomology classes $\alpha_1, \ldots, \alpha_m
\in H^*(X, \Z)$ and $\beta \in H^*(\gdmsp)$ we define the Gromov--Witten invariant
$\Psi^A_{m,g}(\beta;\alpha_1, \ldots, \alpha_k)$ by:
\[ \Psi^A_{m,g}(\beta;\alpha_1, \ldots, \alpha_m) :=
\int_{[\gmodsp]^{\text{vir}}}
\ev^*(\alpha_1\otimes\ldots\otimes\alpha_m)\wedge \pi^*\beta,\]
where $\ev: \gmodsp \rightarrow X^{\otimes m}$ is the $m$--point
evaluation map, and $\pi$  the natural
forgetting (and stabilization) morphism $\pi: \gmodsp \rightarrow \gdmsp$.

In the remaining part of the paper we will restrict ourselves to the genus--zero
case, \ie. when $g=0$, and to the invariants $\Phi^A_m$ where moreover 
the class $\beta=1$ is trivial:
\[ \Phi^A_m(\alpha_1,\ldots,\alpha_m) := \Psi^A_{m,0}(1;\alpha_1,\ldots,\alpha_m).\]
Note that for $m=3$ and $g=0$,
the Deligne--Mumford space of stable curves is a point, hence $\beta=1$ is
the only class that exists.

\section{Torus action and localization formula}\mylabel{sec:gploc} 

In this section we will sketch the construction of Graber and 
Pandharipande's localization formula for the virtual fundamental 
class (see \cite{gp97}). Let $Y$ be a Deligne--Mumford stack with 
a $\C^*$--action, admitting a $\C^*$--equi\-va\-riant perfect 
obstruction theory 
\[ \phi: E^\bullet=[E^{-1}\longrightarrow E^0] \mapright L^\bullet_Y,\] 
as for example the obstruction theory for the Gromov--Witten invariants
constructed in the previous section.

We will fix the perfect obstruction theory once and for all, and will write 
$[Y,E^\bullet]=\viry$ for the virtual fundamental class of $Y$ and $E^\bullet$. 
Let $Y_i$, $i\in {\mathcal I}$ be connected components of the fixed point set of 
the $\C^*$--action on $Y$. Consider the restriction 
of $E^\bullet$ to the fixed point components $Y_i$, 
\[ E^\bullet_i = [E_i^{-1}\otimes \O_{Y_i} \longrightarrow E^0_i\otimes \O_{Y_i}]\] 
that naturally maps to the restriction to $Y_i$ of the cotangent 
complex $L^\bullet_Y$. The restricted cotangent complex $L^\bullet_Y|_{Y_i}$ 
naturally maps to the cotangent complex $L^\bullet_{Y_i}$ of $Y_i$.

For a coherent sheaf ${\mathcal F}$ on $Y_i$ with a 
$\C^*$--action, let ${\mathcal F}=\bigoplus_{k\in\Z} {\mathcal F}^k$ be the character 
decomposition of ${\mathcal F}$ into $\C^*$--eigensheaves of $\O_{Y_i}$--modules. 
We will use the following notation for the fixed 
and the moving subsheaves: 
\[ \begin{array}{rclcl} 
{\mathcal F}^\fix&:=&{\mathcal F}^0&\text{---}&\text{the {\em fixed subsheaf}}\\ 
{\mathcal F}^\move&:=&\oplus_{k\ne0}{\mathcal F}^k&\text{---}&\text{the {\em moving subsheaf}.}\\ 
\end{array} \] 
 
\begin{lemm}[\cite{gp97}] 
The composition $\psi: E_i^{\bullet,\fix} \MRE{\phi_i^\fix} 
L^\bullet_Y|_{Y_i} \longrightarrow L^\bullet_{Y_i}$ is a perfect obstruction 
theory for $Y_i$, where $\phi_i^\fix:E_i^{\bullet,\fix}\longrightarrow 
L^\bullet_Y|_{Y_i}^\fix$ is the fixed map. \proofend
\end{lemm} 
 
\begin{defi} 
Let $Y$ be a Deligne--Mumford stack with a $\C^*$--action and a 
$\C^*$--equivariant perfect obstruction theory $\phi: E^\bullet 
\longrightarrow L^\bullet_Y$. Let $Y_i$, $i\in {\mathcal I}$ be the connected 
fixed point components of the $\C^*$--action, and let 
$\psi_i: E_i^{\bullet, \fix} \longrightarrow L^\bullet_{Y_i}$ be the 
perfect obstruction theory for $Y_i$ constructed above. We will call 
$[Y_i, E_i^{\bullet,\fix}]$ the {\em virtual fundamental class induced 
by $[Y, E^\bullet]$}, and will write $[Y_i]^\virt:= 
[Y_i, E_i^{\bullet,\fix}]$. 
\end{defi} 
 
\begin{defi} 
Let $Y_i$, $E^\bullet_i$ be as above. Let $E_{\bullet,i}=(E^\bullet_i)^\vee$ 
be the dual complex. We define the {\em virtual normal bundle} 
$N_i^\virt$ to $Y_i$ to be the moving part of $E_{\bullet,i}$: 
\[
N_i^\virt := E_{\bullet,i}^\move. 
\]
\end{defi} 
 
Note that $\rk N_i^\virt = \rk E^\bullet|_{Y_i} - \rk E^\bullet_i$, hence the 
rank of the virtual normal bundle is constant on each fixed point component. Since 
moreover the virtual normal bundle has no fixed subbundle under the $\C^*$--action, 
its equivariant Euler class exists: 
\[ e^{\C^*}([N_{0,i}^\virt \longrightarrow N_{1,i}^\virt]):= e^{\C^*}(N_{0,i}^\virt-N_{1,i}^\virt).\]

We are now able to formulate Graber and Pandharipande's 
localization theorem for the virtual fundamental class: 

\begin{theo}[\protect{Localisation formula \cite{gp97}}] 
Let $Y$ be an algebraic stack with a $\C^*$--action that can be $\C^*$--equivariantly embedded 
into a non--singular Deligne--Mumford stack. Let $\phi: E^\bullet \longrightarrow L_Y^\bullet$ 
be a $\C^*$--equivariant perfect obstruction theory for $Y$, and let $[Y,E^\bullet]$ and 
$[Y_i, E_i^\bullet]$ be the virtual fundamental classes of $Y$ and $E^\bullet$, and of the fixed 
point components $Y_i$ and the induced perfect obstruction theories $E_i^\bullet$, respectively. 
Then 
\[ [Y,E^\bullet] = \iota_* \sum_i \frac{[Y_i, E_i^\bullet]}{e^{\C^*}(N_i^\virt)},\] 
where $N_i^\virt$ is the virtual normal bundle to $Y_i$ defined above. 
\proofend
\end{theo} 
 
As a corollary we get the virtual Bott residue formula:
 
\begin{coro}[Virtual Bott residue formula \cite{gp97}] 
Let $G$ be a $\C^*$--equivariant vector bundle on $Y$, of 
rank equal to the virtual dimension of 
$Y$, $\rk G=\dim \viry=\rk E^\bullet$. Then the following 
{\em virtual Bott residue formula} holds: 
\mycomment{eq:bottres} 
\begin{equation}\label{eq:bottres} 
\int_{\viry} e(G) = \sum_{i\in {\mathcal I}}\int_{\virgg{Y_i}} 
\frac{e^{\C^*}(G_i)}{e^{\C^*}(N_i^\virt)} 
\end{equation} 
in the localized ring $A^{\C^*}(Y)\otimes \Q[\mu,\frac{1}{\mu}]$, 
where the bundles $G_i$ are the pullbacks of $G$ under 
$Y_i\hookrightarrow Y$. \proofend 
\end{coro} 
 
\begin{rema} 
Note that the formula indeed makes sense: since $\rk G=\dim \viry$ we actually have 
\[ \int_{\viry} e(G) = \int_\viry e^{\C^*} (G). \] 
In particular, the right hand side of equation (\ref{eq:bottres}) takes values in $\Q$, 
not just in a  polynomial ring over $\Q$. 
\end{rema} 
 
\begin{rema} 
Note that we can replace in all statements above the 
one--di\-men\-sional torus $\C^*$ by a higher dimensional torus 
$(\C^*)^d$. In fact, if we diagonalize the $(\C^*)^d$--action we 
get $d$ commutative $\C^*$--actions. We thus can apply the 
localization formula $d$ times, to get the statement for the 
$(\C^*)^d$--action. 
\end{rema}

\section{Preliminaries on toric varieties}\mylabel{sec:torvar} 
 
This section will mostly serve to  remind the reader of the 
definition and some properties of smooth toric varieties as well 
as to fix the notation.
Of course, everything is already well known, see for example 
(in alphabetic order) \cite{aud91,bat93,cox97,dan78,del88,ful93,oda88}.

\subsection{The algebro--geometric construction of toric varieties} 
 
For all what follows we will fix the following notation: Let $d>0$ be a positive 
integer. Let $N=\Z^d$ be the $d$--dimensional integral lattice, and $M=\Hom(N, \Z)$ 
be its dual space. Moreover, let $N_\R=N\otimes_\Z\R$ and $M_\R=M\otimes_Z\R$ 
be the $\R$--scalar extensions of $N$ and $M$, respectively. 
 
A convex subset $\sigma\subset N_\R$ is called 
a {\em regular $k$--dimensional cone} if there exists 
a $\Z$--basis $v_1, \ldots, v_k, \ldots, v_d$ of $N$ such that the cone $\sigma$ is 
generated by $v_1,\ldots, v_k$. The vectors $v_1,\ldots, v_k\in N$ are the 
{\em integral generators of $\sigma$}. The origin $0\in N_\R$ is the only 
{\em regular zero dimensional cone}. Its set of integral generators is empty.
A {\em face} of a regular cone $\sigma$ is a cone $\sigma'$ generated by a subset of 
the integral generators of $\sigma$. If $\sigma'$ is a (proper) face of $\sigma$, we will 
write $\sigma'\prec\sigma$. 
 
A finite system $\Sigma=\{\sigma_1, \ldots, \sigma_s\}$ of regular cones in 
$N_\R$ is called a regular {\em $d$--dimensional fan} of cones, 
if any face $\sigma'$ of a cone $\sigma\in\Sigma$ in the fan and any intersection
of two cones $\sigma_1, \sigma_2\in\Sigma$ are again in the fan.
A fan $\Sigma$ is called a {\em complete fan} if the (set theoretic) 
union of all cones $\sigma_i$ in $\Sigma$ is all of $N_\R$, \ie. 
$N_\R= \bigcup_i \sigma_i$.
The {\em $k$--skeleton} $\Sigma^{(k)}$ of the fan $\Sigma$ 
is the set of all $k$--dimensional cones in $\Sigma$. 
 
By abuse of language, we will also consider cones $\sigma$ as fans, meaning 
in fact the fan $\Sigma_\sigma$ of $\sigma$ and all its faces: 
$ \Sigma_\sigma = \{ \sigma' \, | \, \sigma'\preccurlyeq\sigma\}$.
 
A subset ${\primcol} \subset \Sigma^{(1)}$ of the $1$--skeleton of a fan $\Sigma$ is 
called a {\em primitive collection} of $\Sigma$ (see \cite{bat91}) 
if ${\primcol}$ is not the set of generators of a cone in $\Sigma$, 
while any proper subset of $\primcol$ is. We will denote the set of primitive 
collections of $\Sigma$ by $\primcolset$. 
 
Let $n=\module{\smash{\Sigma^{(1)}}}$ be the cardinality of the one--skeleton of 
$\Sigma$, and $v_1, \ldots, v_n$ its elements. Let $z_1, \ldots, z_n$ be a set of 
coordinates in $\C^n$ and let 
$\iota: \C^n \longrightarrow N\otimes_\Z\C$ be a linear map such that $\iota(z_i)=v_i$. 
For each primitive collection $\primcol\in\primcolset$, 
$\primcol=\{v_{i_1},\ldots, v_{i_p}\}$, 
we define an $(n-p)$--dimensional affine subspace in $\C^n$ by 
\[ \AA(P) := \{(z_1,\ldots, z_n)\in \C^n \, |\, z_{i_1}=\ldots=z_{i_p}=0\}.\] 
Moreover, we define the set $U(\Sigma)$ to be the open algebraic subset of $\C^n$ 
given by 
\[ U(\Sigma) =\C^n - \bigcup_{\primcol\in\primcolset} \AA(P).\] 
 
The map $\iota: \C^n \longrightarrow N_\C$ induces a map between tori 
$(\C^*)^n \longrightarrow (\C^*)^d$ that we will also call $\iota$. Here, $\C^*=\C-\{0\}$. 
Let $\DD(\Sigma) := \ker (\iota: (\C^*)^n \longrightarrow (\C^*)^d)$ be the kernel of
this map, an $(n-d)$--dimensional subtorus.
 
\begin{defi}\mylabel{def:vt} 
Let $\Sigma$ be a regular $d$--dimensional\footnote{A $d$--dimensional fan is a fan 
in $\Z^d$ containing a cone of dimension $d$.} fan of regular cones. 
The quotient 
$\vt := U(\Sigma) / \DD(\Sigma)$
is called the {\em toric manifold associated with $\Sigma$}.  
\end{defi} 
 
The following proposition provides us with an atlas of charts for toric manifolds.

\begin{prop}\mylabel{thm:homcoords} 
Let $\sigma\in\Sigma^{(k)}$, and let $\{v_{i_1},\ldots, v_{i_k}\}$ be its 
set of generators. Let $\{v_{i_1},\ldots, v_{i_d}$\} be a $\Z$--basis of $N=\Z^d$ completing the 
set of generators of $\sigma$, and let $u_1,\ldots, u_d$ be its dual basis of $M=\Hom(N,\Z)$. 
Define the open subset $V(\sigma)\subset\C^n$ by 
\[ V(\sigma)=\left\{(z_1,\ldots, z_n)\,|\,z_j\neq0 \text{ for $j\notin \{i_1,\ldots, i_k\}$}\right\}.\] 
 
These open sets $V(\sigma)$ satisfy the following properties: 
\begin{enumerate} 
\item $U(\Sigma)=\bigcup_{\sigma\in\Sigma^{(d)}} V(\sigma)$; 
\item if $\sigma'\prec\sigma$, then $V(\sigma')\subset V(\sigma)$; 
\item $V(\sigma)$ is isomorphic to $\C^k\times(\C^*)^{n-k}$, and 
the torus $\DD(\Sigma)$ acts freely on $V(\sigma)$. The quotient $U_\sigma:=V(\Sigma)/\DD(\Sigma)$ 
is the toric subvariety associated to the cone $\sigma\in\Sigma$, whose co--ordinate 
functions $x_1^\sigma,\ldots, x_d^\sigma$ are the following Laurent monomials 
in $z_1,\ldots, z_n$: 
\[ x_j^\sigma=z_1^{\langle v_1, u_j\rangle}\cdots z_n^{\langle v_n,u_j\rangle}.\] 
\end{enumerate} 
\end{prop} 
 
\begin{rema} 
Note that our notation is slightly different to Batyrev's in \cite{bat93}: he defines the open 
sets $U(\Sigma)$ just for (complete) fans, while he calls $U(\sigma)$ what we call $V(\sigma)$. 
\end{rema} 
 
From now on we will consider at complete regular fans $\Sigma$ of regular cones.
 
\subsection{Support functions of a fan and dual polyhedra} 
 
A continuous function $\varphi : N_\R \mapright \R$ is called 
{\em $\Sigma$--piecewise linear}, if $\varphi$ is linear 
on every cone of $\Sigma$. Let $PL(\Sigma)$ be the set of all 
$\Sigma$--piecewise linear functions. Note that $PL(\Sigma)\cong \R^n$ since
$\Sigma$--piecewise linear functions are given by their values on the 
$1$--skeleton of $\Sigma$.

Such a function 
$\varphi\in PL(\Sigma)$ is called {\em upper convex} if for any $x,y\in N_\R$,
$\varphi(x)+\varphi(y)\ge\varphi(x+y)$. If moreover for any two different
$d$--dimensional cones $\sigma_1, \sigma_2\in\Sigma$, the restrictions
$\varphi_{|\sigma_1}$ and $\varphi_{|\sigma_2}$ are different 
linear functions, then $\varphi$ is called {\em strictly upper convex support
function} for $\Sigma$.

\begin{prop}\mylabel{prop:sucsp1} 
A $\Sigma$--piecewise linear function $\varphi$ is a strictly upper convex 
support function if and only if for all primitive collections $\primcol\in\primcolset$, 
$\primcol=\{v_{i_1},\ldots,v_{i_k}\}$, the following inequality holds: 
\[ \varphi(v_{i_1})+\ldots+\varphi(v_{i_k}) > \varphi(v_{i_1}+\ldots+ v_{i_k}).\] 
\end{prop} 
 
We will give another criterion in terms of convex polytopes that will 
be useful in particular for the construction via a moment map: 
 
\begin{prop}\mylabel{prop:sucsp2} 
Let $\Sigma$ be a complete, regular fan in $N=\Z^d$. 
Let $\varphi\in PL(\Sigma)$ be a $\Sigma$--piecewise linear function on $\Sigma$. 
Define a polytope $\Delta_\varphi\in M$  by 
\[ \Delta_\varphi = \{ m \in M_\R\, | \, \langle m,n\rangle \ge -\varphi(n), \,\,\forall n\in N\, \}.\] 
 
Then the function $\varphi$ is a strictly upper convex support function if and 
only if the integral convex polytope $\Delta_\varphi$ is $d$--dimensional 
and has exactly $\{l_\sigma \, |\, \sigma\in\Sigma^{(d)}\}$ as the set of 
its vertices. Here, the $l_\sigma\in M_\R=\Hom(N, \R)$ are given by 
$l_\sigma=\varphi_\sigma$.
\end{prop} 
 
\subsection{Divisors, cohomology and first Chern class} 

The cohomology of a toric manifold $X_\Sigma$ is generated by its $T^N$--invariant
divisors $D_1, \ldots, D_n$ that are given by $D_i=(\{z_i=0\}\cup U(\Sigma))/\DD(\Sigma)$.
To each $\Sigma$--piecewise 
linear functions $\varphi: N_\R \longrightarrow \R$ we can associate a divisor by setting 
$ D_\varphi = \sum_{i=1}^n \varphi(v_i) D_i$, yielding a canonical isomorphism 
$H^2(\vt,\R)\cong PL(\Sigma)/M_\R$.

The cohomology ring of $\vt$ is therefore the quotient of $\R[Z_1, \ldots, Z_n]$ 
by an ideal of relations. As we have seen above, the ideal of linear relations is
$\Lin(\Sigma):=\langle \sum_i u_1(v_i)Z_i, \ldots, \sum_i u_d(v_i) Z_i \rangle$, where
$u_1,\ldots ,u_d$ is some basis of $M=\Hom(N, \Z)$. 
The higher--degree relations in the cohomology ring are given by the so--called
Stanley--Reisner ideal 
$SR(\Sigma):= \left\langle\prod_{v_j\in\primcol} Z_j\right\rangle_{\primcol\in\primcolset}$.

\begin{prop}\mylabel{prop:homcomtor}
The cohomology ring of the compact toric manifold $\vt$ is canonically isomorphic
to the quotient of $\R[Z_1,\ldots, Z_n]$ by the ideal $\Lin(\Sigma)+SR(\Sigma)$:
\[ H^*(\vt, \R) \cong \R[Z_1, \ldots, Z_n]/(\Lin(\Sigma)+SR(\Sigma)).\]
Moreover, the first Chern class $c_1(\vt)$ of $\vt$ is represented by $Z_1+\cdots+Z_n$.
\end{prop} 

Dually, let $R(\Sigma)\subset \Z^n$ be the subgroup of $\Z^n$ defined by 
\[ R(\Sigma) = \{ (\lambda_1,\ldots,\lambda_n) \in \Z^n \,|\, 
\lambda_1 v_1+\ldots+\lambda_n v_n =0\} \cong \Z^{n-d}.\] 
Then the group $R(\Sigma)_\R=R(\Sigma)\otimes_\Z\R$ of $\R$--linear 
extensions of $R(\Sigma)$ is canonically isomorphic to $H_2(\vt,\R)$.
 
The pairing $H^2(\vt,\R)\otimes H_2(\vt,\R)\longrightarrow \R$ lifts to 
$PL(\Sigma)\otimes R(\Sigma)_R$ and is given there by the degree map: 
\[ \deg_\varphi(\lambda) = \sum_{i=1}^n \lambda_i\varphi(v_i).\]

\subsection{Toric manifolds as symplectic quotient}\mylabel{sec:vtsq} 
 
\begin{defi}\mylabel{def:kcone} 
As before let $\Sigma$ be a complete, regular cone in $N$. 
Denote by $K(\Sigma)$ the cone in 
$H^2(\vt, \R)\cong PL(\Sigma)/M_\R$ consisting of the classes of all 
upper convex support function $\varphi$ for $\Sigma$. We denote by 
$K^\circ(\Sigma)$ the interior of $K(\Sigma)$, \ie. the cone consisting 
of the classes of all strictly convex upper support functions in $PL(\Sigma)$. 
\end{defi} 
 
\begin{prop} 
The open cone $K^\circ(\Sigma)\subset H^2(\vt,\R)$ consists of 
classes of K\"ah\-ler $(1,1)$--forms on $\vt$, \ie. $K(\Sigma)$ is 
isomorphic to the closed K\"ahler cone of $\vt$. 
\end{prop} 

If the K\"ahler cone is non--empty, the toric manifold can be constructed
as a symplectic quotient as follows. 
The $n$--dimensional complex space $\C^n$ has a natural 
symplectic structure. 
Remember from above, that $\DD(\Sigma)$ is an algebraic subtorus of $(\C^*)^n$, 
thus acting on $\C^n$. Let $G\cong (S^1)^{n-d}$ be the maximal compact subgroup 
of $\DD(\Sigma)$. Since $\DD(\Sigma)\subset(\C^*)^n$ acts as a subtorus, so does 
$G \subset T^n$. The action of $G\subset T^n$ is naturally Hamiltonian, and we 
obtain its moment map $\mu$ by composing the moment map $\mu_{T^n}$ of the 
$n$--dimensional torus action on $\C^n$ with the restriction map 
$\beta^*: (\gott^n)^* \longrightarrow \gotg^* $: 
\[ \mu: \C^n \MRE{\mu_{T^n}} (\gott^n)^* \MRE{\beta^*} \gotg^*.\] 
 
For almost all $\xi\in \gotg^*$, the moment map is regular. Moreover,
the action of $G$ on the level set $\mu^{-1}(\xi)$ is effective if and only 
if $\mu^{-1}(\xi)\subset U(\Sigma)$, the open subset of $\C^n$ used for the 
algebro--geometric quotient. 
 
\begin{theo}[\cite{del88}] 
Let $\vt$ be a projective simplicial toric variety. Then there exists a regular 
value $\xi\in\gotg^*$ of the moment function $\mu:M\longrightarrow \gotg^*$ such 
that the level set $\mu^{-1}(\xi) \subset U(\Sigma)$ is in the effective 
subset of the action $G$, and there is a diffeomorphism 
\[ \mu^{-1}(\xi)/G \mapright \vt \] 
preserving the cohomology class of the symplectic form. \proofend
\end{theo} 

\begin{prop} 
Let $\varphi$ be the strictly upper convex support function associated
with the symplectic form $\omega_\xi$ of the quotient $\mu^{-1}(\xi)/G$.
Then the polytope $\Delta_\varphi$ is the moment polytope of the induced 
Hamiltonian $T^N$--action on $\vt$. \proofend
\end{prop} 
    
\section{Torus action and its fixed points in $\vt$ and $\gmodsp$}\mylabel{sec:tvaction} 
 
We have seen earlier, that a toric variety $\vt$ has by definition an algebraic torus 
acting on it. In fact, it contains an algebraic torus $K\cong(\C^*)^d$ as open and 
dense subset. This ``big torus'' acts on itself by the usual group multiplication, 
and extends naturally to the rest of $\vt$. 
In general, by pull back through the universal stable map $f:\ucmodsp \longrightarrow X$, 
an action on a manifold $X$ induces an action on the moduli spaces $\modsp$ of stable 
maps to $X$. 
 
In this section, we will study these actions to determine the fixed point components 
in the moduli spaces $\modsp$. Although we will restrict ourselves to genus--zero stable maps, 
it is possible to carry out a similar analysis for higher genus stable maps 
to toric varieties, \cf. Graber and Pandharipande's 
analysis in \cite{gp97} for projective spaces $\P^d$.

\subsection{The torus action on $\vt$ and its fixed points} 
 
As with any set on which a group acts, the toric variety $\vt$ is a disjoint union 
of its orbits (\cf. \cite[chapter 3]{ful93} for details and proofs of the 
following statements). Here again, toric varieties are very nice objects to study: 
for each cone $\sigma$ in a regular fan $\Sigma$, there is exactly one such orbit
$O_\sigma$. Moreover, 
\[ O_\sigma \cong (\C^*)^{n-k} \qquad \text{where $\dim \sigma = k$}. \] 
The orbits $O_\sigma$ are an open subvariety of its closure in $\vt$, which we 
denote by $V_\sigma$. The $V_\sigma$ are closed subvarieties of $\vt$. 
The following proposition expresses the relations between these set; for a 
proof see for example \cite{ful93}. 
 
\begin{prop} 
There are the following relations among orbits $O_\sigma$, orbit closures 
$V_\sigma$, and the affine open sets $U_\sigma$: 
\begin{enumerate} 
\item $U_\sigma = \coprod_{\tau\preccurlyeq\sigma} O_\tau$; 
\item $V_\sigma = \coprod_{\gamma\succcurlyeq\sigma} O_\gamma$; 
\item $O_\sigma = V_\sigma - \bigcup_{\gamma\succ\sigma} V_\gamma$. 
\end{enumerate} \proofend
\end{prop} 
 
In fact, the orbit closures $V_\sigma$ are the $T_N$--invariant divisor $D_i$
defined earlier, or intersections of such divisors. When using the quotient 
construction $\vt=U(\Sigma)/\DD(\Sigma)$ from a (complete) regular fan $\Sigma$, 
one can easily describe the orbit closures $V_\sigma$ as follows: Let 
the $k$--cone $\sigma\in\Sigma$ be given by the set $\{v_{i_1},\ldots,v_{i_l}\}$. 
Then the closed subvariety $V_\sigma$ is the quotient of the set 
\[ Z_\sigma := \left\{ (z_1, \ldots, z_n) \in U(\Sigma)\subset \C^n\,|\, z_{i_1}=\ldots=z_{i_k}=0 \right\} \] 
by the action of the torus $\DD(\Sigma)\cong (\C^*)^{n-d}$. In 
particular, this description gives a useful characterization of 
$V_\sigma$ as subvariety of $\vt$. 
 
In the next section we will be especially interested in such closed subspaces 
$V_\sigma$ that are of dimension zero and one, \ie. fixed points of the $T_N$--action 
on $\vt$, and invariant curves. In a compact toric variety, the latter are 
always isomorphic to $\P^1$, as the closed subvarieties $V_\sigma$ are itself 
toric varieties again, and since $\P^1$ is the only compact one--dimensional 
toric variety. These $T_N$--invariant curves are in a one--to--one correspondence 
to $(d-1)$--dimensional cones, while fixed points are in a one--to--one relation 
to $d$--dimensional cones.

\subsection{Fixed points of the induced torus action on the moduli space} 

To find out how the fixed points of the induced torus action on the 
moduli stack look like, let us consider first a single stable map 
$(C; x_1, \ldots, x_m; f) \in \modsp$, \ie. a stable map 
\[\xymatrix{ 
C \ar[d]^{\pi} \ar[r]^f&\vt\\ 
\Spec \C. \ar@(ul,ur)[u]^-{x_i} 
}\] 
Let $C=C_1\cup\ldots\cup C_k$ be the decomposition of the curve $C$ into 
irreducible and reduced curves $C_i$. Since we only look at rational curves 
$C$, the irreducible and reduced components $C_i$ of $C$ are all rational 
as well, that is, they are isomorphic to $\P^1$. 
 
\begin{lemm} 
The stable map $(C; x_1, \ldots, x_m; f)$ is a fixed point of the induced 
action of $T_N$ on the moduli stack of stable maps $\modsp$ if and only if 
it satisfies all of the following conditions: 
\begin{enumerate} 
\item \mylabel{lem:fp1} All special points of $C$, that is the marked points $x_1, \ldots, 
x_m$ and the intersection points $C_i\cap C_j$, $i\neq j$ of two different 
irreducible and reduced components, are mapped to fixed points of the $T_N$--action 
on $\vt$; 
\item \mylabel{lem:fp2} If $C_i$ is an irreducible and reduced component of $C$ that is mapped 
to a point by $f$, then it is mapped to a fixed point of the $T_N$--action 
on $\vt$; 
\item \mylabel{lem:fp3} If an irreducible and reduced component $C_i$ of $C$ is not mapped 
to a point by $f$, it is mapped to one of the $T_N$--invariant subvarieties 
$V_\sigma\subset \vt$ of dimension one, corresponding to a dimension $(d-1)$ 
cone $\sigma\in\Sigma^{(d-1)}$. 
\end{enumerate} 
\end{lemm} 
 
\begin{rema} The above lemma is a generalisation of similar results by Kontsevich 
\cite{kon95} (also see Graber and Pandharipande's \cite{gp97}) for stable maps to 
a complex projective space $\C\P^n$. 
\end{rema} 
 
\begin{proof} 
For a stable map $(C; x_1,\ldots, x_m; f)$ to be a fixed point of the 
$T_N$--action on $\modsp$ means that for any element $t\in T_N$ in the torus $T_N$, 
the stable map 
$t\cdot (C; \underline{x}; f)$ is isomorphic to the original curve 
$(C; \underline{x}; f)$, \ie. that there exists a morphism $\phi_t: C \longrightarrow C$ 
such that the following diagram is commutative (\cf. definition \ref{def:stablemap}): 
\[\xymatrix{ 
C \ar[r]^{\phi_t} \ar[d]^{\pi} \ar@(ur,ul)[rr]^{t\cdot f} & C \ar[d]^{\pi} \ar[r]^f&X\\ 
\Spec\C \ar@(ul,ur)[u]^-{x_i} \ar[r]^{=} & \Spec\C. \ar@(ul,ur)[u]^-{x_i} 
}\] 
 
Now, it is obvious that a curve $C$ satisfying the three conditions stated in the lemma is 
isomorphic to $t\cdot C$ for any $t\in T_N$, taking for $\phi_t:C\longrightarrow C$ 
the morphism defined on 
the irreducible and reduced components $C_i$ by 
\[ {\phi_t}_{|C_i} = \left\{ \begin{array}{ll} 
                \id_{C_i}\quad&\text{if $f(C_i)=\{pt.\}$}\\[4pt] 
                f^*t \quad&\text{otherwise.} 
            \end{array} \right. 
\] 
 
On the other hand, let $C$ be a fixed point of the $T_N$--action on $\modsp$. We thus 
have to show that $C$ satisfies the three conditions of the lemma. 
 
Let $x_i\in C$ be a marked point of the curve $C$. Then it is obvious that $x_i$ 
has to be mapped to a fixed point in $\vt$: since $\phi_t$ has to be constant on the 
marked points, we have 
\[ \forall t\in T_N:\, t\cdot f(x_i) = f(x_i). \] 
Now, assume that $q$ is a special point of $C$ that is not mapped to a fixed point in $\vt$. 
Then the orbit of $f(q)$ under the $T_N$--action contains certainly a subspace isomorphic 
to $\C^*$. On the other hand, the image of the special points of $C$ by $f$ is a finite 
set. Hence we obtain a contradiction, since the image of a special point under any $\phi_t$ 
is always again a special point. 
 
So if $C_i$ is an irreducible and reduced component of $C$ that is mapped to a point by $f$, 
it has to contain at least three special points by the stability condition, and thus is mapped 
to a fixed point in $\vt$ as well. 
 
Similarly, if $C_i$ is an irreducible and reduced component of $C$ that is not mapped to 
a point by $f$, and the image of which is not contained in the closure of a one--dimensional 
$T_N$--orbit $V_\sigma$, then $C_i$ contains a point whose $T_N$--orbit is at least two--dimensional. 
On the other hand, $t\cdot f(C_i)$ always has to be contained in the image $f(C)$ of $C$ by $f$ that 
is one--dimensional, hence a contradiction. 
\end{proof} 
 
Note that a (general) stable curve to $\vt$ 
\[\xymatrix{ 
C \ar[d]^{\pi} \ar[r]^f&X\\ S. \ar@(ul,ur)[u]^-{x_i} }\] is in a 
fixed point component of the $T_N$--action on the moduli stack 
$\modsp$ if and only if each geometric fiber $C_s$ is a fixed 
point, \ie. satisfies the conditions of the lemma above. 
 
Following Kontsevich's description of the fixed points of the action of $(\C^*)^d$ 
on the moduli space $\M_{g,m}^A(\P^d)$ of stable maps to projective space (\cf. \cite{kon95}), 
we will use decorated graphs to keep track of the different fixed point components in the moduli 
space $\modsp$. 
 
\begin{figure}[ht] 
\begin{center} 
\epsfig{figure=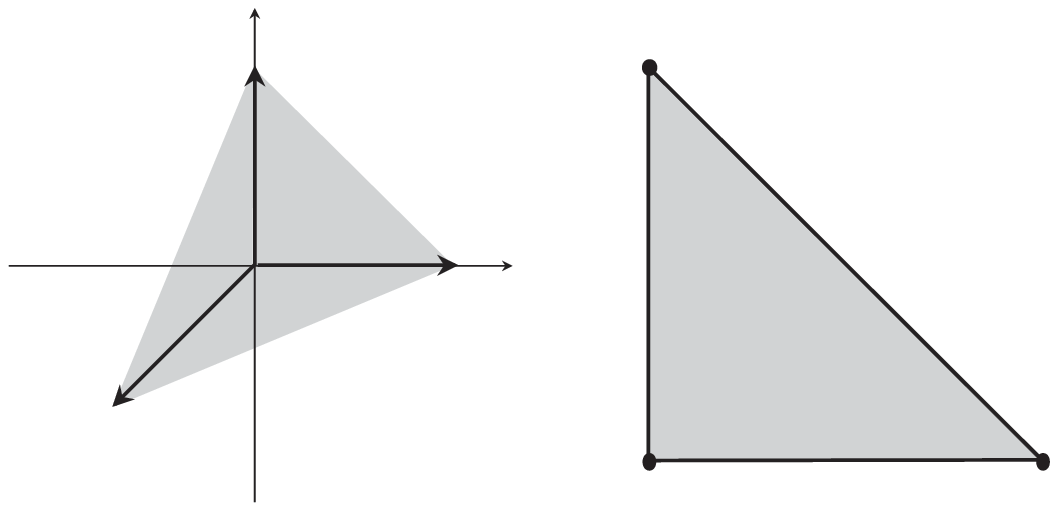,height=35mm} 
\caption{The fan and the convex polyhedron of $\C\P^2$.}\label{fig:grfap} 
\end{center} 
\end{figure}

However, before we will give the definition of the type of graphs we want to consider, let us 
look at an easy example, the moduli space $\cM_{0,m}^{A}(\C\P^2)$ of $m$--pointed stable rational maps of degree $A$ 
to the two--dimensional complex space $\C\P^2$. The fan $\Sigma$ of $\C\P^2$ and the convex polyhedron 
$\Delta_\varphi$ associated to the standard symplectic form $\varphi=c_1(\C\P^2)$ are shown in 
figure \ref{fig:grfap}. 
 
\mycomment{fig:grfap} 
 
By the previous lemma, each fixed point in the moduli space $\cM_{0,m}^{A}(\C\P^2)$ has to ``live 
on the boundary of the polyhedron $\Delta_\varphi$'', since the corners and the (one--dimensional) 
boundary components of the polyhedron correspond to fixed points respectively one--dimensional 
orbits of the torus action on $\C\P^2$. In fact, if one only looks at where the irreducible components 
and the marked points are mapped to in $\C\P^2$, one could abstractly think of such a fixed map as 
a graph that is wrapped around the polyhedron $\Delta_\varphi$. 

In the following we will continue to use the notation of fans and its dual polyhedron
$\Delta_\varphi$, \eg. we will label the vertices of the polyhedron $\Delta_\varphi$ by
maximal cones $\sigma \in \Sigma^{(d)}$, \etc. Remember that $\Delta_\varphi$ is also
the moment polytope of $(\vt, \omega_\varphi)$.

We will define three different kinds of graphs: a {\em topological $\modsp$--graph type} 
$\gtop$ describing the ``image'' of a fixed stable map on the polyhedron $\Delta_\varphi$, a
{\em $\modsp$--graph type} $\gtype$ describing moreover the image of every irreducible
component of a fixed stable map, and a {\em $\modsp$--graph} $\mygraph$ containing all the
data of a $\modsp$--graph type plus the location of the marked points.

\begin{defi} 
Let $\Sigma$ be a complete regular fan in $N\cong \Z^d$, and let $\Delta_\varphi$ be its
dual polyhedron.

A {\em topological $\modsp$--graph type} $\gtop$ is a finite one--dimensional CW--complex
with the following decorations:
\begin{enumerate}
\item A map $\sigma: \Vertex(\gtop) \longrightarrow \Sigma^{(d)}$ mapping each 
vertex\footnote{We will denote vertices with a gothic $\gotv$ to avoid confusion
with generators of cones in a fan.} $\gotv$ of the graph to a vertex $\sigma(\gotv)$
of $\Delta_\varphi$;
\item A map $d: \Edge(\gtop) \longrightarrow \Z_{>0}$, representing multiplicities of maps.
\end{enumerate}
These decorations are subject to the following compatibility conditions:
\renewcommand{\theenumi}{\alph{enumi}}\renewcommand{\labelenumi}{(\theenumi)} 
\begin{enumerate} 
\item The map $\sigma: \Vertex(\gtop) \longrightarrow \Sigma^{(d)}$ is injective;
\item If an edge $e\in \Edge(\gtop)$ connects two vertices $\gotv_1,\gotv_2\in\Vertex(\gtop)$ 
labeled $\sigma(\gotv_1)$ and $\sigma(\gotv_2)$, then the two 
cones must be different and have a common $(d-1)$--dimensional 
face: $\sigma(\gotv_1)\cap\sigma(\gotv_2)\in\Sigma^{(d-1)}$; 
\item There is at most one edge connecting any two vertices: for any two edges 
$e_1, e_2\in \Edge(\gtop)$ with vertices $\gotv_{i,1}$ and $\gotv_{i,2}$, we have
$\{\gotv_{1,1}, \gotv_{1,2}\} \neq \{\gotv_{2,1}, \gotv_{2,2}\}$;
\item The graph represents a stable map with homology class $A$: 
\[ \sum_{\begin{array}{c}\scriptstyle e\in\Edge(\gtop)\\\scriptstyle \partial e=\{\gotv_1(e),\gotv_2(e)\}\end{array}} 
d(e) [V_{\sigma(\gotv_1)\cap\sigma(\gotv_2)}] = A,\] 
where $[ V_{\sigma(\gotv_1)\cap\sigma(\gotv_2)}]$ is the homology class associated to this subvariety,
and $\partial e=\{ \gotv_1(e), \gotv_2(e)\}$ associates to an edge $e$ the 
two vertices $\gotv_1(e), \gotv_2(e)$ it connects. 
\end{enumerate} 

A {\em $\modsp$--graph type $\gtype$} is a finite one--dimensional CW-complex as above
that is subject only to compatibility conditions (b) and (d), and additionally:
\begin{itemize}
\item[(e)] The CW--complex $\gtype$ contains no loops.
\end{itemize}

A {\em $\modsp$--graph $\mygraph$} is a $\modsp$--graph type with an extra decoration:
\begin{itemize}
\item[3.] A map $S: \Vertex(\mygraph) \mapright \gotP(\{1,\ldots, m\})$ associating to each 
vertex a set of marked points; 
\end{itemize}
subject to the following additional compatibility conditions:
\begin{itemize}
\item[(f)] For any two vertices $\gotv_1,\gotv_2\in\Vertex(\mygraph)$, the sets of associated marked 
points are disjoint: $S(\gotv_1)\cap S(\gotv_2)=\emptyset$; 
\item[(g)] Every marked point is associated with some vertex: 
\[\bigcup_{\gotv\in\Vertex(\Sigma)} S(\gotv)=\{1,\ldots, m\}.\] 
\end{itemize}
The are natural maps between the different categories of graphs that we will
denote by $\typemap$ and $\topmap$:
\[
\xymatrix{
\mygraph \ar@{|->}[rr]^\topmap \ar@{|->}[dr]^\typemap && \gtop\\
& \gtype \ar@{|->}[ur]^\topmap
}
\] 
\end{defi}
 
\begin{rema}
Note, that in all cases above, there exists an induced map
\[ \Edge(\cdot) \mapright \Sigma^{(d-1)} \]
from edges of the graph to $(d-1)$--dimensional cones, or dually to edges
of the polyhedron $\Delta_\varphi$. 
\end{rema}

In the remaining part of this section, we will consider
$\modsp$--graphs $\mygraph$ only, which we will simply call graphs when 
the underlying moduli stack is understood. In fact, each of these graphs
describes a fixed point component in $\modsp$, while graph types and
topological graph types describe families of such fixed point components.

\begin{prop}\mylabel{prop:isfixedptcomp}
For a $\modsp$--graph $\mygraph$, let $\deg$ be the function assigning to each
vertex the number of its special points:
\begin{align*} 
\deg : \Vertex(\mygraph)& \mapright \Z_{>0}\\ \gotv&\DMRE{} \#S(\gotv) +
\#\{ e\in\Edge(\mygraph) | \gotv \in \partial e \}.
\end{align*} 
Furthermore, let $\cM_\mygraph$ be the following product of Deligne--Mumford
spaces:
\[ \cM_\mygraph := \prod_{\gotv \in \Vertex(\mygraph)} \dmsp[\deg(\gotv)], \]
where we formally set $\dmsp[0]=\dmsp[1]=\dmsp[2]=\{pt.\}$, the universal curve
above these space being by definition a point as well.

Then there exists a canonical family of $T_N$--fixed stable maps to $\vt$
\[ \pi: \cC_\mygraph \mapright \cM_\mygraph, \]
fitting into the following diagram:
\[
\xymatrix{ 
{\cC}_\mygraph \ar[r] \ar[d]^{\pi} \ar@(ur,ul)[rr]^{f} & 
{\cC}^A_{0,m}(\vt) \ar[d]^{\pi} \ar[r]^f&
\vt\\
{\cM}_\mygraph \ar@(ul,ur)[u]^-{x_i} \ar[r]^-\gamma & 
{\modsp}. \ar@(ul,ur)[u]^-{x_i} & 
}
\]
The image of $\cM_\mygraph$ in $\modsp$ is a fixed point component of the
$T_N$--action on $\modsp$.
\end{prop}

\begin{proof}
First we will describe the family of curves $\cC_\mygraph\longrightarrow \cM_\mygraph$.
For each edge $e\in \Edge(\mygraph)$, let us number its vertices: 
$\partial e=\{\gotv_1(e), \gotv_2(e)\}$. For each such edge $e$ we also  
fix a map of degree one\footnote{Remember that the orbit closures of type
$V_{\sigma(\gotv_1)\cap\sigma(\gotv_2)}$ are isomorphic to $\P^1$!}
$\widetilde{f}_e: \P^1 \longrightarrow V_{\sigma(\gotv_1)\cap\sigma(\gotv_2)}\subset\vt$, 
such that $\widetilde{f}_e(0)=V_{\sigma(\gotv_1)}$ and 
$\widetilde{f}_e(\infty)=V_{\sigma(\gotv_2)}$. We set $f_e:=\widetilde{f}_e\circ z^{d(e)}$
to obtain a map $f_e: \P^1 \longrightarrow V_{\sigma(\gotv_1)\cap\sigma(\gotv_2)}$ of
degree $d(e)$. Note that
up to parametrization such a map is a unique. Also set $C_e:=\P^1$.
Let $\val : \Vertex(\mygraph) \longrightarrow \Z_{>0}$ be the function assigning to 
each vertex the number of outgoing edges:
\[ \val(\gotv) := \#\{ e\in\Edge(\mygraph) | \gotv \in \partial e \}.\]
For each vertex $\gotv\in\Vertex(\gotv)$, chose an ordering of the set 
$\{ e\in\Edge(\mygraph) | \gotv \in \partial e \}$ of the outgoing edges:
\[e_1(\gotv), \ldots, e_{\val(\gotv)}(\gotv).\]

For convenience, let us number the vertices in the graph $\mygraph$: $\gotv_1, \ldots,
\gotv_l$, where $l=\#\Vertex(\mygraph)$ is the number of vertices of the graph $\mygraph$.
We will now glue together the stable map $(C; \underline{x}; f)$ corresponding to a point 
$((C_1, \underline{x_1}), \ldots, (C_l, \underline{x_l}))$ 
in the product $\cM_\mygraph=\prod_{i=1}^l \dmsp[\deg(\gotv_i)]$.
 The curve $C$
is the union 
\[ C = \bigcup_{i=1}^l C_i \cup \bigcup_{e\in\Edge(\mygraph)} C_e, \]
where the different parts are glued together along ordinary nodes according to the
following rule:
\begin{itemize}
\item[ ] If $e\in\Edge(\mygraph)$, $\gotv_i \in \partial e$ then
$e=e_j(\gotv_i)$ for some $j$, and we will glue the $j$--th marked point $x_{i,j}$ of
$C_i$ to $0$ of $C_e$, if $\gotv_i=\gotv_1(e)$, or to $\infty$ of $C_e$ otherwise.
\end{itemize}
The ordered set of marked points $\underline{x}$ is constituted of the ``unused''
marked points of the curves $(C_i, \underline{x_i})$, \ie. the marked points at
which we have not glued curves. The function $f:C \longrightarrow \vt$ is defined
as follows:
\begin{align*}
f|_{C_e} &:= f_e \quad \text{for each edge $e$,}\\
f|_{C_i} &:= V_{\sigma(\gotv_i)} \quad \text{a constant function for each $\gotv_i$.}
\end{align*}

Let us explain this construction in plain English. Remember first, that the moduli
spaces of stable maps to a point are isomorphic to Deligne--Mumford spaces:
\[ \cM^0_{0,m}(pt.) \cong \dmsp[m]. \]
The points in $\cM_\mygraph$ therefore encode the part of the fixed stable map that is send
to fixed points in $\vt$. The parts of the fixed map that are send to 1--dimensional
invariant subspaces are in fact rigid modulo parametrization, and their ``topology'' is
encoded in the graph $\mygraph$.

Let us make a remark to the above construction for vertices $\gotv_i$ with $\deg (\gotv_i)<3$:
these curves $C_i$ are points that are glued to other points, \ie. they do not contribute
irreducible components to the constructed curve $C$. Also, if $\deg (\gotv_i)=2$ and $\val
(\gotv_i)=1$ (otherwise we must have $\val (\gotv_i)=0$ for $\deg(\gotv_i)<3$!), the remaining
marked point of $C_i$ is identified with $0$ or $\infty$ of $C_e$, respectively.

The proof of the lemma is now straightforward. The family $\pi: \cC_\mygraph
\longrightarrow \cM_\mygraph$ is the space constructed above --- essentially
it is the product of the universal curves over the Deligne--Mumford spaces
modulo the constant curves corresponding to the edges of the graph.

Since the image of a fixed stable map is rigid, the constructed family maps onto
a fixed point component of $\modsp$. Also, all fixed point components can be
obtained this way.
\end{proof}

The following notation will be useful in the sequel. A {\em flag} of a graph $\mygraph$
is a the pair of a vertex and an outgoing edge \ie. the set of flags
is
\[
\F(\mygraph):= \left\{ (\gotv, e)\in \Vertex(\mygraph)\times\Edge(\mygraph) | \gotv 
\in \partial e \right\}.
\]
For a graph type $\gtype$ or a topological graph type $\gtop$, flags are defined
the same way.The labeling of the vertices 
$\sigma: \Vertex(\mygraph)\longrightarrow \Sigma^{(d)}$ by $d$--cones 
induces a corresponding labeling of flags by 
\[ \sigma((\gotv, \star)):= \sigma(\gotv) \qquad \text{for $(\gotv,\star)\in {\mathcal F}$}.\] 
We will also use the projections of flags to vertices and edges which 
we will denote by 
\begin{align*} 
\gotv(F)&= \gotv \quad \text{for $(\gotv,\star)=F\in\F$ and $\gotv\in\Vertex$}\\ 
e(F)&= e \quad \text{for $(\star, e)=F\in\F$ and $e\in\Edge$.}\\ 
\end{align*} 
Finally we define the following subsets of $\Vertex$ and $\Edge$:
\begin{align*}
\Vertex_{t,s}&:=\left\{\gotv \in \Vertex | \val(\gotv)=t, \deg(\gotv)=t+s \right\}\\
\Vertex_t&:=\left\{\gotv\in\Vertex | \deg(\gotv)=t \right\}\\
\F_{t,s}&:=\left\{(\gotv,\star) \in \F | \gotv\in\Vertex_{s,t} \right\}\\
\F_t&:= \left\{(\gotv, \star) \in \F | \gotv\in\Vertex_t \right\}.
\end{align*}

\begin{rema}
We do not include marked points as flags, in contrary to 
Kontsevich \cite{kon95} and Graber and Pandharipande \cite{gp97}. However,
in the definition of the subsets $\Vertex_{t,s}$ \etc. the number of marked
points at a vertex $\gotv$ does enter via the function $\deg$.
\end{rema}
 
\begin{exem} 
Let us describe one example in great detail to familiarize with 
the notions defined so far. We will look at the two dimensional 
toric variety that is given by the following fan in $\Z^2$, $e_1$ 
and $e_2$ being a $\Z$--base: 
\begin{gather*} 
v_1=e_1, \quad v_2=e_2, \quad v_3=-e_1+e_2, \quad v_4=-e_2\\ 
{\mathcal P} = \left\{ \{v_1, v_3\}, \{v_2,v_4\}\right\}. 
\end{gather*} 
The fan $\Sigma$ having the $1$--skeleton $\Sigma^{(1)}=\{v_1,\ldots, v_4\}$ and 
the set of primitive collections $\mathcal P$ is shown in figure \ref{fig:polyf1}, as well 
as its polyhedron corresponding to the strictly convex upper support function 
$\varphi=c_1(\vt)$. The toric variety $\vt$ constructed from $\Sigma$ is the 
Hirzebruch surface $F_1=\P(\O_{\P^1}(1)\oplus 1)\cong \widetilde{\P^2}$, which is 
isomorphic to $\P^2$ blown up at one point. 
 
\mycomment{fig:polyf1} 
\begin{figure}[ht] 
\begin{center} 
\unitlength1mm 
\begin{picture}(110,35) 
\put(0,0){\epsfig{figure=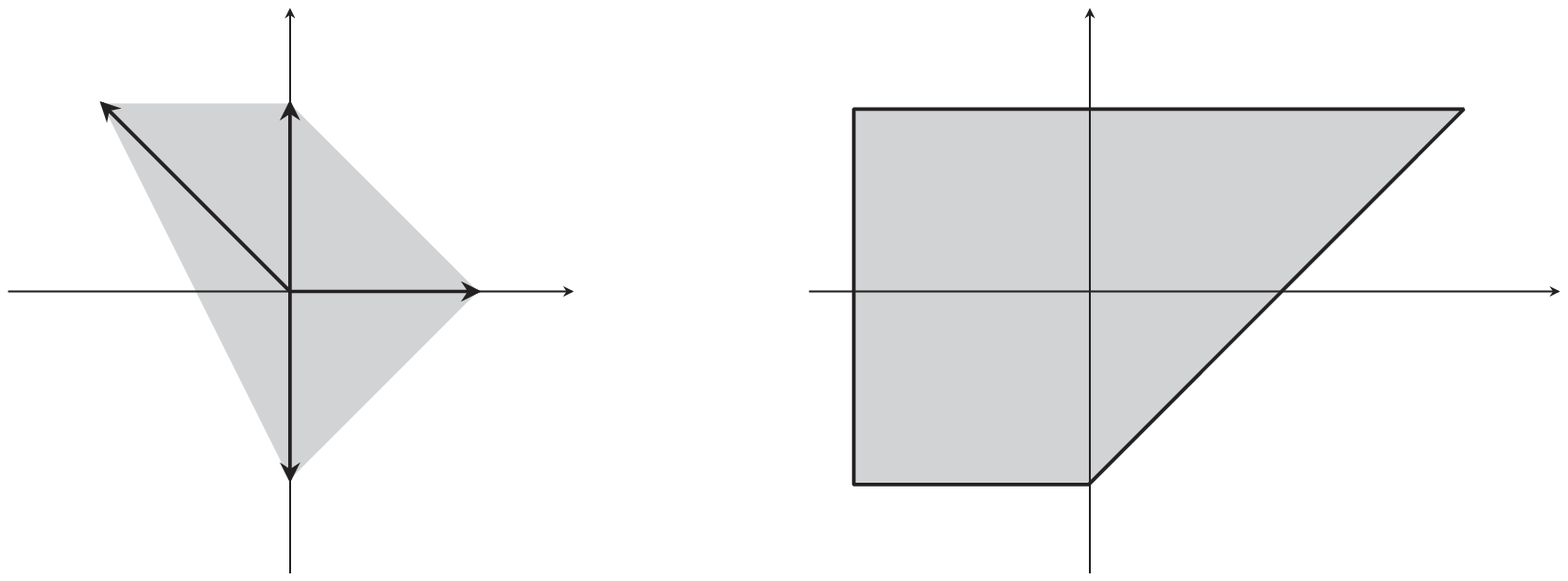,height=35mm}} 
\put(29,18){$v_1$} 
\put(18,30){$v_2$} 
\put(4,31){$v_3$} 
\put(18,3){$v_4$} 
\put(20,20){$\sigma^2_1$} 
\put(12,24){$\sigma^2_2$} 
\put(13.5,14.25){$\sigma^2_3$} 
\put(20,12.5){$\sigma^2_4$} 
\put(70,20){$(0)$} 
\put(49,20){$v_1$} 
\put(57.7,2){$v_2$} 
\put(76,10){$v_3$} 
\put(75,30){$v_4$} 
\put(48,3){$\sigma^2_1$} 
\put(69,3){$\sigma^2_2$} 
\put(92,27.5){$\sigma^2_3$} 
\put(48,27.5){$\sigma^2_4$} 
\end{picture} 
\caption{The fan of $\P(\O_{\P^1}(1)\oplus 1)$, and its polyhedron for $\varphi=c_1$.}\label{fig:polyf1} 
\end{center} 
\end{figure} 
Before we give a graph corresponding to a fixed point in $\modsp$ 
of this toric variety, let us analyze the homology and cohomology 
in degree two of $\vt$. We have seen above, that (integral) 
degree--$2$ cohomology classes are given by $\Sigma$--piecewise 
linear functions, factored out by linear functions $\psi\in 
M=\Hom(N,\Z)$. A function $\varphi\in PL(\Sigma)$ is given by its 
values on the $1$--skeleton, an element $\psi\in M$ by its values 
on $e_1, e_2$. Hence for a $\varphi$ representing an equivalence 
class $[\varphi]\in PL(\Sigma)/M$ we can assume 
\[ \varphi(v_1)=\varphi(v_2)=0,\qquad \varphi(v_3),\varphi(v_4)\in \Z.\] 
Such a class $[\varphi]$ is in the K\"ahler cone if it satisfies 
\begin{gather*} \varphi(v_1)+\varphi(v_3)>\varphi(v_1+v_3) \quad \text{and} \quad 
\varphi(v_2)+\varphi(v_4)>\varphi(v_2+v_4)\\ 
\intertext{that is, with the choices above,} 
\varphi(v_3) > 0 \quad \text{and}\quad \varphi(v_4)>0. 
\end{gather*} 
Note, that this implies in particular, that the first Chern class $c_1(\vt)$ of $\vt$ is indeed 
a K\"ahler class. 
 
For the degree--$2$ homology of $\vt$, notice that the $\Z$--module 
\[ R(\Sigma) = \left\{ (\lambda_1, \ldots, \lambda_n)\,|\, \lambda_1 v_1+\ldots+\lambda_n v_n=0\right\} \] 
is generated by the elements corresponding to the equations 
\begin{gather*} 
\quad v_2+v_4=0 \quad \text{and} \quad v_1+v_3+v_4=0\\ 
\intertext{that is by the elements} 
\lambda^1:=(0,1,0,1) \quad \text{and}\quad \lambda^2:=(1,0,1,1). 
\end{gather*} 
 
To find out the homology classes of the four one dimensional $T_N$--invariant subvarieties 
$V_{\langle v_1\rangle},\ldots, V_{\langle v_4\rangle}$, the Poincar\'e dual cohomology classes 
$[\varphi_1], \ldots, [\varphi_4]$ of which are given by 
\[ \tilde\varphi_i(v_j)=\delta_{ij}= \left\{ \begin{array}{ll} 
                        1 \quad&\text{if i=j}\\[3pt] 
                        0 \quad&\text{otherwise.} 
                         \end{array} \right. 
\] 
Hence, again by Poincar\'e duality, we get 
\[ 
[V_{\langle v_1\rangle}] = \lambda^2, \quad [V_{\langle v_2\rangle}] = \lambda^1, \quad 
[V_{\langle v_3\rangle}] = \lambda^2,\quad [V_{\langle 
v_4\rangle}] = \lambda^1 + \lambda^2.
\] 
Therefore, any 
$\modsp$--graph $\mygraph$ has to ``live'' on the decorated 
$1$--skeleton $\Upsilon_\Sigma$ of the moment polytope $\Delta_\varphi$
shown in figure \ref{fig:rasterf1} 
in the sense that there is a map $f:\mygraph \longrightarrow
\Upsilon_\Sigma$ of one--dimensional CW--complexes such that the 
decorations $\sigma:\Vertex(\mygraph)\longrightarrow \Sigma^{(d)}$ of 
the vertices of $\mygraph$ with fixed points in $\vt$ are induced 
from the decorations of the vertices of $\Upsilon_\Sigma$. Figure 
\ref{fig:graphf1} shows two $\modsp[0]$--graphs for the homology 
class $A=2\lambda^2+\lambda^1$. Note that there are other possible 
graphs for this class. 
 
\mycomment{fig:rasterf1 fig:graphf1} 
\begin{figure}[ht] 
\begin{center} 
\unitlength1mm 
\begin{minipage}[t]{0.45\linewidth} 
\begin{center} 
\begin{picture}(50,35) 
\put(10,2.5){\epsfig{figure=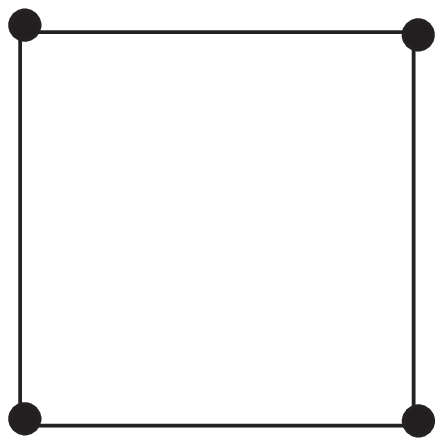,height=30mm}} 
\put(7,15){$\lambda^2$} 
\put(39.5,15){$\lambda^2$} 
\put(21,31.5){$\lambda^1+\lambda^2$} 
\put(24,4){$\lambda^1$} 
\put(6.5,0){$\sigma^2_1$} 
\put(40,0){$\sigma^2_2$} 
\put(40,32.5){$\sigma^2_3$} 
\put(6.5,32.5){$\sigma^2_4$} 
\end{picture} 
\caption{The $1$--skeleton $\Upsilon_\Sigma$ of the moment polytope $\Delta_varphi$.}\label{fig:rasterf1} 
\end{center} 
\end{minipage}\hspace{1cm} 
\begin{minipage}[t]{0.45\linewidth} 
\begin{center} 
\begin{picture}(45,35) 
\put(5,2.5){\epsfig{figure=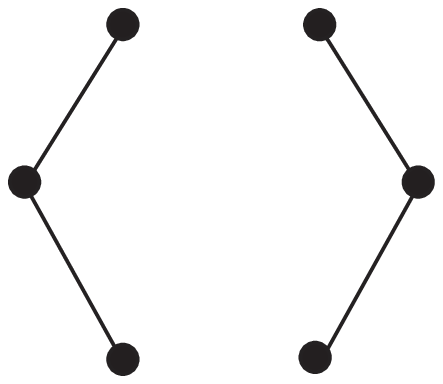,height=30mm}} 
\put(9,32){$\sigma^2_4$} \put(33,32){$\sigma^2_4$} 
\put(8,24.5){$2$} \put(36,24.5){$1$} 
\put(0.5,17.5){$\sigma^2_1$} \put(41,17.5){$\sigma^2_3$} 
\put(12,11.5){$1$} \put(32,11.5){$1$} 
\put(8.5,2){$\sigma^2_2$} \put(32,2){$\sigma^2_2$} 
\end{picture} 
\caption{Two different graph types for $A=2\lambda^2+\lambda^1$.}\label{fig:graphf1} 
\end{center} 
\end{minipage} 
\end{center} 
\end{figure} 
\end{exem} 
 
\subsection{Automorphisms of fixed point components} 
 
For a family $\pi: \cC_\mygraph \longrightarrow \cM_\mygraph$ of 
$T_N$--fixed stable maps to $\vt$ as constructed above, there are two
different sources of $\pi$--equivariant automorphisms we have to consider:
automorphisms of the family $\pi: \cC_\mygraph \longrightarrow \cM_\mygraph$
itself, and automorphisms of this family as substack of $\modsp$.

The first are given by automorphisms of the $\modsp$--graph $\mygraph$
(with its decorations), since by considering products of Deligne--Mumford
spaces we have ordered the nodes. All other possible sources of automorphisms
have already been moded out by taking choices in the proof of Proposition
\ref{prop:isfixedptcomp}.

The extra automorphisms stemming from the map to the stack $\ucmodsp
\longrightarrow \modsp$ is a cyclic permutation of the branches of the
degree--$d(e)$ maps $f_e$. This $\Z_{d(e)}$--action is trivial on $\cM_\mygraph$,
but not on its deformations, underlining the ``orbifold character'' of
the moduli stack $\modsp$.

\begin{lemm}
The automorphism group 
${\bf A}_\mygraph$ of $\pi:\cC_\mygraph \longrightarrow\cM_\mygraph$ 
fits into the following exact sequence of groups 
\[ 1\longrightarrow \prod_{e\in\Edge} \Z_{d(e)} \longrightarrow {\bf A}_\mygraph 
\longrightarrow \Aut(\mygraph) 
\longrightarrow 1,\] where $\Aut(\mygraph)$ acts naturally on $\prod_e 
\Z_{d(e)}$, ${\bf A}_\mygraph$ being the semi--direct product. The 
induced map 
\[ \gamma/{\bf A}_\mygraph: \M_\mygraph/{\bf A}_\mygraph \mapright \modsp \] 
is a closed immersion of Deligne--Mumford stacks. Furthermore, the image is a component 
of the $T_N$--fixed point stack of $\modsp$. 
\proofend
\end{lemm}
 
\subsection{Weights on fixed point components} 
 
At the end of this section, we will compute the weight of the $T_N$--action 
on the irreducible $T_N$--invariant divisors $V_\tau$, $\tau \in \Sigma^{(d-1)}$, and subsequently 
we will derive the weight of the action on a non--constant map 
\[ f: \P^1 \mapright V_\tau\subset\vt, \qquad \tau\in\Sigma^{(d-1)}\] 
represented by an edge $e$ in a $\modsp$--graph $\mygraph$. Let 
$\sigma_1, \sigma_2 \in \Sigma^{(d)}$ be two $d$--cones in 
$\Sigma$ that have a common $(d-1)$--face $\tau\in\Sigma^{(d-1)}$. 
Notice that $V_\tau$ is the closure of a one--dimensional orbit of 
the $T_N$ action, compactified with the two fixed points of this 
action given by the two $d$--cones $\sigma_1$ and $\sigma_2$. So, 
the $T_N$ action reduces to a $\C^*$--action on $V_\tau$, that is 
to the action of a subtorus $\C^*\cong T_\tau \cong T_N$ of $T_N$. 
The torus $T_\tau$ is the image of the map 
\[ \beta_\tau: (\C^*)^n \mapright T_N \mapright T_\tau\] 
given by the quotient map $T_N=(\C^*)^n/\DD(\Sigma)$ followed by restriction to $V_\tau$. Thus, we 
can write elements of $T_\tau$ as equivalence classes of elements in $(\C^*)^n$ by the map $\beta_\tau$. 
 
\begin{lemm}\mylabel{lem:nodeweight} 
Let $\sigma_1,\sigma_2,\tau\in\Sigma$ as above. 
Let $v_{i_1},\ldots,v_{i_{d-1}}$ be the generators of the common face $\tau=\sigma_1\cap\sigma_2$, such that 
\begin{align*} 
\sigma_1&=\langle v_{i_1}, \ldots, v_{i_{d-1}}, v_{l_1(\tau)} \rangle\\ 
\sigma_2&=\langle v_{i_1}, \ldots, v_{i_{d-1}}, v_{l_2(\tau)} \rangle. 
\end{align*} 
 
Let $\we_1, \ldots, \we_n$ be the weights of a diagonal action 
of $(\C^*)^n$ on $\C^n$ with respect to the standard basis.
The induced $\C^*$--action on the subvariety $V_\tau$, 
has weight $\we^{\sigma_1}_{\sigma_2}$ at the point $V_{\sigma_1}$: 
\[
\we^{\sigma_1}_{\sigma_2} := \sum_{j=1}^n \langle v_j, u_d\rangle \we_j, 
\] 
where $u_1,\ldots, u_d$ is the basis of $M=\Hom(N,\Z)$ dual to $v_{i_1},\ldots, v_{i_{d-1}}, v_{l_1(\tau)}$. 
\end{lemm}

\begin{proof} 
The $d$--dimensional cone $\sigma_1$ gives a local chart $U_{\sigma_1}$ of our toric variety $\vt$, 
and the coordinates on $U_{\sigma_1}$ are given 
by (\cf. proposition \ref{thm:homcoords}): 
\[ x_1^{\sigma_1}= \prod_j z_j^{\langle v_j, u_1\rangle}, \quad \ldots\quad, 
x_d^{\sigma_1}=\prod_j z_j^{\langle v_j, u_d\rangle},\] 
 
The $1$--dimensional submanifold corresponding to $\tau$ is given 
by the equations $z_{i_1}=\ldots=z_{i_{d-1}}=0$. In the coordinates 
of $U_{\sigma_1}$, these equations are equivalent to 
\[ x_1^{\sigma_1}=\ldots=x_{d-1}^{\sigma_1}=0.\] 
 
Hence we have to look at the $(\C^*)^d$--action on the $d^{th}$ co--ordinate. Thus the 
action of $(t_1, \ldots, t_n)\in(\C^*)^n$ on $V_\tau$ is given by 
\begin{align*} 
(t_1,\ldots, t_n) \cdot x_d^{\sigma_1}& = \prod_j (t_jz_j)^{\langle v_j, u_d\rangle}\\ 
&= \left(\prod_j t_j^{\langle v_j, u_d\rangle}\right) \cdot 
\left(\prod_j z_j^{\langle v_j, u_d\rangle}\right)\\ 
&= t^{\we^{\sigma_1}_{\sigma_2}} x_d, 
\end{align*} 
using multi--index notation in the last line. Hence the weight of the action on $V_\tau$ is indeed 
$\sum_j \langle v_j, u_d\rangle \we_j$ in the chart $U_{\sigma_1}$. 
\end{proof} 
 
The lemma above gives in particular the $T_N$--action on a the component of a fixed 
stable curve that is mapped to $V_\tau$: 
 
\begin{coro} 
Let $e\in\Edge(\mygraph)$ be an edge of the $\modsp$--graph 
$\mygraph$, and $\gotv_1,\gotv_2\in\partial e$ be the vertices at 
its two ends. Let $\sigma_i=\sigma(\gotv_i)$ be the $d$--cones of 
the vertices $\gotv_i$, and $\tau(e)=\sigma_1\cap\sigma_2$ its 
common $(d-1)$--face, that are generated by 
\begin{align*} 
\sigma_1=\sigma(\gotv_1)&=\langle v_{i_1}, \ldots, v_{i_{d-1}}, v_{l_1(e)} \rangle\\ 
\sigma_2=\sigma(\gotv_2)&=\langle v_{i_1}, \ldots, v_{i_{d-1}}, v_{l_2(e)} \rangle. 
\end{align*} 
For a $T_N$--fixed stable map $(C; x_1,\ldots, x_m; f)\in 
\M_\mygraph\subset\modsp$, let $C_e$ be the irreducible component 
of $C$ corresponding to the edge $e$. Let $F:=(\gotv_1, e)\in\F$ 
be the flag of the edge $e$ at the vertex $\gotv_1$. At the point 
$p_F:=f^{-1}(V_{\sigma(\gotv_1)}) \cap C_e$, the pull back to 
$C_e$ by $f$ of the $T_N$--action on $V_{\tau(e)}$ has the weight 
$\we_F$ at $p_F$: 
\[
\we_F := \frac{1}{d_e} \sum_{j=1}^n \langle v_j, u_d\rangle \we_j, 
\]
where $d_e$ is the multiplicity of the component $C_e$, and $u_1,\ldots, u_d$ is the basis 
of $M=\Hom(N,\Z)$ dual to $v_{i_1},\ldots, v_{i_{d-1}}, v_{l_1(e)}$. 
\end{coro} 
 
\begin{proof} 
The action of $T_N$ on $C_e$ is just the pull back by $f$ of the action on $V_\tau$. Since $f$ 
has multiplicity $d_e$, the formula follows immediately from lemma \ref{lem:nodeweight}. 
\end{proof} 
 
We will introduce some further notation, grouping together certain 
weights on the one--dimensional $T_N$--invariant subvariety of 
$\vt$, or more general on a $\modsp$--graph $\mygraph$. First of 
all, we will write $\sigma_1\face\sigma_2$ for the property of 
$\sigma_1$ and $\sigma_2$ having a common $(d-1)$--dimensional 
proper face: 
\[ \sigma_1\face\sigma_2 \Longleftrightarrow \sigma_1,\sigma_2\in\Sigma^{(d)} \medspace \text{and} \medspace 
\sigma_1\cap\sigma_2 \in\Sigma^{(d-1)}. \] 
 
The {\em total weight of a $d$--dimensional cone $\sigma$} is defined 
to be 
\[ \totwe{\sigma} := \prod_{\gamma\face\sigma} \we_\gamma^\sigma.\] 
Note that $\totwe{\sigma}$ is in fact a polynomial in the generators of $\gott^n$: 
$\totwe{\sigma}\in \Z[\we_1,\ldots,\we_n]$.

\section{The virtual normal bundle for toric varieties}
\mylabel{sec:virtnormbdl}

In this section we analyze Graber and Pandharipande's virtual
normal bundle to the fixed point components of the moduli space of
stable maps for the natural $(\C^*)^n$--action on a toric variety,
hence generalizing Graber and Pandharipande's example for
projective space $\C\P^n$ (\cite{gp97}), and we will derive our
main result. Contrary to their calculations for $\C\P^n$, however,
we will restrict ourselves here to genus zero stable maps.

So let $\vt$ be a smooth projective complex variety.
Remember from section \ref{subsec:obstrcplxGW}, that for the cohomology sheaves of
the dual natural perfect obstruction
theory $E_\bullet$ for our moduli stack $\modsp$ of stable maps
\mycomment{eq:homolobcom}
\begin{equation}\label{eq:homolobcom}
0\mapright \cT^0\mapright E_{0} \MRE{d}
E_{1}\mapright \cT^1\mapright 0,
\end{equation}
the sheaves $\cT^i$ are given by:
\[ \cT^i = \underline\Ext^i_\pi([f^*\Omega^1_X \longrightarrow
\Omega^1_{\ucmodsp/\modsp}(D)], \O_{\ucmodsp}),
 \quad i=0,1. \]

As before we will sometimes write $\cM$ and $\cC$ for the
moduli space $\modsp$ and its universal curve $\ucmodsp$, respectively,
if no confusion can arise.

\begin{lemm}
The sheaves $\cT^i$ fit into the following exact sequence:
\mycomment{eq:gobcomlong}
\begin{equation}\label{eq:gobcomlong}
\begin{split}
0&\mapright \underline\Hom_\pi(\Omega^1_{\cC/\cM}(D),
\O_{\cC}) \mapright
\underline\Hom_\pi(f^*\Omega^1_X, \O_{\cC}) \mapright \cT^0
\mapright\\
&\mapright \underline\Ext^1_\pi(\Omega^1_{\cC/\cM}(D),
\O_{\cC}) \mapright
\underline\Ext_\pi^1(f^*\Omega^1_X, \O_{\cC})
\mapright \cT^1\mapright 0.\\
\end{split}
\end{equation}
\end{lemm}

\begin{proof}
Let $\cK^\bullet$ be the complex
$\cK^\bullet=[f^*\Omega^1_X\longrightarrow \Omega^1_{\ucmodsp/\modsp}(D)]$
indexed at $-1$ and $0$. It fits into the following short exact sequence:
\[
0 \mapright \Omega^1_{\cC/\cM}(D) \mapright \cK^\bullet
\mapright f^*\Omega^1_X[1] \mapright 0.\]
The corresponding long exact sequence of higher direct image sheaves corresponding
to $\pi: \ucmodsp \longrightarrow \modsp$ is then (\ref{eq:gobcomlong} where
exactness on the left, \ie. the injectivity of the map
${\underline\Hom}_\pi(\Omega^1_{\cC/\cM}(D), \O_{\cC})
\longrightarrow {\underline\Hom}_\pi(f^*\Omega^1_X, \O_{\cC})$
induced by the natural map $f^*\Omega^1_X\longrightarrow
\Omega^1_{\cC/\cM}\longrightarrow \Omega^1_{\cC/\cM}(D)$ is
equivalent to the stability of the map $f: \ucmodsp
\longrightarrow X$ (\cf. lemma \ref{lem:stableinj} and the remark
following the lemma).  Exactness on the right of the long
exact sequence follows from the fact of the fibers of $\ucmodsp
\longrightarrow \modsp$ being curves.
\end{proof}

Now, let $\cM_\mygraph$ be a fixed point component in the moduli
stack of stable maps $\modsp$, and $\pi_\mygraph:\cC_\mygraph
\longrightarrow \cM_\mygraph$ its universal curve. By lemma
\ref{lem:restrobth}, we know that the restriction of the long
exact sequence (\ref{eq:gobcomlong}) to the fixed point component
$\cM_\mygraph$ becomes: 
\mycomment{eq:loccomlong}
\begin{equation}\label{eq:loccomlong}
\begin{split}
0&\longrightarrow
\underline\Hom_{\pi_\mygraph}(\Omega^1_{\cC_\mygraph/\cM_\mygraph}(D),
\O_{\cC_\mygraph}) \longrightarrow
\underline\Hom_{\pi_\mygraph}(f_\mygraph^*\Omega^1_X,
\O_{\cC_\mygraph}) \longrightarrow \cT^0|_{\cM_\mygraph} \longrightarrow\\
&\longrightarrow
\underline\Ext^1_{\pi_\mygraph}(\Omega^1_{\cC_\mygraph/\cM_\mygraph}(D_\mygraph),
\O_{\cC_\mygraph}) \longrightarrow
\underline\Ext^1_{\pi_\mygraph}(f_\mygraph^*\Omega^1_X,
\O_{\cC_\mygraph}) \longrightarrow \cT^1|_{\cM_\mygraph}\longrightarrow 0.\\
\end{split}
\end{equation}
In particular, if we restrict to a single stable map $(C;
\underline{x}; f)\in \cM_\mygraph$, we get: 
\mycomment{eq:obcomlong}
\begin{equation}\label{eq:obcomlong}
\begin{split}
0&\longrightarrow \Ext^0(\Omega_C(\underline{x}),\O_C)\longrightarrow 
H^0(C, f^*TX)\longrightarrow {\mathcal T}^0|_{\{C\}}
\longrightarrow\\
&\longrightarrow \Ext^1(\Omega_C(\underline{x}),\O_C)\longrightarrow 
H^1(C, f^*TX)\longrightarrow {\mathcal T}^1|_{\{C\}}
\longrightarrow 0.\\
\end{split}
\end{equation}
To determine the virtual  normal bundle for the torus action on
our toric variety, we will have to fulfill the moving parts of
$\cT^i_\mygraph=\cT^i|_{\cM_\mygraph}$, hence the moving parts in
the sequence (\ref{eq:loccomlong}). We will call the $i^{th}$ term of
this long exact sequence by $B_i$, and the moving part (with
respect to the induced torus action) by $B_i^\move$. Remember,
that the virtual normal bundle $N_\mygraph^\virt$ is the moving
part of the induced complex $E_{\bullet,\mygraph}$, that is
\[ e^{T_N}(N_\mygraph^\virt)=e^{T_N}(E_{0, \mygraph}^\move -E_{1,\mygraph}^\move)
=e^{T_N}({\mathcal T}_\mygraph^{0,\move}-{\mathcal T}_\mygraph^{1,\move}),
\]
where the second equation holds because of the exact sequence (\ref{eq:homolobcom}).
Now, applying the long exact sequence (\ref{eq:loccomlong}) we obtain for the
equivariant Euler class of the
virtual normal class $N^{vir}$ the following formula:
\mycomment{eq:eczerl}
\begin{equation}\label{eq:eczerl}
e^{T_N}(N_\mygraph^{\text{virt}})=\frac{e^{T_N}(B_2^\move)e^{T_N}(B_4^\move)}{e^{T_N}(B_1^\move)e^{T_N}(B_5^\move)}
\in H^*_{T_N}(X, \Q).
\end{equation}
The notation is indeed correct, since the $B_i^\move$,
$i=1,2,4,5$, are vector bundles on $\M_\mygraph$! This does not
apply in general to the fixed parts of these sheaves, or even to
the sheaves ${\mathcal T}^i_\mygraph$. So actually, at least for the
moving parts, we look at a long exact sequence of the kind of
(\ref{eq:obcomlong}).

In the following, we will calculate the contributions of the four bundles
to the equivariant Euler class of the virtual normal bundle.

\subsection{Computation of the equivariant Euler class of $B_1^\move$}

This bundle almost never contributes to the virtual normal bundle --- indeed
we will show the following lemma:

\begin{lemm}
The $T_N$--equivariant Euler class of $B_1^\move$ is given by
\mycomment{eq:ec1}
\begin{equation}\label{eq:ec1}
e^{T_N}(B_1^\move)= \prod_{F\in\F_1} \we_F.
\end{equation}
\end{lemm}

\begin{figure}[ht]
\begin{minipage}[t]{.22\linewidth}
\begin{center}
\epsfig{figure=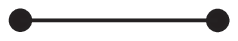,width=\linewidth}
\caption{A graph with two vertices in $\V_{1,0}$.}\label{fig:vnbb1f1}
\end{center}
\end{minipage}\hfill
\begin{minipage}[t]{.22\linewidth}
\begin{center}
\epsfig{figure=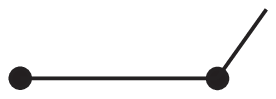,width=\linewidth}
\caption{A graph with the left vertex in $\V_{1,0}$ and the right vertex in $\V_{1,1}$.}\label{fig:vnbb1f2}
\end{center}
\end{minipage}
\hfill
\begin{minipage}[t]{.44\linewidth}
\begin{center}
\epsfig{figure=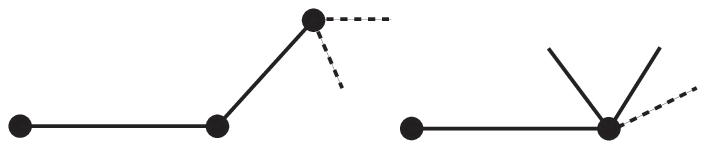,width=\linewidth}
\caption{Two examples for parts of graphs that can occur in $\M_{0,m}(X,A)$ for $m\ge3$.}\label{fig:vnbb1f3}
\end{center}
\end{minipage}
\end{figure}

\begin{proof}
The bundle $B_1=\Ext^0(\Omega_C(D), \O_C)=\Aut_\infty(C)$
parameterizes infinitesimal automorphisms of the pointed domain.
The induced $T_N$--action on $\Aut(C)$ is obviously trivial on all
automorphism $\varphi$ of $C$ that restrict to the identity
$\varphi_{|C_e}=\id_{C_e}$ on all irreducible components $C_e$
corresponding to edges $e\in\Edge(\mygraph)$ in the graph
$\mygraph$. Thus, the moving part of $\Aut_\infty(C)$ splits into
\[ \Aut_\infty^\move(C) = \bigoplus_{e\in\Edge(\mygraph)} \Aut_\infty^\move(C_e).\]
Note in particular, that the bundle $\Aut_\infty^\move(C)$ is
topologically trivial on $\M_\mygraph$ since it only depends on
the irreducible components that are not mapped to a point, \ie.
that are rigid in $\M_\mygraph$.

The $\Aut_\infty^\move(C_e)$ obviously depend on the special points
of the irreducible component $C_e$, \ie. on the way it is ``glued''
at the nodes corresponding to the vertices $\gotv_1(e)$ and $\gotv_2(e)$
of the edge $e$. Since we only look at moduli stacks of stable maps with
at least three marked points, we can exclude the two special cases 
(\cf. figures \ref{fig:vnbb1f1} and \ref{fig:vnbb1f2}) where both
vertices are in $\Vertex_{1,0}$, or where one is in $\Vertex_{1,0}$ and
the other in $\Vertex_{1,1}$. We are left with two cases:
\begin{description}
\item[Case 1] {\bf One edge in $\Vertex_{1,0}$, the other in $\Vertex_{2,0}
\cup \Vertex_{\ge 3}$} --- Without loss of generality we can assume that
$\gotv_1(e)\in \Vertex_{1,0}$. In this case, $C_e$ corresponds to a 
non--contracted $\P^1$ attached to
another non--contracted (or, in the $\Vertex_{\ge 3}$--case, contracted)
component (see figure \ref{fig:vnbb1f3}). Therefore we have to
look at M\"obius transformations that fix one point, infinity say:
\[
[x_1:x_2] \mapsto [ax_1+bx_2:x_2].
\]
Let $F\in\F_1$ be the flag corresponding to $\gotv_1$. We have seen above 
that the induced $T_N$--action
on $C_e$ is given by (using again multi--index notation):
\[ t\cdot[x_1:x_2]=[t^{\we_F} x_1:x_2],\]
since the co--ordinate $x_1$ corresponds to the chart of the flag $F$ 
(while $x_2$ corresponds to the chart
around infinity, \ie. at the vertex $\gotv_2$. To determine the 
$T_N$--action on the group $\Aut_\infty(C_e)$ of infinitesimal automorphisms
of $C_e$, we have to compute:
\begin{align*}
t\cdot(a,b)\cdot t^{-1}\cdot [x_1:x_2]&=t\cdot(a,b)\cdot[t^{-\we_F}x_1:x_2]\\
&=t\cdot[at^{-\we_F}x_1+bx_2:x_2]\\
&=[ax_1+t^{\we_F}bx_2:x_2]\\
&=(a,t^{\we_F}b)\cdot[x_1:x_2],
\end{align*}
hence the $T_N$--action on an automorphism of $C_e$ is given by:
\[
t\cdot(a,b) = (a, t^{\we_F}b).
\]
The moving part of $\Aut_\infty(C_e)$ is thus spanned by the second co--ordinate,
and the weight of the action there is $\we_F$ --- in this case we have:
\[ e^{T_N}(\Aut_\infty(C_e)^\move) = \we_F.\]

\item[Case 2] {\bf Neither vertex is in $\Vertex_{1,0}$} --- 
In this case, any automorphism of $C$ restricts to an automorphism on 
$C_e$ that fixes the two points corresponding to the special points of 
the vertices $\gotv_1$ and $\gotv_2$. Any such automorphism on $C_e$ 
(w.~l.~o.~g.\ we take the two points to be zero and infinity) has to look like
\[ [x_1:x_2] \DMRE{} [a x_1:x_2],\]
where $a\neq 0$ is a non-negative integer. With the same analysis as above 
of the $T_N$--action on such automorphism $a$,
we see that the $T_N$--action on $\Aut_\infty(C_e)$ is trivial, \ie.
\[ \Aut_\infty^\move(C_e)=(0).\]
\end{description}
\end{proof}

\vspace{0.5cm} \mycomment{fig:vnbb1f1-3}

\subsection{The equivariant Euler class of $B_4^\move$}

Here we are looking at the bundle $B_4=\Ext^1(\Omega_C(D),\O_C)=\operatorname{Def}(C)$ of deformations
of the pointed domain, that is deformations that vary some of the special points (varying the isomorphism
class of the curve) or that smooth some
double points.
Again, deformations of contracted components have obviously weight zero,
since the $T_N$--action on these components is trivial, so deformations coming from varying special
points do not contribute to $B_4^\move$.

The other deformations of $C$ come from smoothing nodes of $C$ which join a non--contracted
component and a contracted or non--contracted component. Such a smoothing corresponds to
choosing an element of the tangent bundle at the double point. So let ${\mathcal L}_F$ be the
universal cotangent line (\cf. section \ref{sec:lbdms}) at the double point corresponding
to an $F\in\F_{\ge3}\cup\F_{2,0}$, and
write $e_F=e({\mathcal L}_F)=c_1({\mathcal L}_F)$ for the usual Euler class of this line bundle.

If we look at the smoothing of a double point between a contracted and a non--contracted
component, \ie. if $F\in\F_{\ge3}$, let $F=(\gotv,e)$. In this case, the bundle $\cL_F$ is
the pull back to $\cM_\mygraph$ of the corresponding cotangent line on the Deligne-Mumford space 
$\dmsp[\deg(\gotv)]$ of $\gotv$; on the other components of $\cM_\mygraph$, the bundle $\cL_F$
is trivial. The $T_N$-action on $\cM_\mygraph$ is trivial, and we have seen above that the
$T_N$--action on the bundle $\cL_F$ has weight $\we_F$. Hence in this case:
\[ e^{T_N}({\mathcal L}_F^*) = \we_F - e_F.\]

In the second case, when we look at a vertex $\gotv\in\Vertex_{2,0}$ joining two non--contracted
components, we analogously obtain for the equivariant Euler class of the tangent line at this
node:
\[ e^{T_N}({\mathcal L}_F^*) = \we_{F_1} + \we_{F_2} - e_{F_1} - e_{F_2},\]
where $F_1, F_2\in \F_{2,0}$ are the two flags at $\gotv$.
However, ${\mathcal L}_{F_1}$ and ${\mathcal L}_{F_2}$ are topologically
trivial on $\M_\mygraph$ (since non--contracted components are
rigid in $\M_\mygraph$), so we have proven:

\begin{lemm} The $T_N$--equivariant Euler class of the moving part of the bundle
$B_4$ is equal to:
\mycomment{eq:ec2}
\begin{equation}\label{eq:ec2}
e^{T_N}(B_4^\move)= \prod_{F\in\F_{\ge3}} (\we_F-e_F)
\prod_{\gotv\in\V_{2,0}} (\we_{F_1(\gotv)}+\we_{F_2(\gotv)}),
\end{equation}
where $F_1(\gotv), F_2(\gotv)$ denote the two different flags at
the vertex $\gotv\in\Vertex_{2,0}$.\proofend
\end{lemm}

\subsection{The equivariant Euler class of the quotient $B_2^\move-B_5^\move$}

Like Graber and Pandharipande, we will use the following exact
sequence to calculate the contribution coming from $H^*(f^*TX)$:
\[
0\longrightarrow\O_C\longrightarrow 
\bigoplus_{\gotv\in\V}\O_{C_\gotv}\oplus\bigoplus_{e\in\E}\O_{C_e}
\longrightarrow \bigoplus_{F\in\F} \O_{x_F}\longrightarrow 0.
\]

Note, that a--priori it only makes sense to sum over $\V_{\ge3}$ in the
middle term, and over $\F_{\ge3}$ and $\V_{2,0}$ (instead of $\F$) in
the right term. However, we add the same in both terms, so they
cancel each other. By passing to the pullback under $f$ and taking
cohomology, we obtain: 
\begin{multline*}
0\rightarrow H^0(f^*T\vt)\rightarrow\bigoplus_{\gotv\in\V} H^0(C_\gotv, f^*T\vt)\oplus
\bigoplus_{e\in\E}H^0(C_e, f^*T\vt)\rightarrow\\
\rightarrow \bigoplus_{F\in\F}T_{p_{i(F)}}\vt \rightarrow H^1(f^*T\vt)\rightarrow
\bigoplus_{e\in\E}H^1(C_e, f^*T\vt)\rightarrow 0.\\
\end{multline*}

Note that since we only look at genus zero curves, $H^1(C_v, f^*T\vt)=0$. On the
other hand, $H^1(C_e, f^*T\vt)$ is not necessarily zero for a toric variety $\vt$ as it is in general not
convex. Since $f^*T\vt$ is trivial on $C_\gotv$ for $\gotv\in\Vertex$, 
$H^0(C_\gotv, f^*T\vt)=f_{p_{\sigma(\gotv)}}^*T\vt$. Hence we obtain the following formula:
\mycomment{eq:ec3a}
\begin{equation}\label{eq:ec3a}
\begin{split}
H^0(f^*T\vt)-H^1(f^*T\vt)&=+\, \bigoplus_{\gotv\in\V}T_{p_{\sigma(\gotv)}}\vt \, +\,
\bigoplus_{e\in\E} H^0(C_e, f^*T\vt)\\
&\phantom{{=}}-\,\bigoplus_{F\in\F} T_{p_{\sigma(F)}}\vt
-\,\bigoplus_{e\in\E} H^1(C_e, f^*T\vt).\\
\end{split}
\end{equation}

To compute the equivariant Euler class of $H^0(C_e, f^*T\vt)$, we
again observe, that the bundle is constant.
To determine the weights of the induced action, we look at the following
Euler sequence on $\vt$:
\[ 0 \rightarrow \O^{n-d} \rightarrow \O(Z_1)\oplus\ldots\oplus\O(Z_n)
\rightarrow T\vt \rightarrow 0. \]
Pulling back to $C_e$ and taking cohomology gives
\[ 0 \rightarrow \C^{n-d} \rightarrow H^0(C_e, \O(\cocl_1))\oplus\ldots\oplus
H^0(C_e, \O(\cocl_n))\rightarrow H^0(C_e, f^*T\vt)\rightarrow 0,\]
where the tuple $(\cocl_1, \ldots, \cocl_n)$ describes the homology
class of $f_*[C_e]$.

As in section \ref{sec:tvaction}, let $\partial e=\{\gotv_1,\gotv_2\}$ be the two
nodes at the ends of the edge $e$, $\sigma_i=\sigma(\gotv_i)\in\Sigma^{(d)}$ be the
two $d$--cones in the fan $\Sigma$ corresponding to the two nodes $\gotv_1$ and $\gotv_2$,
and let
\begin{align*}
\sigma_1&=\langle v_{i_1}, \ldots, v_{i_{d-1}}, v_{l_1(e)}\rangle\\
\sigma_2&=\langle v_{i_1}, \ldots, v_{i_{d-1}}, v_{l_2(e)}\rangle.
\end{align*}

Using Proposition \ref{thm:homcoords} it is easy to show that for a one-dimensional
cone 
$v_i\notin \{v_{i_1}, \ldots, v_{i_{d-1}}, v_{l_1(e)}, v_{l_2(e)}\}$, the
$T_N$--action on $\O_{C_e}(\cocl_i)$ is trivial, we can thus disregard it
with respect to moving part of $H^0(C_e, f^*T\vt)$. For an $i_j\in
\{i_1, \ldots, i_{d-1}\}$, $\O_{C_e}(\lambda_{i_j})$ is the direction of
the pull back $f^*T\vt$ of the tangent space that corresponds to $V_{\tau_{i_j}}$,
where $\tau_{i_j}=\langle v_{i_1}, \ldots, \widehat{v_{i_j}}, \ldots, v_{i_{(d-1)}}, v_{l_1(e)}
\rangle$. 

To compute $e^{T_N}(H^{0,\move}(C_e, \O(\cocl_{i_j}))$, let us fix $i_j$ and write
$\tau=\tau_{i_j}$ and $\lambda_e^\gamma=\lambda_{i_j}$.
 If $\gamma\diamond\sigma_1$ such that $\gamma\cap\sigma_1=\tau$, then
the weight at $\sigma_1\in\vt$ of the action on the fiber is just $\we^{\sigma_1}_{\gamma}$, while
the weight at $\sigma_1$ of the action on the base $C_e$ is $\we_{F_1}$, where $F_1$ is the
flag $(\sigma_1(\gotv), e)$.

By analysis of the induced action on holomorphic sections of the bundle $\O_{C_e}(\lambda_e^\gamma)$,
we see that its weights are:
\[ \we^{\sigma_1}_{\gamma_k} -\frac{b}{d_e}\we^{\sigma_1}_{\sigma_2}, 
\quad b=1,\ldots, \lambda_e^{\gamma_k}\]
which are all non--trivial.
For $\O(\lambda_{l_1(e)})$ and $\O(\lambda_{l_2(e)})$, there is in each case
exactly one zero weight among these weights:
\begin{align*}
\O_{C_e}(\lambda_{l_1(e)}):&\quad
\we^{\sigma_1}_{\sigma_2}, \ldots, \frac{1}{d_e}\we^{\sigma_1}_{\sigma_2}, 0\\
\O_{C_e}(\lambda_{l_2(e)}):&\quad
0, -\frac{1}{d_e}\we^{\sigma_1}_{\sigma_2}, \ldots, -\we^{\sigma_1}_{\sigma_2}.
\end{align*}
Hence we have just proven the following lemma:

\begin{lemm}\mylabel{lem:ec3p1}
The $T_N$--equivariant Euler class of the moving part of the trivial bundle
$H^0(C_e, f^*T\vt)$ is:
\begin{multline*}
 e^{T_N} \left(H^{0,\move}(C_e, f^*T\vt)\right)=\\
 (-1)^{d_e}\frac{(d_e!)^2}{d_e^{2d_e}}\left(\we^{\sigma_1}_{\sigma_2}\right)^{2d_e}
\prod_{\sigma_2\neq\gamma\face\sigma_1}
\prod_{b=0}^{\lambda_e^{\gamma}}
\left(\we^{\sigma_1}_\gamma - \frac{b}{d_e}\we^{\sigma_1}_{\sigma_2}\right).
\end{multline*}\proofend
\end{lemm}

\begin{coro}\mylabel{coro:ec3p3}
The $T_N$--equivariant Euler class of the moving part of $H^1(C_e, f^*T\vt)$ is:
\[
e^{T_N}(H^{1,\move}(C_e, f^*T\vt))=\prod_{\sigma_2\neq\gamma\face\sigma_1}
\prod_{b=\lambda_e^\gamma}^{-2}
\left(\we^{\sigma_1}_\gamma - \frac{b+1}{d_e}\we^{\sigma_1}_{\sigma_2}\right).
\]
\end{coro}

\begin{proof}
Just apply Serre duality:
\[ H^1(C_e, \O(\lambda_e^\gamma))=\Hom(\O(\lambda_e^\gamma)=H^0(C_e,\O(\lambda_e^\gamma)\otimes
\omega_{C_e}),\]
where $\omega_{C_e}$ is the dualizing bundle of $C_e\cong\P^1$.
\end{proof}

So it only remains to compute the weights of the (trivial)
bundles $T_{p_{\sigma(\gotv)}}\vt$ and $T_{p_{\sigma(F)}}\vt$ which is
now straightforward:

\begin{lemm}\mylabel{lem:ec3p2}
For a maximal cone $\sigma\in\Sigma^{(d)}$, the equivariant Euler
class of the trivial bundle $T_{p_\sigma}\vt$ is equal to 
\[
e^{T_N}(T_{p_\sigma}\vt)=\prod_{\gamma\face\sigma} \we^\sigma_\gamma = \totwe{\sigma}.
\]
\end{lemm}

\begin{coro}
The $T_N$--equivariant Euler class
of the moving part of the difference bundle $B_2-B_5$ is given by: 
\mycomment{eq:ec3}
\begin{multline}\label{eq:ec3}
e^{T_N} (B_2^\move-B_5^\move) =
\prod_{\gotv\in\V} \left(\totwe{\sigma(\gotv)}\right)^{\val(\gotv)-1}\cdot\\
\prod_{\stackrel{\scriptstyle e\in\E}{\partial e=\{\gotv_1, \gotv_2\}}}
\left(
    \frac{(-1)^d d^{2d}}{(d!)^2 \left(\we^{\sigma_1}_{\sigma_2}\right)^{2d}}
    \prod_{\sigma_2\neq\gamma\face\sigma_1}
    \frac{\displaystyle \prod_{i=\lambda^\gamma_e+1}^{-1}
        \left( \we^{\sigma_1}_\gamma - \frac{i}{d}\cdot \we^{\sigma_1}_{\sigma_2}
        \right)}{\displaystyle
              \prod_{i=0}^{\lambda^\gamma_e}
        \left(\we^{\sigma_1}_\gamma - \frac{i}{d}\cdot \we^{\sigma_1}_{\sigma_2}
        \right)}
\right)_{\begin{array}{l} \scriptstyle \sigma_1=\sigma(\gotv_1)\\
              \scriptstyle \sigma_2=\sigma(\gotv_2)\\
              \scriptstyle d=d_e,
     \end{array}}
\end{multline}
where $\val: \V \longrightarrow {\mathbb N}$ is the number
$\val(\gotv)$ of flags at a vertex $\gotv\in\V$.
\end{coro}

\begin{proof}
Just assemble the long sequence (\ref{eq:ec3a}), Lemma \ref{lem:ec3p1}, 
Corollary \ref{coro:ec3p3}, and Lemma \ref{lem:ec3p2}. To obtain the
exponent of $\totwe{\sigma(\gotv)}$ just observe that there $\val(\gotv)$
flags at the vertex $\gotv$.
\end{proof}

\subsection{The main theorem}

In the previous subsections we have computed all contributions
(\ref{eq:ec1}), (\ref{eq:ec2}) and (\ref{eq:ec3}) entering the
formula (\ref{eq:eczerl}) for the equivariant Euler class of the
virtual normal bundle to $\cM_\mygraph$. Therefore, applying
Graber and Pandharipande's virtual Bott residue formula
(\ref{eq:bottres}) we obtain the following proposition
expressing the genus--zero Gromov--Witten invariants of a smooth
projective toric variety $\vt$ in terms of its $\modsp$--graphs
$\mygraph$ and the fan $\Sigma$:

\mycomment{thm:maintheorem}
\begin{theo}\label{thm:maintheorem}
The genus--zero Gromov--Witten invariants for a toric variety $\vt$ are given
by
\mycomment{eq:mainthm}
\begin{equation}\label{eq:mainthm}
\Phi^{\vt}_{m,A}(Z^{l_1}, \ldots, Z^{l_m}) = \sum_\mygraph
\frac{1}{|{\mathbf A_\mygraph}|}\int_{\M_\mygraph}
\frac{\prod_{j=1}^m\prod_{k=1}^n
\left(\we^{\sigma(j)}_k\right)^{l_{j,k}}}{e^{T_N}(N_\mygraph^{\text{virt}})},
\end{equation}
where
\begin{itemize}
\item we use the convention $0^0=1$;
\item $Z^{l_i}=Z_1^{l_{i,1}}\cdot\ldots\cdot Z_n^{l_{i,n}}$;
\item $\sigma: \{1,\ldots, m\} \rightarrow \Sigma^{(d)}$, the image $\sigma(j)$ of $j$ corresponding
to the fixed point the marked point $j\in\{1, \ldots, m\}$ is mapped to:
\[ \exists \gotv\in\V(\mygraph)\, : \sigma(\gotv)=\sigma(j) \; \wedge \; j\in S(\gotv); \]
\item we define $\we^{\sigma(j)}_k:=\left\{
\begin{array}{ll}
0&\quad \text{if $v_k \notin \Sigma^{(1)}_{\sigma(j)}$}\\
\we_\gamma^{\sigma(j)}&\quad \text{if $\gamma\face\sigma(j)$ and
                $v_k\in\Sigma^{(1)}_{\sigma(j)}\backslash\Sigma^{(1)}_{\gamma};$}
\end{array}\right. $
\item the inverse of the Euler class of the virtual normal bundle is given by
\begin{multline*}
\frac{1}{e^{T_N}(N_\mygraph^{\text{virt}})} = \prod_{\gotv\in\V}
\left(\totwe{\sigma(\gotv)}\right)^{\val(\gotv)-1}\cdot
\prod_{F\in\F_{\ge3}}
\frac{1}{\we_F-e_F}\\\cdot\prod_{F\in\F_1}\we_F\cdot
\prod_{\gotv\in\V_{2,0}}\frac{1}{\we_{F_1(\gotv)}+\we_{F_2(\gotv)}}\cdot\\
\prod_{\stackrel{\scriptstyle e\in\E}{\partial e=\{\gotv_1,
\gotv_2\}}} \left(
    \frac{(-1)^d d^{2d}}{(d!)^2 \left(\we^{\sigma_1}_{\sigma_2}\right)^{2d}}
    \prod_{\sigma_2\neq\gamma\face\sigma_1}
    \frac{\displaystyle \prod_{i=\lambda^\gamma_e+1}^{-1}
        \left( \we^{\sigma_1}_\gamma - \frac{i}{d}\cdot \we^{\sigma_1}_{\sigma_2}
        \right)}{\displaystyle
              \prod_{i=0}^{\lambda^\gamma_e}
        \left(\we^{\sigma_1}_\gamma - \frac{i}{d}\cdot \we^{\sigma_1}_{\sigma_2}
        \right)}
\right)_{\begin{array}{l} \scriptstyle \sigma_1=\sigma(\gotv_1)\\
              \scriptstyle \sigma_2=\sigma(\gotv_2)\\
              \scriptstyle d=d_e.
     \end{array}}
\end{multline*}
\end{itemize}
\end{theo}

\section{Simplifying the formula for the Gromov--Witten invariants}\mylabel{sec:sumgraph}

Theorem \ref{thm:maintheorem} in the previous section provides a formula
for the Gromov--Witten invariants of symplectic toric manifolds. This 
formula, however, is highly combinatorial: in particular, the number
of fixed point components $\M_\Gamma$ can rise very quickly.

In this section, we will give a first combinatorial simplification of
this formula by grouping together the terms of several fixed point 
components. However, the
results obtained are still far away from what one might wish for ---
see the discussion at the end of the chapter. 

First, we will attack the only remaining (classical) cohomology
classes in the formula:

\begin{lemm}\mylabel{lem:simplifya}
Let $P_n(x_1, \ldots, x_k)$ be the homogeneous polynomial of
degree $n$ in the variables $x_1, \ldots, x_k$, that is
\[ P_n(x_1,\ldots, x_k) =
\sum_{\begin{array}{c}\scriptstyle d_1+\cdots+d_k=n\\ \scriptstyle
d_i\ge0\end{array}} \prod_i x_i^{d_i}.
\] 
Further, for a vertex
$\gotv\in\V$, let $\deg(\gotv)$ be the number of flags (including
marked points) at $\gotv$, and let $\val(\gotv)$ the number of
flags at $\gotv$. Then
\[ \int\limits_{\cM_\mygraph} \prod_{F\in\F_3} \frac{1}{\we_F-e_F} = \int\limits_{\cM_\mygraph} \prod_{\gotv\in\V_3}
P_{\deg(\gotv)-3}\left(\frac{e_{F_1}}{\we_{F_1}}, \ldots,
\frac{e_{F_{\val(\gotv)}}}{\we_{F_{\val(\gotv)}}}\right)
\prod_{i=1}^{\val(\gotv)} \frac{1}{\we_{F_i}},\] where as before
the $\we_F$ are the weights of the flag $F$, and $e_F$ is the
universal cotangent line at $F$.
\end{lemm}

\begin{proof}
First we observe that formally
\begin{align*}
\frac{1}{\we_F-e_F}&= \frac{1}{\we_F} + \frac{e_F}{\we_F^2} +
\frac{e_F^2}{\we_F^3} + \cdots\\ &= \frac{1}{\we_F}
\sum_{i=0}^\infty \left(\frac{e_F}{\we_F}\right)^i.
\end{align*}
Since $\cM_\mygraph$ is the product of Deligne--Mumford spaces of
stable curves $\dmsp[m]$ corresponding each to some vertex $\gotv$
of the graph $\mygraph$, the integral above is the product of
integrals over the Deligne--Mumford spaces $\dmsp[m]$ of the
classes corresponding to their associated vertex. So let
$\gotv\in\V_3$ be a vertex with at least three special points,
$k=\val(\gotv)$ be the number of flags at $\gotv$, and $F_1,
\ldots, F_k$ be the flags of $\gotv$. Then
\[ \prod_{j=1}^k \frac{1}{\we_{F_j}-e_{F_j}} = \prod_{j=1}^k
\frac{1}{\we_j}\sum_{i=0}^\infty \left(\frac{e_F}{\we_F}\right)^i.
\]
The dimension of the moduli space $\dmsp[\deg(\gotv)]$ corresponding 
to the vertex $\gotv$ is equal to $\deg(\gotv)-3$, hence
\[
\int\limits_{\dmsp[\deg(\gotv)]} \prod_{j=1}^k
\frac{1}{\we_j}\sum_{i=0}^\infty \left(\frac{e_F}{\we_F}\right)^i
= \int\limits_{\dmsp[\deg(\gotv)]} P_{\deg(\gotv)-3}
\left(\frac{e_{F_1}}{\we_{F_1}}, \ldots, \frac{e_{F_k}}{\we_{F_k}}
\right) \prod_{i=1}^k \frac{1}{\we_{F_i}}.
\]
\end{proof}

\begin{lemm}\mylabel{lem:simplifyb}
Let $\dmsp[m]$ be a Deligne--Mumford space of stable curves, and
let $e_{F_1},\ldots, e_{F_k}$ be universal cotangent lines to
different marked points of $\dmsp[m]$. Then
\[ \int\limits_{\dmsp[m]} P_{m-3}\left(\frac{e_{F_1}}{\we_{F_1}}, \ldots,
\frac{e_{F_k}}{\we_{F_k}}\right) =
\left(\frac{1}{\we_{F_1}}+\cdots+\frac{1}{\we_{F_k}}\right)^{m-3}.
\]
\end{lemm}

\begin{rema}
Note that the Lemma is about Deligne--Mumford spaces, that is
about ``factors'' of a fixed point component $\cM_\Gamma$.
Therefore the parameter $m$ is not the number of marked points in
$\cM_\Gamma$, it is rather the number of special points of a
``factor'' corresponding to a vertex.
\end{rema}

\begin{proof}
Remember that by Corollary \ref{cor:dmexplicit} we have
\[ \int\limits_{\dmsp[m]} e_{F_1}^{d_1} \wedge \ldots \wedge
e_{F_k}^{d_k} = \frac{(m-3)!}{\prod_{i=1}^k d_i!}.
\]
Hence, by the binomial formula, we obtain
\begin{align*}
 \int_{\dmsp[m]} P_{m-3}\left(\frac{e_{F_1}}{\we_{F_1}}, \ldots,
\frac{e_{F_k}}{\we_{F_k}}\right) &= \sum_{\sum_i d_i=m-3}
\frac{(m-3)!}{d_1!\cdots d_k!}
\left(\frac{1}{\we_{F_1}}\right)^{d_1}\cdots
\left(\frac{1}{\we_{F_k}}\right)^{d_k}\\
&=\left(\frac{1}{\we_{F_1}}+\cdots+\frac{1}{\we_{F_k}}\right)^{m-3}.
\end{align*}
\end{proof}

Having totally eliminated integration from the formula for the
calculation of the Gromov--Witten invariants, we will now examine
the part of the sum, that depends not only on the underlying graph
but also on where the marked points are placed on the graph. To
this effect we define the following two terms:
\begin{align*}
T_\mygraph=&
\prod_{t=1}^\infty
\prod_{\gotv\in\V_{t,*}(\mygraph)}
\left(
   \totwe{\sigma(\gotv)}
\right)^{t-1} \cdot  
\left(
   \prod_{i=1}^t
   \frac{1}{\we_{F_i(\gotv)}}
\right) \cdot
\left(
   \frac{1}{\we_{F_1(\gotv)}}+\cdots+
   \frac{1}{\we_{F_t(\gotv)}}
\right)^{t-3} \cdot\\
&\cdot\prod_{\stackrel{\scriptstyle e\in\E}{\partial e=\{\gotv_1,
\gotv_2\}}} \left(
    \frac{(-1)^d d^{2d}}{(d!)^2 \left(\we^{\sigma_1}_{\sigma_2}\right)^{2d}}
    \prod_{\sigma_2\neq\gamma\face\sigma_1}
    \frac{\displaystyle \prod_{i=\lambda^\gamma_e+1}^{-1}
        \left( \we^{\sigma_1}_\gamma - \frac{i}{d}\cdot \we^{\sigma_1}_{\sigma_2}
        \right)}{\displaystyle
              \prod_{i=0}^{\lambda^\gamma_e}
        \left(\we^{\sigma_1}_\gamma - \frac{i}{d}\cdot \we^{\sigma_1}_{\sigma_2}
        \right)}
\right)\\
S_\mygraph=& \left[ \prod_{t,s} \prod_{\gotv\in\V_{t,s}(\mygraph)}
\left(\frac{1}{\we_{F_1(\gotv)}}+\cdots+\frac{1}{\we_{F_t(\gotv)}}\right)^s\right]
\cdot \prod_{j=1}^m \prod_{k=1}^n (\we^{\sigma(j)}_k)^{l_{j,k}}.
\end{align*}
Note that in fact $T_\mygraph$ only depends on the graph
type $\gtype$ of $\mygraph$ and not on where the marked points are placed. 
Before we go on, we will sum up what we have proven so far:

\begin{coro}\mylabel{cor:simplifya}
With the notation as in Theorem \ref{thm:maintheorem}, the genus--zero 
$m$--point Gromov--Witten invariants are given by:
\[ \Phi^{\vt}_{m,A}(Z^{l_1}, \ldots, Z^{l_m}) = \sum_\mygraph
\frac{1}{|{\bf A}_\mygraph|} T_\mygraph \cdot S_\mygraph.
\]
\end{coro}

\begin{proof}
By lemmata \ref{lem:simplifya} and \ref{lem:simplifyb}, the
formula is obviously true for the parts coming from vertices with
at least three flags (\ie. special points). Thus it remains to
show that
\begin{multline*}
\prod_{F\in\F_1}\we_F \cdot \prod_{\gotv\in\V_{2,0}}
\frac{1}{\we_{F_1(\gotv)}+\we_{F_2(\gotv)}} = \left[
\prod_{\gotv\in\V_{1,0}} \frac{1}{\we_{F(\gotv)}} \cdot
\left(\frac{1}{\we_{F(\gotv)}}\right)^{-2} \right]\\ 
\cdot \left[
\prod_{\gotv\in\V_{2,0}}
\frac{1}{\we_{F_1(\gotv)\we_{F_2(\gotv)}}} \cdot
\left(\frac{1}{\we_{F_1(\gotv)}}+\frac{1}{\we_{F_2(\gotv)}}\right)^{-1}
\right]\\
\cdot \left[
\prod_{\gotv\in\V_{1,1}} \frac{1}{\we_{F(\gotv)}}\cdot
\left(\frac{1}{\we_{F(\gotv)}}\right)^{-2}\cdot
\frac{1}{\we_{F(\gotv)}} \right] .
\end{multline*}
The first term on the right hand side of the equation obviously
coincides with the first term on the left hand side. The third
factor on the right hand side is trivial, while the second term on
the right hand side is equal to the second term on the left, since
\[ \frac{1}{a+b} = \frac{1}{a b}
\left(\frac{a+b}{a b}\right)^{-1}=\frac{1}{a
b}\left(\frac{1}{a}+\frac{1}{b}\right)^{-1}.
\]
\end{proof}

Remember that the $T_{\gtype}$ term only depends on the graph type.
Hence we could further simplify
the formula for the Gromov--Witten invariants if we were able to
explicitly compute the sum of all $S_\mygraph$ over all $\mygraph\in
\typemap^{-1}(\gtype)$ corresponding to a graph type $\gtype$.

Since the Gromov--Witten invariants are commutative with respect
to the cohomology classes, we can suppose that these classes are
\[ \underbrace{Z^{l_1}, \ldots, Z^{l_1}}_{m_1-times}, \ldots,
\underbrace{Z^{l_q}, \ldots, Z^{l_q}}_{m_q-times}, \qquad
m_1+\cdots+m_q=m,
\]
such that the multi--indices $l_i$ are pairwise different. The
positions of the marked points on the graph are independent from
each other, so we can consider each of the $Z^{l_i}$ on its own.
Let $\we_{\sigma(\gotv)}^{l_i}$ be the weight of the class at the
vertex $\gotv$:
\[\we_{\sigma(\gotv)}^{l_i} :=
\prod_{k=1}^n \left(\we_k^{\sigma(\gotv)}\right)^{l_{i,k}}. \]
Note that as before we use the convention $0^0=1$, but of course
$0^n=0$ for a positive integer $n$. We will also use the notation
\[ \totwe{F(\gotv)}:=\left(\frac{1}{\we_{F_1(\gotv)}}+\cdots+
\frac{1}{\we_{F_{\val(\gotv)}(\gotv)}}\right)
\]
for the {\em total weight} of the flags at the vertex $\gotv$ that
is not to be confused with the total weight $\totwe{\sigma}$ at a
fixed point $\sigma$.

\begin{rema}
It is worthwhile to notice that the weight polynomial
$\we_{\sigma(\gotv)}^{l_i}$ only depends on $\topmap(\mygraph)$,
\ie. only on the fixed point
$\sigma(\gotv)$ in $X_\Sigma$, but neither on the graph $\Gamma$
nor on a vertex $\gotv$ in such a graph.

On the other hand, the weight polynomial $\totwe{F(\gotv)}$
depends on both the topology of the graph type and the
multiplicity decorations.
\end{rema}

Let $S_{\gtype, l_i}$ be the part of the sum
$\sum\limits_{\mygraph\in\typemap^{-1}(\gtype)} S_{\mygraph}$
corresponding to the $m_i$ classes $Z^{l_i}$:
\mycomment{eq:sgammadef}
\begin{equation}\label{eq:sgammadef}
S_{\gtype, l_i}:=\sum_{(\gotv_1, \ldots,
\gotv_{m_i})} \prod_{t=1}^{m_i} \totwe{F(\gotv_t)} \cdot
\we_{\sigma(\gotv_t)}^{l_i},
\end{equation}
where the sum is running over all ordered $m_i$--tuples $(\gotv_1,
\ldots, \gotv_{m_i})$ of vertices of the graph type
$\gtype$.

\begin{lemm}\mylabel{lem:simplifyc}
With the notation from above, the following holds:
\begin{enumerate}
\item  Let ${\mathbf A}_{\gtype}$ be the automorphism group of the
graph type $\gtype$
\[ \sum\limits_{\mygraph\in\typemap^{-1}(\gtype)}
\frac{S_{\mygraph}}{|{\mathbf A}_{\mygraph}|}
= \frac{\prod\limits_{i=1}^q S_{\gtype, l_i} }{|{\mathbf
A}_{\gtype}|};
\]
\item Let $\{\gotv_1, \ldots, \gotv_r\}$ be the set of vertices of
the graph type $\gtype$. Then
\mycomment{eq:movea}
\begin{equation}\label{eq:movea}
S_{\gtype, l_i} = \left( \sum_{j=1}^r \totwe{F(\gotv_j)} \cdot
\we_{\sigma(\gotv_j)}^{l_i}\right)^{m_i}.
\end{equation}
\end{enumerate}
\end{lemm}

\begin{rema}
In the second part of the lemma, of course we only have to
consider vertices $\gotv$ for which the weight
$\we_{\sigma(\gotv)}^{l_i}$ of $Z^{l_i}$ at $\gotv$ is non--zero.
Thus, in practice, the number $r$ of vertices to consider will
vary with the class $Z^{l_i}$ and be considerably lower than the
number of all vertices in the graph type $\gtype$.
\end{rema}

\begin{proof}
Given a graph type $\gtype$, a graph $\mygraph$ in this graph type
is given by the positions
\[ \gotv^1_1, \ldots, \gotv^1_{m_1},\ldots, \gotv^q_{m_q}
\]
of the $m$ marked points (where the classes $Z^{l_1},
\ldots, Z^{l_1}, \ldots, Z^{l_q}, \ldots, Z^{l_q}$ are attached), and vice versa
any such tuple $(\gotv^1_1, \ldots, \gotv^q_{m_q})$ of vertices of
$\gtype$ yields a graph $\mygraph$. Hence, up to automorphisms,
the first part of the lemma is obvious. For the automorphism
groups, first remark that ${\mathbf A}_\mygraph$ injects into
${\mathbf A}_{\gtype}$. The cokernel consists of automorphisms of
$\gtype$ that change the location of the marked points on the
graph $\mygraph$ and hence reflect multiply counted instances of
the same graph in $S_{\gtype, l_i}$.

To prove the second part, we observe that
\[ \sum_{(\gotv_1, \ldots,
\gotv_{m_i})} \prod_{t=1}^{m_i} \totwe{F(\gotv_t)} \cdot
\we_{\sigma(\gotv_t)}^{l_i} = \sum_{\begin{array}{c}\scriptstyle
(d_1,\ldots,d_r)\\ \scriptstyle d_1+\cdots+d_r=m_i\end{array}}
\frac{m_i!}{d_1!\cdots d_r!} \prod_{j=1}^r
\left(\totwe{F(\gotv_j)}\cdot
\we_{\sigma(\gotv_j)}^{l_i}\right)^{d_j},
\]
yielding the desired expression by applying the binomial formula.
\end{proof}

We will now further analyze the $S_{\gtype,l_i}$ terms. In fact,
we will show that these terms are equal for two graph types as
long as their image on the moment polytope is the same.

\begin{coro}\mylabel{cor:sgtop}
Let $\gtype[1]$ and $\gtype[2]$ be two graph types of a smooth
projective toric variety $X_\Sigma$ that have the same topological
graph type: $\topmap(\gtype[1])=\topmap(\gtype[2])=\gtop$. Then the
``S--terms'' of the two graph types coincide:
\[ S_{\gtype[1], l_i} = S_{\gtype[2], l_i}.\]
In particular we define $S_{\gtop, l_i}= S_{\gtype[1], l_i}$.
\end{coro}

\begin{rema}
Note that the corollary only states the equality of the
``S--terms'' for the two graph types, but not of their
automorphism groups that might well be different.
\end{rema}

\begin{proof}
Let $\gtype$ be a graph type of $X_\Sigma$. Remember, that the
weights $\totwe{F(\gotv_j)}$ in the sum of the right hand side of
(\ref{eq:movea}) are given by
\[
\totwe{F(\gotv_j)}=\frac{1}{\we_{F_1(\gotv_j)}} + \cdots
 +\frac{1}{\we_{F_{n_j}(\gotv_j)}}
\]
where $F_1(\gotv_j),\ldots,F_{n_j}(\gotv_j)$ are the flags at
$\gotv_j$ corresponding to its $n_j$ edges. Hence we can write
\[ \sum_{j=1}^r \totwe{F(\gotv_j)} \cdot \we^{l_i}_{\sigma(\gotv_j)} =
\sum_{F\in\F(\gtype)} \frac{\we^{l_i}_{\sigma(F)}}{\we_F}.
\]
Note that as the notation suggests, $\we^{l_i}_{\sigma(F)}$ only
depends on the class $Z^{l_i}$ and the vertex $\sigma(F)$ of the
polytope $P_\Sigma$, hence it is well defined for the image flag
$\topmap_*(F)$ on $\gtop=\topmap(\gtype)$.

For a flag $F\in\F(\gtype)$, let $F=(\gotv_1, e)$ and $\partial
e=\{\gotv_1,\gotv_2\}$. Then the weight $\we_F$ is defined by
\[ \we_F = \frac{\we^{\sigma(\gotv_1)}_{\sigma(\gotv_2)}}{d_e},\]
hence is given by the two vertices $\sigma(\gotv_1),
\sigma(\gotv_2)$ of the polytope $P_\Sigma$, and the multiplicity
of the edge $e$. Accordingly, for a flag $(\gotv_{\text{top},1},
e_{\text{top}}) = F_{\text{top}}\in\F(\gtop)$ in a topological
graph type $\gtop$ corresponding to an edge
$e_{\text{top}}\in\E(\gtop)$ with multiplicity
$d_{e_{\text{top}}}$ and ends $\partial e_{\text{top}} =
\{\gotv_{\text{top},1}, \gotv_{\text{top},2}\}$, let us define the
weight
\[ \we_{F_{\text{top}}} := \frac {\we^{\sigma(\gotv_{\text{top},1})}_{\sigma(\gotv_{\text{top},2})}}{d_{e_{\text{top}}}}.\]
Then for a graph type $\gtype$ such that $\topmap(\gtype)=\gtop$
one easily sees that
\[ \frac{1}{\we_{F_{\text{top}}}} = \sum_{F\in \F(\gtype)\cap\topmap^{-1}_*(F_{\text{top}})} \frac{1}{\we_F}.
\]
Hence we obtain that indeed
\begin{multline*}
S_{\gtype, l_i} = \left(\sum_{j=1}^r \totwe{F(\gotv_j)} \cdot 
\we^{l_i}_{\sigma(\gotv_j)} \right)^{m_i}=\\
= \left(\sum_{\gotv\in\V(\topmap(\gtype))} \totwe{F(\gotv)} \cdot
\we^{l_i}_{\sigma(\gotv)}\right)^{m_i} = S_{\gtop, l_i}
\end{multline*}
only depends on the topological graph type $\gtop$.
\end{proof}

\begin{rema}
Note that although that topological graphs are no longer required
to be without loops, we can use the ``same'' formulas to compute
$S_{\gtop, l_i}$, \ie. we do not need to choose a graph type
$\gtype$ and then compute $S_{\gtype, l_i}$!
\end{rema}

\subsection{Concluding remarks}

\noindent Before we go on to the examples, however, let us make some concluding remarks:
\begin{enumerate}
\item The formula in the Main Theorem \ref{thm:maintheorem} 
does not yield directly the $m$--point Gromov--Witten invariants
(for $m>3$) needed in the computation of quantum products with more than two factors, \ie. the
invariants
\mycomment{eq:gwimpt}
\begin{equation}\label{eq:gwimpt}
\Psi^\vt_{m,A}([pt]; \alpha_1, \ldots, \alpha_m).
\end{equation}
It does compute the invariants
\[ \Psi^\vt_{m,A}(1; \alpha_1, \ldots, \alpha_m).\]
Since for $m=3$, the Deligne--Mumford space of stable curves is just a point, $\dmsp[3]=[pt]$,
the theorem gives the three--point Gromov--Witten invariants needed for computing quantum products
of two factors: $\alpha\star\beta$.

This, however, is no real disadvantage, since the decomposition law for Gromov--Witten invariants
expresses the $m$--point invariants in (\ref{eq:gwimpt}) with the help of the three--point invariants.
\item In this article, we have not tackled the case of higher genus Gromov--Witten invariants for two
reasons:
\begin{itemize}
\item First of all we have been interested in a better understanding of the quantum cohomology of
projective toric manifolds with the hope of eventually computing it for non--Fano manifolds as well.
\item Secondly, even for genus--zero invariants the formula for the virtual normal bundle becomes
combinatorically quite complicated.
\end{itemize}
However, generalizing the above theorem to higher--genus
Gromov--Wit\-ten invariants should essentially work the same way as in the case of
complex projective space, that has been studied by Graber and
Pandharipande (see \cite{gp97}). Note, however, that the fixed
point components $\cM_\mygraph$ will then contain higher genus
Deligne--Mumford spaces $\dmsp$ as well, complicating the
computation of the integrals in (\ref{eq:mainthm}): there is no
longer an explicit formula such as in corollary \ref{cor:dmexplicit},
but only an recursive formula.
\item The simplifications to the combinatorial nature of the formula (\ref{eq:mainthm})
can only be a first step if one really wants to work with the Gromov--Witten invariants
and the quantum cohomology ring of symplectic toric manifolds. To obtain a closed 
formula of the quantum cohomology ring using this approach, what one really would
like to have is a formula reading off the Gromov--Witten invariants say directly
from the moment polytope.
\end{enumerate}


\newcommand{\ml}[1][1]{\unitlength1mm\hspace{-0.35cm}
                       \begin{picture}(10,5)
                           \put(0,1.65){\line(1,0){13.5}}
                           \put(4,2.25){\parbox{0.4cm}{\centering\small #1}}
                       \end{picture}
                       }
\newcommand{\mls}[1][1]{\unitlength1mm\hspace{-0.35cm}
			 \begin{picture}(6,5)
			     \put(0,1.65){\line(1,0){9.5}}
			     \put(2,2.25){\parbox{0.4cm}{\centering\small #1}}
			 \end{picture}
			 }
\newcommand{\bbull}[1]{\mbox{\begin{minipage}[t]{0.7cm}
                \begin{center}\rule{0pt}{0.7cm}
                    \huge$\bullet$\\[-0.45cm]
                    \parbox{0.7cm}{\centering\tiny$#1$}
                \end{center}
                 \end{minipage}
                 }
              }
\newcommand{\fork}[4]{\mbox{\unitlength1mm
                            \begin{picture}(33,16)
                              \put(0,0){\bbull{#1}}
                              \put(13,0){\bbull{#2}}
                              \put(26,0){\bbull{#3}}
                              \put(12,13){
                                \mbox{\begin{minipage}[t]{0.7cm}
                                        \begin{center}\rule{0pt}{0.7cm}
                                          \huge$\bullet$\\[-0.45cm]
                                          \parbox{0.7cm}{\hfill\tiny$#1$}
                                        \end{center}
                                      \end{minipage}
                                     }
                              }
                              \put(3.5,1.5){\line(1,0){26}}
                              \put(16.5,0){\line(0,1){13.5}}
                            \end{picture}
                           }
                      }

\section{Examples}\mylabel{sec:examples}

In this section we give three examples of actual computations using
the localization formula applied to toric manifolds (Theorem
\ref{thm:maintheorem}). We first compute  the quantum cohomology ring
of the projective space --- this ring is of course well known, but
this makes it also a good example to experiment with our formula.
Next we compute the invariants and the quantum cohomology of the
Fano threefold $\P(\O_{\C\P^2}(2)\oplus 1)$. The reason why we
have chosen this specific example is that, as far as we know, it
is the simplest smooth projective toric variety that does not
appear in previous work, \eg. this example is not accessible by
the methods used in \cite{qr98}. Finally, we will compute some invariants
of the non--Fano threefold $\P(\O_{\C\P^2}(3)\oplus 1)$ to show that
the quantum cohomology ring of this manifold is not the one defined
in \cite{bat93}.

\subsection{The projective space}

Before we actually start computing the Gromov--Witten invariants of any complex projective
space $\P^n$ using the above formula, we will reduce the number of invariants
for which we actually have to use the formula. Remember that for Gromov--Witten invariants
the so--called composition law holds:
\[ \Phi_4^{A,X} ([pt];\alpha_1,\alpha_2,\alpha_3,\alpha_4)= \sum_{A=A_1+A_2} \sum_{i=1}^N
\Phi_3^{A_1,X}(\alpha_1,\alpha_2, \beta_i)\Phi_3^{A_2,X}(\beta^i, \alpha_3, \alpha_4), \]
where $(\beta_1,\ldots,\beta_N)$ is a basis of $H^*(X,\Z)$, $(\beta^1,\ldots, \beta^N)$ its
dual basis in the same space, and $[pt]\in\dmsp[4]$ the cohomology class Poincar\'e dual to
a point.

For the complex projective space $\P^n$, let $H\in H^2(\P^n, \Z)$
be the generator of degree--2 cohomology. Then $(1,H, H^2, \ldots, H^n)$ is a basis of $H^*(X,\Z)$
whose dual is $(H^n,\ldots,1)$. So any class $A\in H_2(\P^n, \Z)$ is necessarily a multiple of
the class $H$, $A=kH$, and if $A$ contains holomorphic curves $k$ needs to be positive. Remember
that the virtual dimension of the moduli stack $\cM_{0,3}(\P^n, A)$ is equal to
\[ \dim_\virt \cM_{0,3}(\P^n, A) = \langle c_1(\P^n), A\rangle +n + 3 -3 = k(n+1) +n. \]
Hence, for $k>1$, the virtual dimension of the moduli space is bigger than $3n$. Therefore
there can only be non--trivial Gromov--Witten invariants for the class\footnote{And of course
for the trivial class $A=0$, given by the intersection numbers (of the usual cup product).} $A=H$.

So let us look at the composition law for $A=H$. First of all, the
dimension of the (virtual) fundamental class of the moduli stack
$\M_{0,m}(\P^n, H)$ is equal to $\langle c_1(\P^n), H\rangle+ n+ m-3= 2n +m -2$. Also, we can not
decompose the class $H$ into effective classes, \ie. classes that contain again holomorphic curves.
Suppose that $p\ge q\ge r\ge 2$ and $p+q+r=2n+1$. Hence we obtain:
\begin{align*}
\Phi_3^H(H^p, H^q, H^r)&=\Phi_3^H(H^p,H^q,H^r)\cdot\Phi_3^0(H^{n-r}, H^{r-1}, H)\\
&=\Phi_4^H([pt];H^p, H^q, H^{r-1}, H) \quad \text{(since $p+q>n$)}\\
&=\Phi_4^H([pt];H^p, H^{r-1}, H^q, H) \quad \text{(since GWI are commutative)}\\
&=\Phi_3^H(H^p, H^{r-1}, H^{q+1}) \cdot \Phi_3^0(H^{n-q-1}, H^q, H)\\
&=\Phi_3^H(H^p, H^{q+1}, H^{r-1}).
\end{align*}
Therefore, by induction on $r$ we get
\[ \Phi_3^H(H^p, H^q, H^r) = \Phi_3^H(H^n, H^n, H),\]
that is we only have to use the fixed--point formula to compute one single Gromov--Witten invariant
for each complex projective space $\P^n$.

So let $(e_1,\ldots, e_n)$ be a basis of $\Z^n$, and $\Sigma$ be the fan given by the following
$1$--skeleton and set of primitive collections:
\begin{gather*}
v_1=e_1, \ldots, v_n=e_n, v_{n+1}=-e_1-\ldots-e_n\\
{\mathcal P} = \left\{\{v_1,\ldots, v_{n+1}\}\right\}.
\end{gather*}
We will denote the $n+1$ different $n$-dimensional cones in the fan $\Sigma$ as follows:
\[ \sigma_i = \langle v_1, \ldots, \hat{v_i}, \ldots, v_{n+1}\rangle, \]
where the element with the hat has to be omitted. One easily sees that the weights at $\sigma_i$
on the edge connecting to $\sigma_j$ is given by
\[ \omega_{\sigma_j}^{\sigma_i} = \we_j-\we_i.\]
As usual we will denote by $Z_1,\ldots, Z_{n+1}$ the $(\C^*)^n$--divisors in $\P^n$ coming from
the hyperplanes $\{z_i=0\}\subset\C^{n+1}$. So let us compute the invariant
$\Phi_3^H(H^n,H^n,H)$:

\begin{multline*}
\Phi_3^H(Z_1Z_2\cdots Z_n, Z_2Z_3\cdots Z_{n+1}, Z_1)
=\bbull{\sigma_{n+1}} \ml[\ \ 1] \bbull{\sigma_1}\\
=-\frac{(\we_{1}-\we_{n+1})\cdots(\we_n-\we_{n+1})\cdot(\we_2-\we_1)\cdots(\we_{n+1}-\we_1)\cdot(\we_1-\we_{n+1})}%
{(\we_1-\we_{n+1})^3\cdot(\we_2-\we_{n+1})
(\we_2-\we_1)\cdots(\we_n-\we_{n+1})(\we_n-\we_1)}\\ 
\end{multline*}
Summing up, we get the following results for the projective space $\C\P^n$:
\begin{itemize}
\item The only non--trivial genus--zero three--point Gromov--Witten invariants are
\begin{enumerate}
\item $\Phi_3^0(H^p,H^q, H^r)=1$ if $p+q+r=n$; and
\item $\Phi_3^H(H^p,H^q, H^r)=1$ if $p+q+r=2n+1$.
\end{enumerate}
\item Its (small) quantum cohomology ring is given by
\[ QH^*(\C\P^n, \C) = \C[H, q]/_{\langle H^{n+1}-q\rangle}. \]
\end{itemize}

\subsection{The $3$--folds $\P(\O_{\P^2}(m) \oplus 1)$: some general computations}

Before we actually compute Gromov--Witten invariants of two of these projective
bundles, we do some preparatory computations for the general case, \eg. we describe
their fan, derive their moment polytope and compute their cohomology ring and weight
table.

The manifold $\P(\O_{\P^2}(m) \oplus 1)$ is $3$--dimensional, hence its fan
lives in the lattice $N:={\Z}^3$. Let $(e_1, e_2, e_3)$ be the standard basis of
${\Z}^3$. Then the one--dimensional cones of the fan $\Sigma$ that
corresponds to the manifold $\P(\O_{\P^2}(m) \oplus 1)$ are generated by the
following vectors:
\begin{alignat*}{3}
\gotv_1&=e_1&\qquad \gotv_2&=-e_1\\
\gotv_3&=e_2&\qquad\gotv_4&=e_3&\qquad\gotv_5&=-e_2-e_3-me_1.
\end{alignat*}
The set of primitive collections of the fan $\Sigma$ is given by:
\[ {\mathfrak P} = \left\{ \{ \gotv_1, \gotv_2 \}, \{ \gotv_3, \gotv_4, \gotv_5\}\right\}.
\]
In fact, it is easy to recover the projective bundle structure over the
two--dimensional projective space from this description: $\gotv_1$ and $\gotv_2$
describe the fiber, a ${\C\P}^1$, while the other three generators
describe the ${\C\P}^2$ of the base space where the twisting of
the bundle is reflected in the $(-me_1)$--term in $\gotv_5$.

From the description by the set of primitive collections, 
it easy to get the set of maximal cones
in $\Sigma$. It consists of the following $3$--dimensional cones:
\begin{alignat*}{2}
\sigma_1&=\langle \gotv_1, \gotv_3, \gotv_4\rangle&\qquad
\sigma_4&=\langle \gotv_2, \gotv_3, \gotv_4\rangle\\
\sigma_2&=\langle \gotv_1, \gotv_3, \gotv_5\rangle&\qquad
\sigma_5&=\langle \gotv_2, \gotv_3, \gotv_5\rangle\\
\sigma_3&=\langle \gotv_1, \gotv_4, \gotv_5\rangle&\qquad
\sigma_6&=\langle \gotv_2, \gotv_4, \gotv_5\rangle.
\end{alignat*}

The degree--$2$ homology $H_2(X_\Sigma, {\Z})$
of $X_\Sigma=\P(\O_{\P^2}(m) \oplus 1)$ can be identified
with the group $R(\Sigma)\subset {\Z}^5$ given by
(Proposition \ref{prop:homcomtor}):
\[ R(\Sigma) := \left\{ \lambda=(\lambda^1, \ldots, \lambda^5) \in {\Z}^5
\, | \, \lambda^1\gotv_1+\cdots+\lambda^5\gotv_5=0 \right\}.
\]
The two vectors $\lambda_1:=(1,1,0,0,0)$ and
$\lambda_2=(0,-m,1,1,1)$ generate the ring $R(\Sigma)$. Moreover,
it is straightforward to check that they generate the effective
cone of $\Sigma$, \ie. the cone of homology classes that contain
holomorphic curves. Also note that by the same Proposition,
we have that
\[ \langle c_1(X_\Sigma), \lambda_1\rangle = 2 \qquad \text{and}
\qquad \langle c_1(X_\Sigma), \lambda_2\rangle = 3-m,
\]
where $c_1(X_\Sigma)$ is the first Chern class of the tangent bundle
of $X_\Sigma$.

Remember that the cohomology ring $H^*(X_\Sigma, {\C})$ is equal to the
quotient
\[ H^*(X_\Sigma, {\C}) = {\C}[Z_1, \ldots, Z_5]\left/_{
\textstyle \langle{\rm SR}(\Sigma) + {\Lin}(\Sigma)\rangle,}\right.
\]
where ${\rm SR}(\Sigma)$ is the {\em Stenley--Reisner ideal} of $\Sigma$, and
${\Lin}(\Sigma)$ is the ideal generated by linear relations. The former
is generated by monomials given by the set of primitive collections:
\[ {\rm SR}(\Sigma) = \langle Z_1 Z_2, Z_3 Z_4 Z_5 \rangle. \]
For the latter, the ideal ${\rm Lin}(\Sigma)$, let $u_1, \ldots, u_3$ be
any $\Z$--basis of the lattice $M$ dual to $N$, \eg. $u_i=e_i^*$.
Then, the ideal is generated the three terms $\sum_j \langle \gotv_j, u_i
\rangle Z_j$, \ie.
\[ Z_1-Z_2-mZ_5, \qquad Z_3-Z_5, \qquad Z_4-Z_5.\]
Therefore the cohomology ring of $X_\Sigma=\P(\O_{\P^2}(m) \oplus 1)$ is
equal to:
\begin{align*} H^*(X_\Sigma, {\C})& = {\C}[Z_1, \ldots, Z_5]\left/_{\langle
\textstyle \begin{array}{l}Z_1-Z_2-mZ_5, Z_3-Z_4,\\ Z_3-Z_5, Z_1Z_2, Z_3Z_4Z_5\end{array} \rangle}\right.\\
& = {\C}[Z_1, Z_3]\left/_{\langle
\textstyle Z_1^2-mZ_1Z_3, Z_3^3 \rangle.}\right.
\end{align*}

\begin{figure}[ht]
\begin{center}
\begin{minipage}[b]{.37\linewidth}
\begin{center}
\epsfig{figure=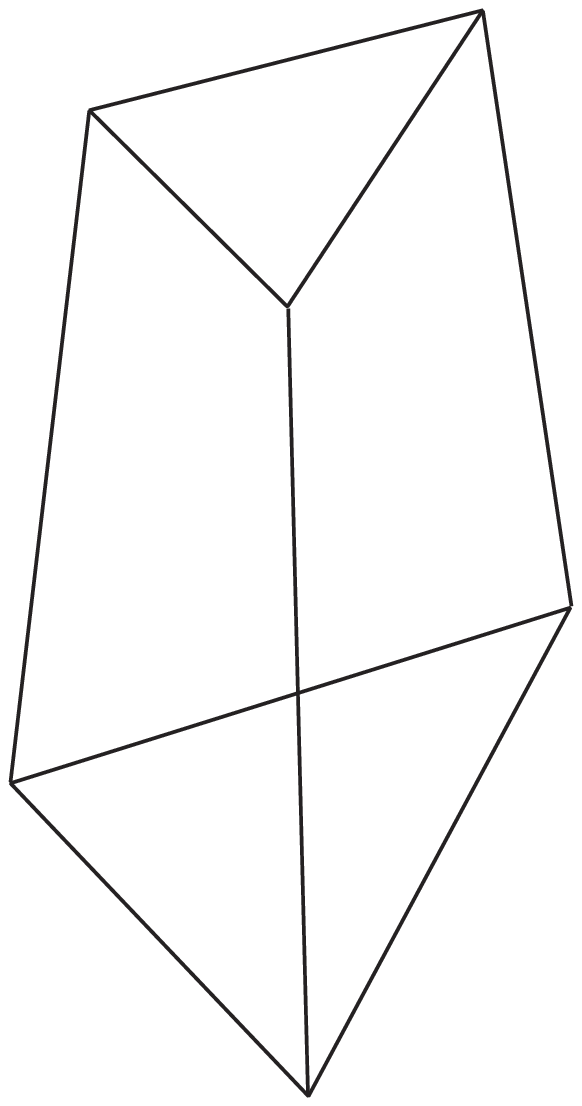,height=3.5cm}
\caption{The image of the moment map: the bottom triangle has side length $m\lambda_1+\lambda_2$, and the top one $\lambda_2$.}\label{fig:polyebdl}
\end{center}
\end{minipage}\hfill
\begin{minipage}[b]{72mm}
\begin{center}
\unitlength0.8mm
\begin{picture}(80,40)
\put(10,5){\circle*{2}} \put(10,35){\circle*{2}} \put(10,5){\line(0,1){30}}
\put(8,1){$\sigma_3$} \put(8,37){$\sigma_1$} \put(-5,18){$\scriptstyle m\lambda_1+\lambda_2$}
\put(70,5){\circle*{2}} \put(10,5){\line(1,0){60}} \put(70,1){$\sigma_6$} \put(38,6){$\lambda_1$}
\put(25,20){\circle*{2}} \put(10,5){\line(1,1){15}} \put(25,20){\line(-1,1){15}} \put(19,20){$\sigma_2$}
\put(17.5,10.5){$\scriptstyle m\lambda_1+\lambda_2$} \put(17.5,28){$\scriptstyle m\lambda_1+\lambda_2$}
\put(55,20){\circle*{2}} \put(25,20){\line(1,0){30}} \put(55,20){\line(1,1){15}} \put(55,20){\line(1,-1){15}}
\put(70,35){\circle*{2}} \put(70,5){\line(0,1){30}} \put(10,35){\line(1,0){60}}
\put(38,21){$\lambda_1$} \put(38,36){$\lambda_1$} \put(57,20){$\sigma_5$}
\put(62.5,13){$\scriptstyle\lambda_2$} \put(62.5,25.5){$\scriptstyle\lambda_2$}
\put(71,18){$\scriptstyle\lambda_2$}
\put(70,37){$\sigma_4$}
\end{picture}
\caption{Invariant 1--dim.\ subvarieties of the variety $\P_{\P^2}(\O(m)\oplus 1)$ and their homology classes.\newline \newline}\label{fig:cod1divex}
\end{center}
\end{minipage}
\end{center}
\end{figure}

From the data above, it is now easy to determine the moment polytope
of $X_\Sigma$ (see figure \ref{fig:polyebdl}, and using Poincar\'e duality,
to determine the homology classes of the invariant one--dimensional irreducible
subvarieties. Note that these invariant subvarieties are mapped to the edges
of the moment polytope --- their respective homology classes are shown in
figure \ref{fig:cod1divex}.

\noindent For our calculations of the Gromov--Witten invariants, it will be convenient to have the following weight
table at hand, giving the weights of the torus action at the different charts (remember that each maximal
cone gives a chart of the toric manifold):

\renewcommand{\arraystretch}{1.3}
\begin{center}
\begin{figure}
\begin{center}
\noindent\begin{tabular}{|l@{$=$}l@{\hspace{0.3cm}}l@{$=$}l|l@{$=$}l@{\hspace{0.3cm}}l@{$=$}l|}\hline
\multicolumn{4}{|c|}{$\sigma_1=\langle v_1, v_3, v_4\rangle$}&
\multicolumn{4}{c|}{$\sigma_4=\langle v_2, v_3, v_4\rangle$}\\\hline
$u_1$&$e_1^*$&$\we_{\sigma_4}^{\sigma_1}$&$\we_{12}-m\we_5$&
$u_1$&$-e_1^*$&$\we_{\sigma_1}^{\sigma_4}$&$\we_{21}+m\we_5$\\
$u_2$&$e_2^*$&$\we_{\sigma_3}^{\sigma_1}$&$\we_3-\we_5$&
$u_2$&$e_2^*$&$\we_{\sigma_6}^{\sigma_4}$&$\we_3-\we_5$\\
$u_3$&$e_3^*$&$\we_{\sigma_2}^{\sigma_1}$&$\we_4-\we_5$&
$u_3$&$e_3^*$&$\we_{\sigma_5}^{\sigma_4}$&$\we_4-\we_5$\\\hline\hline
\multicolumn{4}{|c|}{$\sigma_2=\langle v_1, v_3, v_5\rangle$}&
\multicolumn{4}{c|}{$\sigma_5=\langle v_2, v_3, v_5\rangle$}\\\hline
$u_1$&$e_1^*-me_3^*$&$\we_{\sigma_5}^{\sigma_2}$&$\we_{12}-m\we_4$&
$u_1$&$-e_1^*+me_3^*$&$\we_{\sigma_2}^{\sigma_5}$&$\we_{21}+m\we_4$\\
$u_2$&$e_2^*-e_3^*$&$\we_{\sigma_3}^{\sigma_2}$&$\we_3-\we_4$&
$u_2$&$e_2^*-e_3^*$&$\we_{\sigma_6}^{\sigma_5}$&$\we_3-\we_4$\\
$u_3$&$-e_3^*$&$\we_{\sigma_1}^{\sigma_2}$&$\we_5-\we_4$&
$u_3$&$-e_3^*$&$\we_{\sigma_4}^{\sigma_5}$&$\we_5-\we_4$\\\hline\hline
\multicolumn{4}{|c|}{$\sigma_3=\langle v_1, v_4, v_5\rangle$}&
\multicolumn{4}{c|}{$\sigma_6=\langle v_2, v_4, v_5\rangle$}\\\hline
$u_1$&$e_1^*-me_2^*$&$\we_{\sigma_6}^{\sigma_3}$&$\we_{12}-m\we_3$&
$u_1$&$-e_1^*+me_2^*$&$\we_{\sigma_3}^{\sigma_6}$&$\we_{21}+m\we_3$\\
$u_2$&$e_3^*-e_2^*$&$\we_{\sigma_2}^{\sigma_3}$&$\we_4-\we_3$&
$u_2$&$e_3^*-e_2^*$&$\we_{\sigma_5}^{\sigma_6}$&$\we_4-\we_3$\\
$u_3$&$-e_2^*$&$\we_{\sigma_1}^{\sigma_3}$&$\we_5-\we_3$&
$u_3$&$-e_2^*$&$\we_{\sigma_4}^{\sigma_6}$&$\we_5-\we_3$\\\hline
\end{tabular}
\caption{The weight table for $\P_{P^2}(\O(m)\oplus 1)$. Substitute $\we_{12}$ 
by $\we_1-\we_2$ and $\we_{21}$ by $\we_2-\we_1$.}
\end{center}
\end{figure}
\end{center}

\renewcommand{\arraystretch}{1}

\subsection{The Gromov--Witten invariants of $\P(\O_{\P^2}(2) \oplus 1)$}

We will now compute all genus--$0$ three--point Gromov--Witten invariants 
of the Fano threefold $\P(\O_{\P^2}(2) \oplus 1)$. The Fano property of
this manifold implies that there is only a finite number of non--zero
invariants. 

We will give full details for the computations of some of the invariants,
but in some cases, due to the length of the computations, we will just give 
the different graph types one need to consider and the result. However, even
those computations that are omitted here are straightforward.

Remember from above that the degree--$2$ homology classes $\lambda_1=(1,1,0,0,0)$ and 
$\lambda_2=(0,-2,1,1,1)$ generate the effective cone.
Their pairing with the Chern classes $c_1(\vt)$ is two respectively one.
Any degree--$2$ homology class $A\in H_2(\vt, \Z)$ having non--zero Gromov--Witten
invariants satisfies $\langle c_1(\vt), A\rangle \le 6$; otherwise 
the virtual dimension of the moduli stack were negative. 

However, it is easy to see that for homology classes $A$ such that $\langle c_1(\vt, A\rangle)=6$,
all Gromov--Witten invariants are equal to zero: let $A$ be such a class, then the only
possibly non--trivial GW invariant would be
\[ \Phi^A(Z_1Z_3Z_4, Z_1Z_3Z_5, Z_1Z_4Z_5). \]
However, a graph $\mygraph$ such that the integral of these
classes over the corresponding fixed point moduli space
$\cM_\mygraph$ is non--zero has to contain the nodes $\sigma_1$,
$\sigma_2$ and $\sigma_3$. By looking at figure
\ref{fig:cod1divex} one immediately sees that such a graph
$\mygraph$ would have to have homology class $A_\mygraph$ with
$\langle c_1(\vt), A_\mygraph\rangle \geq 8$ (in this case
$A_\mygraph=3\lambda_1+2\lambda_2$). Hence, all non--zero
Gromov--Witten invariants of $\vt$ have $\langle c_1(\vt), A
\rangle \leq 5$.

We can equally exclude all classes $A=a_1\lambda_1 + a_2\lambda_2$ with $a_2>3$. For if
$\langle c_1(\vt), A \rangle \leq 5$ and $a_2>3$ we had $a_1=0$. So let us consider the
Gromov--Witten invariant
\[ \Phi^{a_2\lambda_2}(\alpha_1, \alpha_2, \alpha_3).\]
Since $\dim \vt = \langle c_1(\vt, a_2\lambda_2\rangle) = 3+
a_2>6$, at least one of the $\alpha_i$'s is of degree six, say
$\alpha_1=Z_1Z_3Z_4$. But there is no graph $\mygraph$ of homology
class a multiple of $\lambda_2$ that contains $\sigma_1$.

So we only have to compute the 3--point genus--0 Gromov--Witten invariants for the
following classes:
\[ \lambda_2, \lambda_1, 2\lambda_2, 3\lambda_2, \lambda_1+\lambda_2, 
2\lambda_1, \lambda_1+2\lambda_2, \lambda_1+3\lambda_2, 2\lambda_1+\lambda_2.
\]
\begin{list}{}{}
\item[\bf $\lambda_2$--invariants:] The graph type for these invariants live on the ``upper''
triangle in the moment polytope, \ie. they can only contain the vertices $\sigma_4$, $\sigma_5$
and $\sigma_6$.
\begin{list}{}{}
\item[\bf $\Phi^{\lambda_2}(Z_1, Z_i, Z_jZ_k)=0$.] This is because the weight of $Z_1$ on the 
upper triangle is zero.
\item[\bf $\Phi^{\lambda_2}(Z_i, Z_j, Z_1Z_k)=0$.] Again, the weight of $Z_1$ on the upper
triangle is zero, hence the same applies to $Z_1Z_k$.
\item[\bf $\Phi^{\lambda_2}(Z_3, Z_3,Z_4Z_5)=-1$.] In this case we get the following two possible graph types:
\[ \gtype[a] = \bbull{\sigma_4} \ml \bbull{\sigma_6}, \qquad 
\gtype[b] = \bbull{\sigma_5} \ml \bbull{\sigma_6}.
\]
For both graph types, the automorphism group is trivial.
The $S_{\gtype}$ and $T_{\gtype}$ terms are as follows:
\begin{align*}
S_{\gtype[a]}&= 1^2 \cdot \frac{(\we_4-\we_3)(\we_5-\we_3)}{(\we_5-\we_3)} = (\we_4-\we_3)\\
S_{\gtype[b]}&= 1^2 \cdot \frac{(\we_4-\we_3)(\we_5-\we_3)}{(\we_4-\we_3)} = (\we_5-\we_3)\\
T_{\gtype[a]}&= \frac{(-1) \cdot 1^2}{1^2 (\we_3-\we_5)^2} \frac{(\we_{21}+\we_3\we_5)}{(\we_4-\we_5)(\we_4-\we_3)}
\cdot \frac{(\we_3-\we_5)^2}{(\we_3-\we_5)} \cdot \frac{(\we_5-\we_3)^2}{(\we_5-\we_3)}\\
&= \frac{(\we_{21} +\we_3 +\we_5)}{(\we_4-\we_5)(\we_4-\we_3)}\\
T_{\gtype[b]}&= \frac{(-1) \cdot 1^2}{1^2(\we_3-\we_4)^2} \frac{(-\we{12}+\we_3+\we_4)}{(\we_5-\we_4)(\we_4-\we_3)}
\cdot \frac{(\we_3-\we_4)^2}{(\we_3-\we_4)} \cdot \frac{(\we_4-\we_3)^2}{(\we_4-\we_3)}\\
&=\frac{(\we_{21}+\we_3+\we_4)}{(\we_5-\we_4)(\we_5-\we_3)}
\end{align*}
\end{list}
\item[\bf $\lambda_1$--invariants:] The graph types of this homology class live on the side edges,
\ie. either connect $\sigma_1$ and $\sigma_4$, $\sigma_2$ and $\sigma_5$, or $\sigma_3$ and
$\sigma_6$.
\begin{list}{}{}
\item[\bf $\Phi^{\lambda_1}(Z_i, Z_3, Z_2Z_4Z_5)=0$.] For the weight of $Z_2Z_4Z_5$ not to be zero,
the graph has to contain $\sigma_6$. However, the weight of $Z_3$ at $\sigma_3$ and $\sigma_6$
is zero.
\item[\bf $\Phi^{\lambda_1}(Z_i, Z_jZ_3, Z_4Z_5)=0$.] By a similar argument.
\item[\bf $\Phi^{\lambda_1}(Z_1, Z_1, Z_2Z_3Z_4)=1$.] The only graph type is
\[\gtype = \bbull{\sigma_1} \ml \bbull{\sigma_4}\]
Its automorphism group is trivial, and the two other terms are given by:
\begin{align*}
S_{\gtype}&=1^2\cdot(\we_3-\we_5)(\we_4-\we_5)\\
T_{\gtype}&=\frac{(-1)}{(\we_{12}-2\we_5)} \frac{1}{(\we_3-\we_5)(\we_4-\we_5)}
\frac{(\we_{12}-2\we_5)^2}{(\we_{12}-2\we_5)} \frac{(\we_{21}+2\we_5)^2}{(\we_{21}+2\we_5)}.
\end{align*}
\item[$\Phi^{\lambda_1}(Z_1, Z_1Z_4, Z_1Z_3)=1$] We get the same graph type and
the same terms $\AA_{\gtype}$, $S_{\gtype}$ and $T_{\gtype}$ as in the previous case.
\end{list} 
\item[\bf $2\lambda_2$--invariants:] Again these invariants have to live on the
upper $\sigma_4$--$\sigma_5$--$\sigma_6$ triangle of the polytope, so the weight of
any class containing $Z_1$ will cause an invariant to be zero.
\begin{list}{}{}
\item[\bf $\Phi^{2\lambda_2}(Z_i, Z_j, pt)=0$.] The class $\lambda_2$ has negative self
intersection with the divisor $D_2$, hence any holomorphic curve has to lie in $D_2$,
so it can not pass trough a point in general position.
\item[\bf $\Phi^{2\lambda_2}(Z_1, Z_iZ_j, Z_kZ_l)=0$.] Simply because of $Z_1$.
\item[\bf $\Phi^{2\lambda_2}(Z_i, Z_1Z_j, Z_kZ_l)=0$.] Again because of $Z_1Z_j$.
\item[\bf $\Phi^{2\lambda_2}(Z_3, Z_3Z_4,Z_4Z_5)=-2$.] We will just give the graph types
one has to consider and leave the actual computations to the reader. The graph types
are:
\begin{alignat*}{3}
\gtype[a]&=\bbull{\sigma_4} \ml[2] \bbull{\sigma_6}&
\gtype[b]&=\bbull{\sigma_4} \ml \bbull{\sigma_6} \ml \bbull{\sigma_5}\\
\gtype[aa]&=\bbull{\sigma_4} \ml \bbull{\sigma_6} \ml \bbull{\sigma_4}&
\gtype[c]&=\bbull{\sigma_5} \ml \bbull{\sigma_4} \ml \bbull{\sigma_6} \\
\gtype[ab]&=\bbull{\sigma_6} \ml \bbull{\sigma_4} \ml \bbull{\sigma_6}&
\gtype[d]&=\bbull{\sigma_4} \ml \bbull{\sigma_5} \ml \bbull{\sigma_6}
\end{alignat*}
\end{list}
\item[\bf $3\lambda_2$--invariants:] Again, all graphs have to live on the upper
triangle of the polytope. Therefore we get:
\begin{list}{}{}
\item[\bf $\Phi^{3\lambda_2}(Z_i,Z_jZ_k,Z_1Z_lZ_p)=0$.] Because of the $Z_1$ in the third class.
\item[\bf $\Phi^{3\lambda_2}(Z_1Z_i, Z_jZ_k, Z_lZ_p)=0$.] Because of the $Z_1$ in the first class.
\item[\bf $\Phi^{3\lambda_2}(Z_3Z_4,Z_4Z_5,Z_3Z_5)=-4$.] Again, we only give the different 
graph types here --- they are:
\begin{alignat*}{2}
\gtype[a]&=\bbull{\sigma_5} \ml[2] \bbull{\sigma_4} \ml \bbull{\sigma_6} &
\gtype[b]&=\bbull{\sigma_5} \ml \bbull{\sigma_4} \ml[2] \bbull{\sigma_6}\\
\gtype[aa]&=\bbull{\sigma_4} \mls \bbull{\sigma_5} \mls \bbull{\sigma_4} \mls \bbull{\sigma_6}&
\gtype[ba]&=\bbull{\sigma_4} \mls \bbull{\sigma_6} \mls \bbull{\sigma_5} \mls \bbull{\sigma_4}\\
\gtype[ab]&=\fork{\sigma_5}{\sigma_4}{\sigma_6}{\sigma_5}&
\gtype[bb]&=\fork{\sigma_5}{\sigma_4}{\sigma_6}{\sigma_6}\\
\gtype[c]&=\bbull{\sigma_4} \ml[2] \bbull{\sigma_5} \ml \bbull{\sigma_6}&
\gtype[d]&=\bbull{\sigma_4} \ml \bbull{\sigma_5} \ml[2] \bbull{\sigma_6}\\
\gtype[ca]&=\bbull{\sigma_5} \mls \bbull{\sigma_4} \mls \bbull{\sigma_5} \mls \bbull{\sigma_6}&
\gtype[da]&=\bbull{\sigma_4} \mls \bbull{\sigma_5} \mls \bbull{\sigma_6} \mls \bbull{\sigma_5}\\
\gtype[cb]&=\fork{\sigma_4}{\sigma_5}{\sigma_6}{\sigma_4}&
\gtype[db]&=\fork{\sigma_4}{\sigma_5}{\sigma_6}{\sigma_6}\\
\gtype[e]&=\bbull{\sigma_4} \ml[2] \bbull{\sigma_6} \ml \bbull{\sigma_5}&
\gtype[f]&=\bbull{\sigma_4} \ml \bbull{\sigma_6} \ml[2] \bbull{\sigma_5}\\
\gtype[ea]&=\bbull{\sigma_6} \mls \bbull{\sigma_4} \mls \bbull{\sigma_6} \mls \bbull{\sigma_5}&
\gtype[fa]&=\bbull{\sigma_4} \mls \bbull{\sigma_6} \mls \bbull{\sigma_5} \mls \bbull{\sigma_6}\\
\gtype[eb]&=\fork{\sigma_4}{\sigma_6}{\sigma_5}{\sigma_4}&
\gtype[fb]&=\fork{\sigma_4}{\sigma_5}{\sigma_5}{\sigma_5}\\
\gtype[ga]&=\bbull{\sigma_4} \mls \bbull{\sigma_5} \mls \bbull{\sigma_6} \mls \bbull{\sigma_4}&
\gtype[gb]&=\bbull{\sigma_5} \mls \bbull{\sigma_4} \mls \bbull{\sigma_6} \mls \bbull{\sigma_5}\\
\gtype[gc]&=\bbull{\sigma_6} \mls \bbull{\sigma_4} \mls \bbull{\sigma_5} \mls \bbull{\sigma_6}
\end{alignat*}
\end{list} 
\item[\bf $(\lambda_1+\lambda_2)$--invariants:] The graph types have one component in the upper
triangle and one on a side edge --- they can not have an edge on the lower triangle
$\sigma_1$--$\sigma_2$--$\sigma_3$.
\begin{list}{}{}
\item[\bf $Phi^{\lambda_1+\lambda_2}(Z_i, Z_1Z_5,Z_1Z_3Z_4)=0$.] The third marked point has to be
at $\sigma_1$, and the second at $\sigma_2$ or $\sigma_3$ which would require an edge on the
lower triangle.
\item[\bf $\Phi^{\lambda_1+\lambda_2}(Z_3Z_4,Z_4Z_5,Z_3Z_5)=0$.] Needs two edges in the upper
triangle.
\item[\bf $\Phi^{\lambda_1+\lambda_2}(Z_1Z_3,Z_1Z_4,Z_1Z_5)=0$.] Needs an edge in the lower
triangle.
\item[\bf $\Phi^{\lambda_1+\lambda_2}(Z_1,Z_4Z_5,Z_1Z_3Z_4)=1$.] The only graph type to consider is
\[ \gtype = \bbull{\sigma_1} \ml \bbull{\sigma_4} \ml \bbull{\sigma_6}
\]
The automorphism group is trivial, $S_{\gtype}=(\we_4-\we_3)(\we_3-\we_5)(\we_4-\we_5)$, and
$T_{\gtype}$ the invers of $S_{\gtype}$.
\item[\bf $\Phi^{\lambda_1+\lambda_2}(Z_3,Z_4Z_5,Z_1Z_3Z_4)=1$.] Same graph and terms
as in the previous case.
\item[\bf $\Phi^{\lambda_1+\lambda_2}(Z_1Z_3,Z_4Z_5,Z_3Z_4)=1$.] Again the same.
\item[\bf $\Phi^{\lambda_1+\lambda_2}(Z_1Z_3,Z_4Z_5,Z_1Z_4)=1$.] Again the same.
\end{list}
\item[\bf $2\lambda_1$--invariants are all zero:]
The virtual dimension of the corresponding moduli space is seven, so for the
invariants $\Phi^{2\lambda_1}(\alpha_1,\alpha_2,\alpha_3)$ we have two cases for
the classes $\alpha_i$:
\begin{itemize}
\item $\deg \alpha_1=\deg \alpha_2=6$ and $\deg \alpha_3=2$\newline
In this cases we can set $\alpha_1=Z_1Z_3Z_4$ and
$\alpha_2=Z_1Z_4Z_5$. There is obviously no graph $\mygraph$ in
the homology class $2\lambda_1$ containing both $\sigma_1$ and
$\sigma_3$.
\item $\deg \alpha_1=6$ and $\deg \alpha_2=\deg \alpha_3=4$\newline
We can set $\alpha_1=Z_1Z_3Z_4$. For $\alpha_2$ we have to
choices: $\alpha_2=Z_4Z_5$ or $\alpha_2=Z_1Z_5$. Again, there is
no graph $\mygraph$ with homology class $2\lambda_1$ that contains
the necessary nodes ($\sigma_1$ and one of the following:
$\sigma_3$ or $\sigma_6$ respectively $\sigma_2$ or $\sigma_3$).
\end{itemize}
\item[\bf $(\lambda_1+2\lambda_2)$--invariants:]
Since the homology class $A=\lambda_1+2\lambda_2$ contains only
one $\lambda_1$, all graphs $\mygraph$ in this homology class
contain exactly one of the following nodes: $\sigma_1$, $\sigma_2$
and $\sigma_3$. Therefore the following invariants are all zero:
\begin{list}{}{}
\item[$\Phi^{\lambda_1+2\lambda_2}(Z_1, Z_1Z_3Z_4, Z_1Z_4Z_5)=0$.]
\item[$\Phi^{\lambda_1+2\lambda_2}(Z_3, Z_1Z_3Z_4, Z_1Z_4Z_5)=0$.]
\item[$\Phi^{\lambda_1+2\lambda_2}(Z_1Z_3, Z_3Z_4, Z_1Z_4Z_5)=0$.]
\item[$\Phi^{\lambda_1+2\lambda_2}(Z_1Z_3, Z_1Z_3, Z_1Z_4Z_5)=0$.]
\end{list}
Only Gromov--Witten invariant in this class remains to be computed:
\begin{list}{}{}
\item[\bf $\Phi^{\lambda_1+2\lambda_2}(Z_4Z_5,Z_3Z_5, Z_1Z_3Z_4)=1$.]
Here we have to consider three graph types:

\noindent\begin{tabular}{cc}
$\bbull{\sigma_1} \mls \bbull{\sigma_4} \mls \bbull{\sigma_5} \mls \bbull{\sigma_6}$&
\\
$\bbull{\sigma_1} \mls \bbull{\sigma_4} \mls \bbull{\sigma_6} \mls \bbull{\sigma_5}$&
\raisebox{0ex}[-5ex]{$\fork{\sigma_1}{\sigma_4}{\sigma_6}{\sigma_5}$}
\end{tabular}

\noindent The details of the computations are left to the reader.
\end{list}
\item[\bf $(\lambda_1+3\lambda_2)$--invariants are all zero:]
The virtual dimension of the corres\-pon\-ding moduli space is eight,
so we can set $\alpha_1=Z_1Z_3Z_4$ and $\alpha_2=Z_1Z_4Z_5$. A
graph $\mygraph$ that could give a non--zero integral on
$\cM_\mygraph$ had to contain $\sigma_1$ and $\sigma_3$, which is
impossible since the class $A=\lambda_1+3\lambda_2$ contains only
one $\lambda_1$.
\item[\bf $(2\lambda_1+\lambda_2)$--invariants:] Graph types for this class
have either on edge in the ``lower triangle'', or two side edges and one
on the ``upper triangle''.
\begin{list}{}{}
\item [\bf $\Phi^{2\lambda_1+\lambda_2}(Z_3Z_5, Z_1Z_3Z_4, Z_1Z_4Z_5)=0$.] Would need
two edges in the triangles.
\item[\bf $\Phi^{2\lambda_1+\lambda_2}(Z_2^2, Z_1Z_3Z_4, Z_1Z_4Z_5)=-2$.]
The only graph type we need to consider is:
\[ \gtype = \bbull{\sigma_1} \mls \bbull{\sigma_4} \mls \bbull{\sigma_6}
\mls \bbull{\sigma_3} \]
By the usual computations we get $|\AA_{\gtype}|=1$ and
\begin{align*}
S_{\gtype}&= -2 (\we_{12}-\we_3-\we_5)(\we_3-\we_5)^2(\we_4-\we_5)(\we_4-\we_3)\\
T_{\gtype}&= -\left( (-\we_12+\we_3+\we_5)(\we_3-\we_5)^2(\we_4-\we_5)(\we_4-\we_3)\right)^{-1}.
\end{align*}
\item[\bf $\Phi^{2\lambda_1+\lambda_2}(Z_1Z_3, Z_1Z_3Z_4, Z_1Z_4Z_5)=1$.]
Here we need to consider the following two graph types for which we leave
the computations to the reader:
\[ \bbull{\sigma_1} \ml \bbull{\sigma_3} \qquad
\bbull{\sigma_1} \mls \bbull{\sigma_4} \mls \bbull{\sigma_6} \mls \bbull{\sigma_3}
\]
\end{list}
\end{list}

\subsection{The quantum cohomology ring of $\P(\O_{\P^2}(2) \oplus 1)$}

For completeness, we will also compute the quantum cohomology ring of this manifold,
although it has already been known thanks to Givental's work \cite{giv97} where he uses
techniques from mirror symmetry to compute the quantum cohomology ring for Fano toric
varieties, obtaining the same formula postulated by Batyrev in \cite{bat93}.

Since the usual cohomology ring is given by
\[ H^*(\vt, \Q) = \Q[Z_1, Z_2, Z_3, Z_4, Z_5]/_{\left\langle 
\displaystyle \begin{array}{l} Z_1-Z_2-2Z_5, Z_3-Z_4, Z_3-Z_5,\\
Z_1Z_2, Z_3Z_4Z_5 \end{array} \right\rangle,}\]
it suffices to calculate the quantum products $Z_1\star Z_2$ and $Z_3\star Z_3\star Z_3$
to find a representation of the quantum cohomology ring. Remember that given the
Gromov--Witten invariants, the quantum product satisfies the following determining
equalities:
\begin{equation}
\begin{split}
\langle (\alpha\star\beta)_\lambda\cap\gamma, [X]\rangle=&
\Phi^{X,\lambda}_{3,0}(\alpha,\beta,\gamma)\\
\alpha\star\beta=&\sum_{\lambda\in H^2(X,\Z)}(\alpha\star\beta)_\lambda q^\lambda,
\end{split}
\end{equation}
where $\alpha,\beta,\gamma\in H^*(X,\Q)$ are cohomology classes of the manifold $X$.
Thus if $\theta_1, \ldots, \theta_r$ is a basis of $H^*(X,\Q)$ and
$\vartheta_1, \ldots, \vartheta_r$ its
dual basis with respect to the cap product plus integration, we obtain
\[ (\alpha\star\beta)_\lambda = \sum_{i=1}^r \Phi^{X,\lambda}_{3,0}(\alpha,\beta, \theta_i)
\vartheta_i.\]
Now, for our particular example $X=X_\Sigma$, we will take the following basis with
its dual basis:
\[\begin{array}{l|c|c|c|c|c|c}
\mbox{basis}&1&Z_1&Z_3&Z_1Z_3&Z_3^2&Z_1Z_3^2\\\hline
\mbox{dual basis}&Z_1Z_3^2&Z_3^2&Z_1Z_3-2Z_3^2&Z_3&Z_1-2Z_3&1.
\end{array}
\]
So, for the first product $Z_1\star Z_2$ we obtain
\begin{eqnarray*}
Z_1\star Z_2&=&Z_1 \star Z_1 - 2Z_1\star Z_3\\
&=&(Z_1\star Z_1)_{\lambda_1} q^{\lambda_1}\\
&=& q^{\lambda_1}, \mbox{ since}\\
Z_1\star Z_1&=&Z_1^2 + q^{\lambda_1}\\
&=&\underbrace{Z_1Z_2}_{=0} + 2Z_1Z_3 + q^{\lambda_1}.
\end{eqnarray*}
For $Z_3^{\star3}$ we first calculate the quantum square product of $Z_3\star Z_3$:
\begin{eqnarray*}
Z_3\star Z_3&=&Z_3\cup Z_3 +
\underbrace{(Z_3\star Z_3)_{\lambda_2}}_{\Phi^{\lambda_2}(Z_3, Z_3, Z_3^2)(Z_1-2Z_3)}q^{\lambda_2}
+\underbrace{(Z_3\star Z_3)_{\lambda_1}}_{\Phi^{\lambda_1}(Z_3, Z_3, pt)=0}q^{\lambda_1}+\\
&&+ \underbrace{(Z_3\star Z_3)_{2\lambda_2}}_{\Phi^{2\lambda_2}(Z_3, Z_3, pt)=0}q^{2\lambda_2}\\
&=&Z_3^2 - (Z_1-2Z_3) q^{\lambda_2}.
\end{eqnarray*}

So we will also need the products $Z_1\star Z_3$ and $Z_3^2\star Z_3$. For the first product
notice that all Gromov--Witten invariants $\Phi^A(Z_1, Z_3, \alpha)$ are zero for $A\neq 0$.
So
\[ Z_1 \star Z_3 = Z_1\cup Z_3.\]
For the second product we obtain:
\begin{eqnarray*}
Z_3^2\star Z_3&=& Z_3^3 + (Z_3^2\star Z_3)_{\lambda_1}q^{\lambda_1} +
(Z_3^2\star Z_3)_{\lambda_2}q^{\lambda_2}+\\
&&+ (Z_3^2\star Z_3)_{2\lambda_2}q^{2\lambda_2} +(Z_3^2\star Z_3)_{3\lambda_2}q^{3\lambda_2}+
(Z_3^2\star Z_3)_{\lambda_1+\lambda_2}q^{\lambda_1+\lambda_2}\\
&=& -(Z_1Z_3-2Z_3^2)q^{\lambda_2}-2(Z_1-2Z_3)q^{2\lambda_2} + q^{\lambda_1+\lambda_2}
\end{eqnarray*}

\noindent Summing all up we thus get the following expression for $Z^{\star3}$:
\begin{eqnarray*}
Z_3\star Z_3\star Z_3&=&Z_3^2\star Z_3 + Z_1\star Z_3q^{\lambda_2}\\
&=& (-2Z_1Z_3+4Z_3^2)q^{\lambda_2}-4(Z_1-2Z_3)q^{2\lambda_2}+q^{\lambda_1+\lambda_2}\\
&=& (-2Z_1\star Z_3 +4 Z_3\star Z_3 + q^{\lambda_1})q^{\lambda_2}\\
&=&  Z_2\star Z_2 q^{\lambda_2}.
\end{eqnarray*}

Thus we obtain  the same result as Givental in \cite{giv97}, the quantum cohomology
ring defined by Batyrev in \cite{bat93}:

\begin{prop}
The quantum cohomology ring of the smooth toric Fano variety $\P_{\C\P^2}(\O(2)\oplus 1)$
is equal to
\begin{align*}
\lefteqn{QH^*(\P_{\C\P^2}(\O(2)\oplus 1), \C)=}\\ &= 
\C[Z_1, Z_2,Z_3,Z_4, Z_5, q_1, q_2]/_{\left\langle
\begin{array}{l}Z_1-Z_2-2Z_5, Z_3-Z_4, Z_3-Z_5,\\ Z_1Z_2-q_1, Z_3Z_4Z_5-Z_2^2q_2\end{array}\right\rangle}\\
&= \C[Z_2, Z_3, q_1, q_2]/{\langle Z_2^2+2Z_2Z_3-q_1, Z_3^3-Z_2^2q_2\rangle}.
\end{align*}
\end{prop}

\subsection{On the quantum cohomology ring of the $3$--fold $\P(\O_{\P^2}(3) \oplus 1)$}

In this section, we will give the details of the computation of some
Gromov--Witten invariants of the (non--Fano) threefold
$\P(\O_{\P^2}(3) \oplus 1)$. As a corollary, we will show that --- in
the contrary to the previous example in the Fano case --- the
quantum cohomology ring of this variety defined by the Gromov--Witten
invariants of \cite{rt95,bf97,beh97,lt96,lt98a,sie96} does not
coincide with the ring defined formally in \cite{bat93}.

\subsubsection{Some Gromov--Witten invariants of $\P(\O_{\P^2}(3) \oplus 1)$}

In this section we will compute the genus--$0$ Gromov--Witten invariants
$\Phi^{p\lambda_2}(Z_3, Z_3, Z_3)$ and $\Phi^{p\lambda_2}(Z_2, Z_2, Z_2)$
for $p\leq 2$. Note first of all that since $Z_3^3=0$ in the cohomology
ring,
\[ \Phi^0(Z_3, Z_3, Z_3)=0.\]

\paragraph{The invariant $\Phi^{\lambda_2}(Z_3, Z_3, Z_3)$}

The $\lambda_2$--graphs have to live on the triangle in the moment
polytope spanned by $\sigma_4$, $\sigma_5$ and $\sigma_6$. Hence,
there are three $\lambda_2$--graphs:
\[ \gtype[1] = \bbull{\sigma_4} \ml[1] \bbull{\sigma_5} \qquad
\gtype[2] = \bbull{\sigma_4} \ml[1] \bbull{\sigma_6} \quad
\gtype[3] = \bbull{\sigma_5} \ml[1] \bbull{\sigma_6}
\]

First of all we note that none of the graphs has any non--trivial automorphisms.
Since the only edge of any of the graphs has multiplicity one we get
\[ |\AA_{\gtype[i]} | = 1. \]

We will now compute the $S_\Gamma$--terms for the three graphs. In fact,
there is only one such for each of the graphs, since the cohomology classes
are all equal to $Z_3$. By applying formula \ref{eq:sgammadef} we obtain:
\begin{align*}
S_{\gtype[1]}&= \left( \frac{1\cdot (\we_3-\we_5)}{(\we_4-\we_5)} + \frac{1\cdot (\we_3-\we_4)}{(\we_5-\we_4)}\right)^3\\
&=1^3=1.
\end{align*}
\begin{align*}
S_{\gtype[2]}&= \left( \frac{1\cdot (\we_3-\we_5)}{(\we_3-\we_5)} + 0 \right)^3 = 1\\
S_{\gtype[3]}&= \left( \frac{1\cdot (\we_3-\we_4)}{(\we_3-\we_4)} + 0 \right)^3 = 1.
\end{align*}

Let us now compute the $T_\Gamma$--terms. Remember, that $\lambda_e=\lambda_2$.
Hence we get:
\begin{align*}
T_{\gtype[1]}&= \left( \frac{(-1) \cdot 1^2}{1^2 \cdot (\we_4-\we_5)^2} \cdot
\frac{(-\we_1+\we_2+\we_4+2\we_5)(-\we_1+\we_2+2\we_4+\we_5)}{(\we_3-\we_5)(\we_3-\we_4)} \right)\cdot\\
&\phantom{=}\cdot \frac{1\cdot\left( \frac{1}{(\we_4-\we_5)} \right)^{-2}}{(\we_4-\we_5)}
\cdot \frac{1\cdot\left( \frac{1}{(\we_5-\we_4)} \right)^{-2}}{(\we_5-\we_4)}\\
&= \frac{(-\we_1+\we_2+\we_4+2\we_5)(-\we_1+\we_2+2\we_4+\we_5)}{(\we_3-\we_5)(\we_3-\we_4)}\\
T_{\gtype[2]}&= \left( \frac{(-1) \cdot 1^2}{1^2 \cdot (\we_3-\we_5)^2} \cdot
\frac{(-\we_1+\we_2+\we_3+2\we_5)(-\we_1+\we_2+2\we_3+\we_5)}{(\we_4-\we_5)(\we_4-\we_3)} \right)\cdot\\
&\phantom{=}\cdot \frac{1\cdot\left( \frac{1}{(\we_3-\we_5)} \right)^{-2}}{(\we_3-\we_5)}
\cdot \frac{1\cdot\left( \frac{1}{(\we_5-\we_3)} \right)^{-2}}{(\we_5-\we_3)}\\
&= \frac{(-\we_1+\we_2+\we_3+2\we_5)(-\we_1+\we_2+2\we_3+\we_5)}{(\we_4-\we_5)(\we_4-\we_3)}\\
T_{\gtype[3]}&= \left( \frac{(-1) \cdot 1^2}{1^2 \cdot (\we_3-\we_3)^2} \cdot
\frac{(-\we_1+\we_2+\we_3+2\we_4)(-\we_1+\we_2+2\we_3+\we_4)}{(\we_5-\we_4)(\we_5-\we_3)} \right)\cdot\\
&\phantom{=}\cdot \frac{1\cdot\left( \frac{1}{(\we_3-\we_4)} \right)^{-2}}{(\we_3-\we_4)}
\cdot \frac{1\cdot\left( \frac{1}{(\we_4-\we_3)} \right)^{-2}}{(\we_4-\we_3)}\\
&= \frac{(-\we_1+\we_2+\we_3+2\we_4)(-\we_1+\we_2+2\we_3+\we_4)}{(\we_5-\we_4)(\we_5-\we_3)}
\end{align*}

Therefore computing the sum $\sum_i (S_{\Gamma_i}\cdot T_{\Gamma_i})/|\AA_{\Gamma_i}|$ (for example
with the Maple package) yields
\[ \Phi^{\lambda_2}(Z_3, Z_3, Z_3) = 3.\]

\paragraph{The invariant $\Phi^{2\lambda_2}(Z_3, Z_3, Z_3)$}

As for the $\lambda_2$--invariant, this invariant has to live in the triangle with
corners $\sigma_4$, $\sigma_5$ and $\sigma_6$. In this case, we get twelve different graph
types:
\begin{xalignat*}{3}
\gtype[a]&=\bbull{\sigma_4} \mls[2] \bbull{\sigma_5}&
\gtype[aa]&=\bbull{\sigma_4} \mls \bbull{\sigma_5} \mls \bbull{\sigma_4}&
\gtype[ab]&=\bbull{\sigma_5} \mls \bbull{\sigma_4} \mls \bbull{\sigma_5}\\
\gtype[b]&=\bbull{\sigma_4} \mls[2] \bbull{\sigma_6}&
\gtype[ba]&=\bbull{\sigma_4} \mls \bbull{\sigma_6} \mls \bbull{\sigma_4}&
\gtype[bb]&=\bbull{\sigma_6} \mls \bbull{\sigma_4} \mls \bbull{\sigma_6}\\
\gtype[c]&=\bbull{\sigma_5} \mls[2] \bbull{\sigma_6}&
\gtype[ca]&=\bbull{\sigma_5} \mls \bbull{\sigma_6} \mls \bbull{\sigma_5}&
\gtype[cb]&=\bbull{\sigma_6} \mls \bbull{\sigma_5} \mls \bbull{\sigma_6}\\
&&\gtype[d]&=\bbull{\sigma_4} \mls \bbull{\sigma_5} \mls \bbull{\sigma_6}\\
&&\gtype[e]&=\bbull{\sigma_4} \mls \bbull{\sigma_6} \mls \bbull{\sigma_5}\\
&&\gtype[f]&=\bbull{\sigma_5} \mls \bbull{\sigma_4} \mls \bbull{\sigma_6}
\end{xalignat*}

The order of the automorphism groups $\AA_{gtype[i]}$ is as follows:

\medskip

\begin{center}
\begin{tabular}{|l||c|c|c|}\hline
&&aa, ba, ca &\\
&a, b, c&ab, bb, cb&d, e, f\\\hline\hline
$|\Aut_{\gtype}|$&1&2&1\\\hline
$|\bigoplus \Z_{d_e}|$&2&1&1\\\hline\hline
$|\AA_{\gtype}|$&2&2&1\\\hline
\end{tabular}
\end{center}

\medskip

We will now compute the $S_\Gamma$--terms. Note that the images on the moment
polytope of the graphs $\gtype[a]$, $\gtype[aa]$ and $\gtype[ab]$ are all
the same, and hence there $S_{\gtype}$--terms coincide. The same applies to the
graphs $\gtype[b]$, $\gtype[ba]$ and $\gtype[bb]$, respectively
$\gtype[c]$, $\gtype[ca]$ and $\gtype[cb]$. Hence we obtain:
\[
S_{\gtype[a]}=S_{\gtype[aa]}=S_{\gtype[ab]}=
\left( \frac{2(\we_3-\we_5)}{(\we_4-\we_5)} + \frac{2(\we_3-\we_4)}{(\we_5-\we_4)}\right)^3
=2^3=8;
\]
\[
S_{\gtype[b]}=S_{\gtype[ba]}=S_{\gtype[bb]}=
\left(\frac{2(\we_3-\we_5)}{(\we_3-\we_5)} + 0 \right)^3 = 2^3 = 8;
\]
\[
S_{\gtype[c]}=S_{\gtype[ca]}=S_{\gtype[cb]}=
\left(\frac{2(\we_3-\we_4)}{(\we_3-\we_4)} + 0 \right)^3 = 2^3 = 8;
\]
\[
S_{\gtype[d]}=\left(\frac{(\we_3-\we_5)}{(\we_4-\we_5)} + \left(\frac{1}{(\we_5-\we_4)} +
\frac{1}{(\we_3-\we_4)}\right)(\we_3-\we_4) + 0 \right)^3= 2^3=8;
\]
\[
S_{\gtype[e]}=\left( \frac{(\we_3-\we_5)}{(\we_3-\we_5)} + 0 +
\frac{(\we_3-\we_4)}{(\we_3-\we_4)} \right)^3 = 2^3 = 8;
\]
\[
S_{\gtype[f]}=\left(\frac{(\we_3-\we_4)}{(\we_5-\we_4)} + \left(\frac{1}{(\we_3-\we_5)} +
\frac{1}{(\we_4-\we_5)}\right)(\we_3-\we_5) + 0 \right)^3= 2^3=8.
\]

All we are left with is computing the $T_{\gtype}$--terms and summing up everything. To
save some space and writing, we will write again $\we_{12}=\we_1-\we_2=-\we_{21}$. So
let us compute the $T_{\gtype}$--terms:
\begin{align*}
T_{\gtype[a]}&= \frac{(-1)^2 \cdot 2^4}{2^2 \cdot (\we_4-\we_5)^4} \cdot
\frac{(\we_{21}+\frac{5}{2}\we_5+\frac{1}{2}\we_4)\cdots
(\we_{21}+\frac{1}{2}\we_5+\frac{5}{2}\we_4)}{(\we_3-\we_5)(\we_3-\frac{1}{2}\we_4-\we_5)
(\we_3-\we_4)}\cdot\\
&\phantom{= }\cdot\frac{1\cdot\left(\frac{2}{(\we_4-\we_5)}\right)^{-2}}{\frac{(\we_4-\we_5)}{2}} \cdot
\frac{1\cdot\left(\frac{2}{(\we_5-\we_4)}\right)^{-2}}{\frac{(\we_5-\we_4)}{2}}\\
&= -
\frac{(\we_{21}+\frac{5}{2}\we_5+\frac{1}{2}\we_4)\cdots(\we_{21}+\frac{1}{2}\we_5+\frac{5}{2}\we_4)}{
(\we_4-\we_5)^2(\we_3-\we_5)(\we_3-\we_4)(\we_3-\frac{1}{2}\we_4-\frac{1}{2}\we_5)}.
\end{align*}
\begin{align*}
T_{\gtype[aa]}&=\left[ \frac{(-1)\cdot 1^2}{1^2\cdot
(\we_4-\we_5)^2} \cdot
\frac{(\we_{21}+2\we_5+\we_4)(\we_{21}+\we_5+2\we_4)}{(\we_3-\we_5)(\we_3-\we_4)}
\right]^2 \cdot\\
&\phantom{= }\cdot\left[ \frac{1\cdot (\we_4-\we_5)^2}{(\we_4-\we_5)}\right]^2
 \cdot
\frac{(\we_{21}+3\we_4)(\we_3-\we_4)(\we_5-\we_4)\cdot
\left(\frac{2}{(\we_5-\we_4)}\right)^{-1}}{(\we_5-\we_4)^2}\\
&=\frac{1}{2}
\frac{(\we_{21}+3\we_4)(\we_{21}+2\we_5+\we_4)^2(\we_{21}+\we_5+2\we_4)^2}{(\we_4-\we_5)^2(\we_3-\we_5)^2(\we_3-\we_4)}.
\end{align*}
\begin{align*}
T_{\gtype[ab]}&=\left[ \frac{(-1)\cdot 1^2}{1^2\cdot
(\we_4-\we_5)^2} \cdot
\frac{(\we_{21}+2\we_5+\we_4)(\we_{21}+\we_5+2\we_4)}{(\we_3-\we_5)(\we_3-\we_4)}
\right]^2 \cdot\\
&\phantom{=} \cdot\left[ \frac{1\cdot (\we_5-\we_4)^2}{(\we_5-\we_4)}\right]^2
\cdot
\frac{(\we_{21}+3\we_5)(\we_3-\we_5)(\we_4-\we_5)\cdot
\left(\frac{2}{(\we_4-\we_5)}\right)^{-1}}{(\we_4-\we_5)^2}\\
&=\frac{1}{2}
\frac{(\we_{21}+3\we_5)(\we_{21}+2\we_5+\we_4)^2(\we_{21}+\we_5+2\we_4)^2}{(\we_4-\we_5)^2(\we_3-\we_4)^2(\we_3-\we_5)}.
\end{align*}
\begin{align*}
T_{\gtype[b]}&= \frac{(-1)^2 \cdot 2^4}{2^2 \cdot (\we_3-\we_5)^4} \cdot
\frac{(\we_{21}+\frac{5}{2}\we_5+\frac{1}{2}\we_3)\cdots
(\we_{21}+\frac{1}{2}\we_5+\frac{5}{2}\we_3)}{(\we_4-\we_5)(\we_4-\frac{1}{2}\we_3-\we_5)
(\we_4-\we_3)}\cdot\\
&\phantom{= }\cdot\frac{1\cdot\left(\frac{2}{(\we_3-\we_5)}\right)^{-2}}{\frac{(\we_3-\we_5)}{2}} \cdot
\frac{1\cdot\left(\frac{2}{(\we_5-\we_3)}\right)^{-2}}{\frac{(\we_5-\we_3)}{2}}\\
&= -
\frac{(\we_{21}+\frac{5}{2}\we_5+\frac{1}{2}\we_3)\cdots(\we_{21}+\frac{1}{2}\we_5+\frac{5}{2}\we_3)}{
(\we_3-\we_5)^2(\we_4-\we_5)(\we_4-\we_3)(\we_4-\frac{1}{2}\we_3-\frac{1}{2}\we_5)}.
\end{align*}
\begin{align*}
T_{\gtype[ba]}&=\left[ \frac{(-1)\cdot 1^2}{1^2\cdot
(\we_5-\we_3)^2} \cdot
\frac{(\we_{21}+2\we_3+\we_5)(\we_{21}+\we_3+2\we_5)}{(\we_4-\we_3)(\we_4-\we_5)}
\right]^2 \cdot\\
&\phantom{=} \cdot\left[ \frac{1\cdot (\we_3-\we_5)^2}{(\we_3-\we_5)}\right]^2
\cdot
\frac{(\we_{21}+3\we_3)(\we_4-\we_3)(\we_5-\we_3)\cdot
\left(\frac{2}{(\we_5-\we_3)}\right)^{-1}}{(\we_5-\we_3)^2}\\
&=\frac{1}{2}
\frac{(\we_{21}+3\we_3)(\we_{21}+2\we_3+\we_5)^2(\we_{21}+\we_3+2\we_5)^2}{(\we_5-\we_3)^2(\we_4-\we_5)^2(\we_4-\we_3)}.
\end{align*}
\begin{align*}
T_{\gtype[bb]}&=\left[ \frac{(-1)\cdot 1^2}{1^2\cdot
(\we_5-\we_3)^2} \cdot
\frac{(\we_{21}+2\we_3+\we_5)(\we_{21}+\we_3+2\we_5)}{(\we_4-\we_3)(\we_4-\we_5)}
\right]^2 \cdot\\
&\phantom{=} \cdot\left[ \frac{1\cdot (\we_5-\we_3)^2}{(\we_5-\we_3)}\right]^2
\cdot
\frac{(\we_{21}+3\we_5)(\we_3-\we_5)(\we_4-\we_5)\cdot
\left(\frac{2}{(\we_3-\we_5)}\right)^{-1}}{(\we_3-\we_5)^2}\\
&=\frac{1}{2}
\frac{(\we_{21}+3\we_5)(\we_{21}+2\we_3+\we_5)^2(\we_{21}+\we_3+2\we_5)^2}{(\we_5-\we_3)^2(\we_4-\we_3)^2(\we_4-\we_5)}.
\end{align*}
\begin{align*}
T_{\gtype[c]}&= \frac{(-1)^2 \cdot 2^4}{2^2 \cdot (\we_3-\we_4)^4} \cdot
\frac{(\we_{21}+\frac{5}{2}\we_4+\frac{1}{2}\we_3)\cdots
(\we_{21}+\frac{1}{2}\we_4+\frac{5}{2}\we_3)}{(\we_5-\we_4)(\we_5-\frac{1}{2}\we_4-\we_3)
(\we_5-\we_3)}\cdot\\
&\phantom{=}\frac{1\cdot\left(\frac{2}{(\we_3-\we_4)}\right)^{-2}}{\frac{(\we_3-\we_4)}{2}} \cdot
\frac{1\cdot\left(\frac{2}{(\we_4-\we_3)}\right)^{-2}}{\frac{(\we_4-\we_3)}{2}}\\
&= -
\frac{(\we_{21}+\frac{5}{2}\we_4+\frac{1}{2}\we_3)\cdots(\we_{21}+\frac{1}{2}\we_4+\frac{5}{2}\we_3)}{
(\we_3-\we_4)^2(\we_5-\we_4)(\we_5-\we_3)(\we_5-\frac{1}{2}\we_4-\frac{1}{2}\we_3)}.
\end{align*}
\begin{align*}
T_{\gtype[ca]}&=\left[ \frac{(-1)\cdot 1^2}{1^2\cdot
(\we_3-\we_4)^2} \cdot
\frac{(\we_{21}+2\we_4+\we_3)(\we_{21}+\we_4+2\we_3)}{(\we_5-\we_4)(\we_5-\we_3)}
\right]^2 \cdot\\
&\phantom{=} \cdot
\left[ \frac{1\cdot (\we_3-\we_4)^2}{(\we_3-\we_4)}\right]^2
\cdot
\frac{(\we_{21}+3\we_3)(\we_4-\we_3)(\we_5-\we_3)\cdot
\left(\frac{2}{(\we_4-\we_3)}\right)^{-1}}{(\we_4-\we_3)^2}\\
&=\frac{1}{2}
\frac{(\we_{21}+3\we_3)(\we_{21}+2\we_4+\we_3)^2(\we_{21}+\we_4+2\we_3)^2}{(\we_3-\we_4)^2(\we_5-\we_4)^2(\we_5-\we_3)}.
\end{align*}
\begin{align*}
T_{\gtype[cb]}&=\left[ \frac{(-1)\cdot 1^2}{1^2\cdot
(\we_3-\we_4)^2} \cdot
\frac{(\we_{21}+2\we_4+\we_3)(\we_{21}+\we_4+2\we_3)}{(\we_5-\we_4)(\we_5-\we_3)}
\right]^2 \cdot\\
&\phantom{=} \cdot
\left[ \frac{1\cdot (\we_4-\we_3)^2}{(\we_4-\we_3)}\right]^2
\cdot
\frac{(\we_{21}+3\we_4)(\we_3-\we_4)(\we_5-\we_4)\cdot
\left(\frac{2}{(\we_3-\we_4)}\right)^{-1}}{(\we_3-\we_4)^2}\\
&=\frac{1}{2}
\frac{(\we_{21}+3\we_4)(\we_{21}+2\we_4+\we_3)^2(\we_{21}+\we_4+2\we_3)^2}{(\we_3-\we_4)^2(\we_5-\we_3)^2(\we_5-\we_4)}.
\end{align*}
\begin{align*}
T_{\gtype[d]}&=\left[ \frac{(-1)\cdot 1^2}{1^2\cdot
(\we_4-\we_5)^2} \cdot
\frac{(\we_{21}+2\we_5+\we_4)(\we_{21}+\we_5+2\we_4)}{(\we_3-\we_5)(\we_3-\we_4)}
\right]\cdot\\
&\phantom{=}\cdot \left[ \frac{(-1)\cdot 1^2}{1^2\cdot
(\we_3-\we_4)^2} \cdot
\frac{(\we_{21}+2\we_4+\we_3)(\we_{21}+\we_4+2\we_3)}{(\we_5-\we_4)(\we_5-\we_3)}
\right] \cdot \\
&\phantom{=}\cdot
\frac{(\we_{21}+3\we_4)(\we_3-\we_4)(\we_5-\we_4)\cdot
\left(\frac{1}{(\we_3-\we_4)}+\frac{1}{(\we_5-\we_4)}\right)^{-1}}{(\we_3-\we_4)(\we_5-\we_4)}
\cdot \\
&\phantom{=} \cdot \frac{1\cdot
\left(\frac{1}{(\we_4-\we_5)}\right)^{-2}}{(\we_4-\we_5)} \cdot
\frac{1\cdot \left(
\frac{1}{(\we_4-\we_3)}\right)^{-1}}{(\we_4-\we_3)}\\
&=-
\frac{(\we_{21}+2\we_5+\we_4)(\we_{21}+\we_5+2\we_4)(\we_{21}+2\we_4+\we_3)
}{(\we_3-\we_5)^2
(\we_3-\we_4)(\we_5-\we_4)(\we_5+\we_3-2\we_4)}\cdot\\
&\mbox{\hspace{15em}} \cdot(\we_{21}+\we_4+2\we_3)(\we_{21}+3\we_4).
\end{align*}
\begin{align*}
T_{\gtype[e]}&=\left[ \frac{(-1)\cdot 1^2}{1^2\cdot
(\we_3-\we_5)^2} \cdot
\frac{(\we_{21}+2\we_5+\we_3)(\we_{21}+\we_5+2\we_3)}{(\we_4-\we_5)(\we_4-\we_3)}
\right]\cdot\\
&\phantom{=}\cdot \left[ \frac{(-1)\cdot 1^2}{1^2\cdot
(\we_4-\we_3)^2} \cdot
\frac{(\we_{21}+2\we_3+\we_4)(\we_{21}+\we_3+2\we_4)}{(\we_5-\we_3)(\we_5-\we_4)}
\right] \cdot \\
&\phantom{=}\cdot
\frac{(\we_{21}+3\we_3)(\we_4-\we_3)(\we_5-\we_3)\cdot
\left(\frac{1}{(\we_4-\we_3)}+\frac{1}{(\we_5-\we_3)}\right)^{-1}}{(\we_4-\we_3)(\we_5-\we_3)}
\cdot \\
&\phantom{=}\cdot
\frac{1\cdot
\left(\frac{1}{(\we_3-\we_5)}\right)^{-2}}{(\we_3-\we_5)} \cdot
\frac{1\cdot \left(
\frac{1}{(\we_3-\we_4)}\right)^{-1}}{(\we_3-\we_4)}\\
&=-
\frac{(\we_{21}+2\we_5+\we_3)(\we_{21}+\we_5+2\we_3)(\we_{21}+2\we_3+\we_4)}{(\we_4-\we_5)^2
(\we_4-\we_3)(\we_5-\we_3)(\we_4+\we_5-2\we_3)}\cdot\\
&\mbox{\hspace{15em}} \cdot (\we_{21}+\we_3+2\we_4)(\we_{21}+3\we_3).
\end{align*}
\begin{align*}
T_{\gtype[f]}&=\left[ \frac{(-1)\cdot 1^2}{1^2\cdot
(\we_5-\we_4)^2} \cdot
\frac{(\we_{21}+2\we_4+\we_5)(\we_{21}+\we_4+2\we_5)}{(\we_3-\we_4)(\we_3-\we_5)}
\right]\cdot\\
&\phantom{=}\cdot \left[ \frac{(-1)\cdot 1^2}{1^2\cdot
(\we_3-\we_5)^2} \cdot
\frac{(\we_{21}+2\we_5+\we_3)(\we_{21}+\we_5+2\we_3)}{(\we_4-\we_5)(\we_4-\we_3)}
\right] \cdot \\
&\phantom{=}\cdot
\frac{(\we_{21}+3\we_5)(\we_3-\we_5)(\we_4-\we_5)\cdot
\left(\frac{1}{(\we_3-\we_5)}+\frac{1}{(\we_4-\we_5)}\right)^{-1}}{(\we_3-\we_5)(\we_4-\we_5)}
\cdot \\
&\phantom{=} \cdot
\frac{1\cdot
\left(\frac{1}{(\we_5-\we_4)}\right)^{-2}}{(\we_5-\we_4)} \cdot
\frac{1\cdot \left(
\frac{1}{(\we_5-\we_3)}\right)^{-1}}{(\we_5-\we_3)}\\
&=-
\frac{(\we_{21}+2\we_4+\we_5)(\we_{21}+\we_4+2\we_5)(\we_{21}+2\we_5+\we_3)}{(\we_3-\we_4)^2
(\we_3-\we_5)(\we_4-\we_5)(\we_4+\we_3-2\we_5)} \cdot\\
&\mbox{\hspace{15em}} \cdot(\we_{21}+\we_5+2\we_3)(\we_{21}+3\we_5) .
\end{align*}

Now, we just have to add up\footnote{For this fastidious exercise we have used the Maple package.} 
all the twelve terms $S_\Gamma\cdot
T_\Gamma/|\AA_\Gamma|$ --- the result is:
\[ \Phi^{2\lambda_2}(Z_3, Z_3, Z_3) = -45.\]

\paragraph{The invariants $\Phi^{p\lambda_2}(Z_2, Z_2, Z_2)$}

In the cohomology ring, we have the relation $Z_2=Z_1-3Z_3$. Since
the Gromov--Witten invariant is linear and commutative in its arguments, we obtain
\begin{align*}
\Phi^{p\lambda_2}(Z_2, Z_2, Z_2) =& \Phi^{p\lambda_2}(Z_1,Z_1,Z_1)
- 9\Phi^{p\lambda_2}(Z_1,Z_1,Z_3) +\\
&+27\Phi^{p\lambda_2}(Z_1,Z_3,Z_3) -
27\Phi^{p\lambda_2}(Z_3,Z_3,Z_3).
\end{align*}
Now observe that the weight $\we_{\sigma(\gotv)}^{l_i}$ is zero for $l_i$
corresponding to $Z_1$ and $\gotv\in \{\gotv_4, \gotv_5,
\gotv_6\}$ in the triangle of the moment polytope where
$p\lambda_2$--graphs live. This yields
\[ \Phi^{p\lambda_2}(Z_2, Z_2, Z_2) = -27\Phi^{p\lambda_2}(Z_3, Z_3, Z_3),\]
and by the results from the previous subsections we get
\begin{align*}
\Phi^0(Z_2,Z_2,Z_2)&=Z_2^3=9Z_1Z_3^2=9\\
\Phi^{\lambda_2}(Z_2,Z_2,Z_2)&=-27\cdot 3 = -81\\
\Phi^{2\lambda_2}(Z_2,Z_2,Z_2)&=-27\cdot (-45)=1215.
\end{align*}

\subsubsection{Conclusions for the quantum cohomology ring of $\P(\O_{\P^2}(3) \oplus 1)$}

In \cite[Definition 5.1]{bat93}, Batyrev defines
the quantum cohomology ring of a symplectic toric manifold
$X_\Sigma$ with symplectic form\footnote{Remember that two--forms
on $X_\Sigma$ can be represented by piecewise linear functions on the fan $\Sigma$}
$\varphi$ by (using our notation)
\[ QH^*_\varphi(X_\Sigma, \C) := \C[Z_1, \ldots, Z_n]
\left/_{\textstyle ({\rm Lin}(\Sigma) + Q_\varphi(\Sigma)),}\right.\]
where the ideal $Q_\varphi(\Sigma)$ is generated by the monomials
$B_\varphi(P)$:
\[ B_\varphi(P) = Z_{i_1} \cdots Z_{i_k} - E_\varphi(P) Z_{j_1}^{c_1}\cdots Z_{j_l}^{c_l},\]
where $P \in {\mathfrak P}$ runs over all primitive collections of
$\Sigma$, and
\begin{enumerate}
\item for a primitive collection $P=\{\gotv_{i_1}, \ldots
\gotv_{i_k}\}$, let $\gotv_P:=\gotv_{i_1}+\cdots+\gotv_{i_k}$ be
the sum of these vectors, $\sigma_P$ be the minimal cone that
contains the vector $\gotv_P$, $\gotv_{j_1}, \ldots, \gotv_{j_l}$
be the generators of $\sigma_P$, and $c_1, \ldots, c_l$ be the
integral coefficients of $\gotv_P$ with respect to these
generators:
\[ \gotv_P = c_1\gotv_{j_1} + \cdots + c_l\gotv_{j_l};\]
\item $E_\varphi(P):=
\exp(c_1\varphi(\gotv_{j_1})+\cdots+c_l\varphi(\gotv_{j_l}) -
\varphi(\gotv_{i_1})-\cdots-\varphi(\gotv_{i_l})).$
\end{enumerate}

Batyrev shows in \cite{bat93} that the structure constants of
this quantum cohomology ring are some intersection products on the
moduli space of holomorphic mappings $f: \C\P^1 \longrightarrow
X_\Sigma$. In fact, these intersection products are almost the
Gromov--Witten invariants (that where not properly defined at the
time): the only difference is that Batyrev does not compactify his
moduli space of maps.

It is therefore very interesting to know whether the quantum
cohomology ring obtained from Gromov--Witten invariants as
structure constants coincided with Batyrev's ring, \ie. to know
whether the boundary components of the moduli space have any
effect on the quantum product.

Givental has given in \cite{giv97}
an affirmative answer for Fano toric manifolds: there the two
rings are the same. For non--Fano manifolds, nothing was known so
far to the author's knowledge.

Let us suppose the quantum cohomology ring of $X_\Sigma$ was as
proposed by Batyrev\footnote{Instead of the $\exp(\ldots)$--terms we
will introduce formal variables. This is equivalent to Batyrev's notion of
the ring, since his formulas do not depend on the K\"ahler form $\varphi$.}:
\[
QH^*(X_\Sigma)=\C[Z_1, \ldots,
Z_5, q^{\lambda_1}, q^{\lambda_2}]\left/_{\left\langle
\begin{array}{l} Z_3-Z_5, Z_4-Z_5,
Z_1-Z_2-3Z_5,\\
Z_1Z_2-q^{\lambda_1},
Z_3Z_4Z_5-Z_2^3q^{\lambda_2}\end{array}\right\rangle.}\right.
\]
The relation
$Z_3Z_4Z_5-Z_2^3q^{\lambda_2}$ implies on the level of
Gromov--Witten invariants that
\[ \forall \alpha\in H^*(X_\Sigma),
A\in H_2(X_\Sigma) : \, \Phi^{A+\lambda_2}(Z_3, Z_4, Z_5, \alpha)
= \Phi^A(Z_2, Z_2, Z_2, \alpha).
\]
Taking $A=p\lambda_2$ to be a
multiple of the class $\lambda_2$, and $\alpha=P.D.(point)$ to be
the Poincar\'e dual of a point, this boils down to 
$\Phi^{(p+1)\lambda_2}(Z_3, Z_3, Z_3) =
\Phi^{p\lambda_2}(Z_2,Z_2,Z_2)$ for all $p\in {\mathbb N}$, which is obviously not satisfied
for $p$ as small as $p=0$ and $p=1$.

\begin{coro}
Batyrev's definition as in \cite{bat93} of the
quantum cohomology ring of a symplectic toric manifold does not
coincide with the one obtained from Gromov--Witten invariants as
defined in \cite{bf97} for the manifold $X_\Sigma=P_{{\C\P}^2}({\O}(3) \oplus 1)$.
\end{coro}

\begin{rema}
The example given above is not the simplest one can give for a manifold
where the quantum cohomology ring defined by Gromov--Witten invariants does
not coincide with Batyrev's ring. Consider the class of Fano surfaces
$F_k:=\P_{\C\P^1}(O(k)\oplus 1)$. All Fano surfaces are K\"ahler manifolds, and
for two different parameters $k_1, k_2 \in \Z$ $F_{k_1}$ and $F_{k_2}$ are in
the same deformation class of symplectic manifolds if and only if $k_1 \equiv
k_2 \mod 2$ (otherwise they are topologically different). Hence, via the
diffeomorphism identifying $F_{k_1}$ and $F_{k_2}$, the Gromov--Witten invariants
are the same, and so is the quantum cohomology ring.

Batyrev's rings, however, are not identified by this diffeomorphism.  
\end{rema}

\bibliographystyle{alpha}
\bibliography{mathbib}

\newcommand{\etalchar}[1]{$^{#1}$}
\begin{thebibliography}{Knu83b}

\bibitem[Aud91]{aud91}
Mich\`ele Audin.
\newblock {\em The topology of torus actions on symplectic manifolds}.
\newblock Number~93 in Progress in Math. Birkh\"auser Verlag, 1991.

\bibitem[Aud97]{aud96}
Mich\`ele Audin.
\newblock Cohomologie quantique, {S\'eminaire} {Bourbaki}, vol.\ 1995/1996.
\newblock {\em Ast\'erisque}, 241:29--58, 1997.

\bibitem[Bat91]{bat91}
Victor~V. Batyrev.
\newblock On the classification of smooth projective toric varieties.
\newblock {\em T\^ohoku Math.~J.}, 43:569--585, 1991.

\bibitem[Bat93]{bat93}
Victor~V. Batyrev.
\newblock Quantum cohomology rings of toric manifolds.
\newblock {\em Ast\'erisque}, 218:9--34, 1993.

\bibitem[Beh97]{beh97}
Kai Behrend.
\newblock {Gromov--Witten} invariants in algebraic geometry.
\newblock {\em Invent. Math.}, 127:601--627, 1997.

\bibitem[BF97]{bf97}
Kai Behrend and Barbara Fantechi.
\newblock The intrinsic normal cone.
\newblock {\em Invent. Math.}, 128(1):45--88, 1997.

\bibitem[BM96]{bm96}
Kai Behrend and Yuri Manin.
\newblock Stacks of stable maps and {Gromov--Witten} invariants.
\newblock {\em Duke Math. J.}, 85(1):1--60, 1996.

\bibitem[BS77]{bs77}
C.~B\v{a}nic\v{a} and O.~St\v{a}n\v{a}\c{s}il\v{a}.
\newblock {\em M\'ethodes alg\`ebriques dans la th\'eorie globales des espaces
  complexes}.
\newblock Varia Mathematica. Gauthier-- Villars, 1977.

\bibitem[Buc81]{buc81}
Ragnar~O. Buchweitz.
\newblock {\em Contributions \`a la th\'eorie des singularit\'es}.
\newblock PhD thesis, Universit\'e Paris 7, 1981.

\bibitem[Cox97]{cox97}
David~A. Cox.
\newblock Recent developments in toric geometry.
\newblock In Koll\'ar et~al. \cite{kol97b}, pages 389--436.

\bibitem[Dan78]{dan78}
V.~I. Danilov.
\newblock The geometry of toric varieties.
\newblock {\em Russian Math.\ Surveys}, 33(2):97--154, 1978.

\bibitem[Del88]{del88}
Thomas Delzant.
\newblock Hamiltoniens p\'eriodiques et images convexes de l'application
  moment.
\newblock {\em Bull.\ Soc.\ math.\ France}, 116:315--339, 1988.

\bibitem[DM69]{dm69}
P.~Deligne and D.~Mumford.
\newblock Irreducibility of the space of curves of given genus.
\newblock {\em Publ. Math. IHES}, 36:75--110, 1969.

\bibitem[FO99]{fo96}
K.~Fukaya and K.~Ono.
\newblock Arnold conjecture and {Gromov--Witten} invariant.
\newblock {\em Topology}, 38(5):933--1048, 1999.

\bibitem[FP97]{fp97}
W.~Fulton and R.~Pandharipande.
\newblock Notes on stable maps and quantum cohomology.
\newblock In Koll\'ar et~al. \cite{kol97b}, pages 45--96.

\bibitem[Ful93]{ful93}
William Fulton.
\newblock {\em Introduction to toric varieties}.
\newblock Number 131 in Annals of Mathematics Studies. Princeton Univ.~Press,
  1993.

\bibitem[Giv96]{giv96}
Alexander~B. Givental.
\newblock Equivariant {Gromov-Witten} invariants.
\newblock {\em Int.\ Math.\ Res.\ Not.}, 13:613--663, 1996.

\bibitem[Giv98]{giv97}
Alexander~B. Givental.
\newblock A mirror theorem for toric complete intersections.
\newblock In {\em Topological field theory, primitive forms and related topics,
  Kyoto, 1996}, pages 141--175. Birkh\"auser Verlag, 1998.

\bibitem[GP99]{gp97}
T.~Graber and R.~Pandharipande.
\newblock Localization of virtual classes.
\newblock {\em Invent. math.}, 135(2):487--518, 1999.

\bibitem[Gro85]{gro85}
M.~Gromov.
\newblock Pseudo--holomorphic curves in symplectic geometry.
\newblock {\em Invent. math.}, 82:307--347, 1985.

\bibitem[Har66]{har66}
Robin Hartshorne.
\newblock {\em Residues and Duality}.
\newblock Number~20 in Lect. Notes in Math. Springer Verlag, 1966.

\bibitem[HM98]{hm98}
Joe Harris and Ian Morrison.
\newblock {\em Moduli of curves}.
\newblock Number 187 in Graduate Texts in Math. Springer Verlag, 1998.

\bibitem[Ill71]{ill71}
Luc Illusie.
\newblock {\em Complexe Cotangent et D\'eformations I}.
\newblock Number 239 in Lect. Notes in Math. Springer Verlag, 1971.

\bibitem[K{\etalchar{+}}97]{kol97b}
J\'anos Koll\'ar et~al., editors.
\newblock {\em {Algebraic geometry. Proceedings of the Summer Research
  Institute, Santa Cruz, CA, USA, July 9--29, 1995}}, volume 62:2 of {\em Proc.
  Symp. Pure Math.}, 1997.

\bibitem[Kee92]{kee92}
Sean Keel.
\newblock Intersection theory of moduli space of stable $n$--pointed curves of
  genus zero.
\newblock {\em Trans.\ Amer.\ Math.\ Soc.}, 330(2):545--574, 1992.

\bibitem[KM76]{knu76a}
Finn~F. Knudsen and David Mumford.
\newblock The projectivity of the moduli space of stable curves, {I}:
  {Preliminaries} on "det" and "{Div}".
\newblock {\em Math. Scand.}, 39:19--55, 1976.

\bibitem[Knu83a]{knu83b}
Finn~F. Knudsen.
\newblock The projectivity of the moduli space of stable curves, {II}: {The}
  stacks {$M_{g,n}$}.
\newblock {\em Math. Scand.}, 52:161--199, 1983.

\bibitem[Knu83b]{knu83c}
Finn~F. Knudsen.
\newblock The projectivity of the moduli space of stable curves, {III}: {The}
  line bundles on {$M_{g,n}$} and a proof of the projectivity of {$M_{g,n}$} in
  characteristic 0.
\newblock {\em Math. Scand.}, 52:200--212, 1983.

\bibitem[Kon92]{kon92}
Maxim Kontsevich.
\newblock Intersection theory on the moduli space of curves and the matrix
  {Airy} function.
\newblock {\em Comm.\ Math.\ Phys.}, 147:1--23, 1992.

\bibitem[Kon95]{kon95}
Maxim Kontsevich.
\newblock Enumeration of rational curves via torus actions.
\newblock In R.~Dijkgraaf, C.~Faber, and G.~van~der Geer, editors, {\em The
  moduli space of curves}, pages 335--368. Birkh\"auser Verlag, 1995.

\bibitem[LLY97]{lly97}
Bong~H. Lian, Kefeng Liu, and Shing-Tung Yau.
\newblock Mirror principles {I}.
\newblock {\em Asian J. Math.}, 1(4):729--763, 1997.

\bibitem[LMB92]{lm92}
G\'erard Laumon and Laurent Moret-Bailly.
\newblock Champs alg\'ebriques.
\newblock Preprint 92--42, Universit\'e de Paris--Sud, Math\'ematiques, 91405
  Orsay, France, 1992.

\bibitem[LT98a]{lt98b}
Jun Li and Gang Tian.
\newblock Comparison of the algebraic and the symplectic {Gromov--Witten}
  invariants.
\newblock Preprint alg--geom/9712035, arXiv preprint server, January 1998.

\bibitem[LT98b]{lt98a}
Jun Li and Gang Tian.
\newblock Virtual moduli cycles and {Gromov--Witten} invariants of algebraic
  varieties.
\newblock {\em J.\ Amer.\ Math.\ Soc.}, 11(1):119--174, 1998.

\bibitem[LT98c]{lt96}
Jun Li and Gang Tian.
\newblock Virtual moduli cycles and {Gromov--Witten} invariants of general
  symplectic manifolds.
\newblock In {\em Topics in symplectic $4$--manifolds, Irvine, 1996}, pages
  47--83. International Press, 1998.

\bibitem[Oda88]{oda88}
Tadao Oda.
\newblock {\em Convex bodies and algebraic geometry}.
\newblock Number~15 in Ergebnisse der Mathematik und ihrer Grenzgebiete,
  3.~Folge. Springer Verlag, 1988.

\bibitem[QR98]{qr98}
Zhenbo Qin and Yongbin Ruan.
\newblock Quantum cohomology of projective bundles over {${\mathbb P}^n$}.
\newblock {\em Trans.\ Amer.\ Math.\ Soc.}, 350(9):3615--3638, September 1998.

\bibitem[RT95]{rt95}
Yongbin Ruan and Gang Tian.
\newblock A mathematical theory of quantum cohomology.
\newblock {\em J.\ Diff.\ Geom.}, 42(2):259--367, 1995.

\bibitem[Rua96]{rua96}
Yongbin Ruan.
\newblock Virtual neighborhoods and pseudo--holomorphic curves.
\newblock Preprint alg--geom/9611021, arXiv preprint server, November 1996.

\bibitem[Sie96]{sie96}
Bernd Siebert.
\newblock {Gromov--Witten} invariants for general symplectic manifolds.
\newblock Preprint dg-ga/9608005, arXiv preprint server, August 1996.

\bibitem[Sie97]{sie97}
Bernd Siebert.
\newblock An update on (small) quantum cohomology.
\newblock Preprint, Ruhr--Universit\"at Bochum, March 1997.
\newblock To appear in: Proceedings of the conference on Geometry and Physics,
  Montreal, 1995.

\bibitem[Sie98]{sie98}
Bernd Siebert.
\newblock Algebraic and symplectic {Gromov--Witten} invariants coincide.
\newblock Preprint math/9804018, arXiv preprint server, April 1998.

\bibitem[Spi99a]{spi99}
Holger Spielberg.
\newblock A formula for the {Gromov--Witten} invariants of toric varieties.
\newblock PhD thesis/Preprint 1999/11, IRMA, Universit\'e Louis Pasteur,
  Strasbourg, February 1999.

\bibitem[Spi99b]{spi99b}
Holger Spielberg.
\newblock The {Gromov--Witten} invariants of symplectic toric manifolds, and
  their quantum cohomology ring.
\newblock {\em {C.~R.~Acad.~Sci.~Paris}, S\'erie I}, 329(8):699--704, 1999.

\bibitem[Vis89]{vis89}
Angelo Vistoli.
\newblock Intersection theory on algebraic stacks and on their moduli spaces.
\newblock {\em Invent. math.}, 97:613--670, 1989.

\bibitem[Wit91]{wit91}
Edward Witten.
\newblock Two--dimensional gravity and intersection theory on moduli spaces.
\newblock {\em Surveys Diff.\ Geom.}, 1:243--310, 1991.

\end{thebibliography}

\end{document}